\documentclass[12pt]{article}
\usepackage{latexsym, amssymb, enumerate, amsmath}
\usepackage[dvipdfmx]{graphicx}
\begin{document}
\newcommand{\vare}{\varepsilon} 
\newtheorem{tth}{Theorem}[section]
\newtheorem{dfn}[tth]{Definition}
\newtheorem{lem}[tth]{Lemma}
\newtheorem{prop}[tth]{Proposition}
\newtheorem{coro}[tth]{Corollary}
\renewcommand{\quad}{\hspace*{1em}}
\begin{center}
{\Large {\bf Multi-bifurcations of 
Wavefronts on $r$-corners } }
\vspace*{0.4cm}\\
{\large Takaharu Tsukada}
\footnote{Higashijujo 3-1-16 
Kita-ku, Tokyo 114-0001
JAPAN. e-mail : tsukada@math.chs.nihon-u.ac.jp}
\vspace*{0.2cm}\\
{\large  College of Humanities \& Sciences, Department of Mathematics,\\
 Nihon University}\end{center}
\begin{abstract}
We extend the notion of reticular Legendrian unfoldings in order to
investigate multi-time bifurcations of wavefronts generated by an $r$-corner.
We give a classification list of generic and stable bifurcations with two time parameter
and give all generic figures in the plane and the space.
\end{abstract}

\section{Introduction}
\quad
Legendrian singularity can be found in many problems of differential geometry, 
calculus of variations and mathematical physics. 
One of the most successful their applications is 
the study of singularity of wavefronts. 
Bifurcation of wavefronts generated by a hypersurface without boundary  in a smooth manifold
is investigated as the theory of Legendrian unfoldings by S.Izumiya \cite{izumiya1}.
We investigated the theory of reticular Legendrian unfoldings 
in order to describe bifurcations of wavefronts generated by  a hypersurface with an $r$-corner in  \cite{bifsemi}.
These theories are investigated on one-parameter bifurcations of wavefronts.
In this paper we investigate $m$-parameter bifurcations of wavefronts on an $r$-corner. 
Since almost theory can be proved by the parallel methods of \cite{bifsemi},
we give our theory along the paper and omit the details the parts which can be prove by the parallel methods.\\

Let us consider a $m$-parameter family 
$\{ L_{\sigma,t} \}_{\sigma\subset I_r,t\in ({\mathbb R}^m,0)}$ of contact regular $r$-cubic configurations on 
$J^1({\mathbb R}^n,{\mathbb R})$ defined by contact embedding germs $C_t:
(J^1({\mathbb R}^n,{\mathbb R}),0)
\rightarrow J^1({\mathbb R}^n,{\mathbb R})$ depending 
smoothly on $t\in ({\mathbb R}^m,0)$ such that $C_0(0)=0$, $L_{\sigma,t} =
C_t(L^0_\sigma) $
for all $\sigma\subset I_r,\ t\in ({\mathbb R}^m,0)$.
We investigate bifurcations of wavefronts of $\{ L_{\sigma,t} \}_{\sigma\subset I_r}$ around time $0$.
In order to realize this, we shall need to extend the notion of {\em reticular Legendrian unfoldings} which is 
defined in \cite{bifsemi}.

\section{Stabilities of unfoldings}\label{unfold:sec} 
\quad
In this section we recall  the theory of function germs  with respect to 
{\it the reticular $t$-${\cal P}$-${\cal K}$-equivalence relation}
which is developed  in \cite{tPKfunct}.

Let ${\mathbb H}^r=\{ (x_1,\ldots,x_r)\in {\mathbb R}^r|x_1\geq 0,\ldots,x_r\geq 0\}$ 
be an $r$-corner.
We denote by ${\cal E}(r;k_1,r;k_2)$ the set of all germs at $0$ of
smooth maps ${\mathbb H}^r\times {\mathbb R}^{k_1} \rightarrow 
{\mathbb H}^r\times {\mathbb R}^{k_2}$ and set ${\mathfrak M}(r;k_1,r;k_2)=
\{ f\in {\cal
E}(r;k_1,r;k_2)|f(0)=0 \}$.
We denote ${\cal E}(r;k_1,k_2)$ for ${\cal E}(r;k_1,0;k_2)$ and 
denote ${\mathfrak M}(r;k_1,k_2)$ for ${\mathfrak M}(r;k_1,0;k_2)$.

 If $k_2=1$ we write simply ${\cal E}(r;k)$ for 
${\cal E}(r;k,1)$
and ${\mathfrak M}(r;k)$ for ${\mathfrak M}(r;k,1)$. 
Then ${\cal E}(r;k)$ is an ${\mathbb R}$-algebra in the usual
way and ${\mathfrak M}(r;k)$ is its unique maximal ideal. 
We also denote by ${\cal E}(k)$ for 
${\cal E}(0;k)$
and ${\mathfrak M}(k)$ for ${\mathfrak M}(0;k)$.
We remark that ${\cal E}(r;k,p)$ is 
an ${\cal E}(r;k)$-module generated by $p$-elements.

We
denote by $J^l(r+k,p)$ the set of $l$-jets at $0$ of germs in ${\cal
E}(r;k,p)$. There are natural projections:
\[ \pi_l:{\cal E}(r;k,p)\longrightarrow J^l(r+k,p),\ 
\pi^{l_1}_{l_2}:J^{l_1}(r+k,p)\longrightarrow J^{l_2}(r+k,p)\ (l_1 > l_2).  \]
We write $j^lf(0)$ for $\pi_l(f)$ for each $f\in {\cal E }(r;k,p)$.

Let $(x,y)=(x_1,\cdots,x_r,y_1,\cdots,y_k)$ be a fixed
coordinate system of $({\mathbb H}^r\times {\mathbb R}^k,0)$. 
We denote by 
${\cal B}(r;k)$ the group of
diffeomorphism germs $({\mathbb H}^r\times {\mathbb R}^{k},0)\rightarrow 
({\mathbb H}^r\times {\mathbb R}^{k},0)$ of the form:
\[ \phi(x,y)=(x_1\phi_1^1(x,y),\cdots,x_r\phi_1^r(x,y),\phi_2^1(x,y),\cdots,\phi_2^k(x,y)
). \]

We denote by 
${\cal B}_n(r;k+n)$ the group of
diffeomorphism germs $({\mathbb H}^r\times {\mathbb R}^{k+n},0)\rightarrow 
({\mathbb H}^r\times {\mathbb R}^{k+n},0)$ of the form:
\[
\phi(x,y,u)=(x_1\phi_1^1(x,y,u),\cdots,x_r\phi_1^r(x,y,u),\phi_2^1(x,y,u),\cdots,\phi_2^k(x,y,u)
,\phi_3^1(u),\ldots,\phi_3^n(u)).\]
We denote $\phi(x,y,u)=(x\phi_1(x,y,u),\phi_2(x,y,u),\phi_3(u))$, 
$\frac{\partial f_0}{\partial y}=(\frac{\partial f_0}{\partial y_1},$ $\cdots,\frac{\partial f_0}{\partial y_k})$, and denote other notations analogously.
\begin{lem}\label{gw1.8:cor}{\rm (cf., \cite[Corollary 1.8]{spsing}), also see \cite[Lemma2.1]{tPKfunct})}
Let $B$ be a submodule of ${\cal E}(r;k+n+m)$,
$A_1$ be a finitely generated 
${\cal E}(m)$-submodule of ${\cal E}(r;k+n+m)$ generated $d$-elements, and 
$A_2$ be a finitely generated ${\cal E}(n+m)$ submodule of ${\cal E}(r;k+n+m)$.
Suppose 
\[ {\cal E}(r;k+n+m)=B+A_2+A_1+{\mathfrak M}(m){\cal E}(r;k+n+m)
+{\mathfrak M}(n+m)^{d+1}{\cal E}(r;k+n+m). \]
Then 
\[{\cal E}(r;k+n+m)=B+A_2+A_1,\]
\[ {\mathfrak M}(n+m)^d{\cal E}(r;k+n+m)\subset B+A_2+{\mathfrak M}(m){\cal E}(r;k+n+m).\]
\end{lem}

We recall the stabilities of $n$-dimensional unfolding under 
{\it reticular ${\cal P}$-${\cal K}$-equivalence} which is developed in 
\cite{retLeg}.

We say that $f_0,g_0\in{\cal E}(r;k)$ are {\it reticular ${\cal K}$-equivalent} if
there exist $\phi\in{\cal B}(r;k)$ and a unit $a\in {\cal E}(r;k)$
such that $g_0=a\cdot f_0\circ \phi$. 

We say that a function germ $f_0\in {\mathfrak M}(r;k)$ is {\it reticular ${\cal K}$-$l$-determined} 
if all function germ which has same $l$-jet of $f_0$ is 
reticular ${\cal K}$-equivalent to $f_0$.
If $f_0$ is reticular ${\cal K}$-$l$-determined for some $l$, then we say that
$f_0$ is reticular ${\cal K}$-finitely determined.\\
\begin{lem}\label{findetc:lm}{\rm (see \cite[Lemma 2.3]{tPKfunct})}
Let $f_0(x,y)\in {\mathfrak M}(r;k)$ and let 
\[ {\mathfrak M}(r;k)^{l+1}\subset {\mathfrak M}(r;k)(\langle f_0,
x\frac{\partial f_0}{\partial
x}\rangle +{\mathfrak M}(r;k)\langle
\frac{\partial f_0}{\partial y}\rangle )
+{\mathfrak M}(r;k)^{l+2},\]
then $f_0$ is reticular ${\cal K}$-$l$-determined. 
Conversely if $f_0(x,y)\in {\mathfrak M}(r;k)$ is reticular ${\cal K}$-$l$-determined, then 
\[ {\mathfrak M}(r;k)^{l+1}\subset \langle f_0,x\frac{\partial f_0}{\partial x}
\rangle_{ {\cal E}(r;k) }  +{\mathfrak M}(r;k)\langle
\frac{\partial f_0}{\partial y}\rangle. \]
\end{lem}   

We say that $f,g\in{\cal E}(r;k+n)$ are {\it reticular ${\cal P}$-${\cal K}$-equivalent} if
there exist $\Phi\in{\cal B}_n(r;k+n)$ and a unit $\alpha\in {\cal E}(r;k+n)$
such that $g=\alpha\cdot f\circ \Phi$.\\

We say that $f(x,y,u)\in{\mathfrak M}(r;k+n)$ is {\it reticular ${\cal P}$-${\cal K}$-infinitesimally stable} if
\[
{\cal E}(r;k+n)= \langle f,x\frac{\partial
 f}{\partial x},\frac{\partial f}{\partial y}\rangle_{ {\cal
 E}(r;k+n) }+\langle \frac{\partial f}{\partial u}\rangle_{{\cal E}(n)}. 
\]

We define the several stabilities of unfolding of function germ under the reticular ${\cal P}$-${\cal K}$-equivalence
in ${\mathfrak M}(r;k)$ in \cite{tPKfunct}.
We have the following theorem:
\begin{tth}\label{pk:tth}{\rm (see \cite[Theorem 2.5]{tPKfunct})}
Let $f\in {\mathfrak M}(r;k+n)$ be an unfolding of $f_0\in {\mathfrak M}(r;k)$. 
Then the following are equivalent. \\
{\rm (1)} $f$ is reticular ${\cal P}$-${\cal K}$-stable.\\
{\rm (2)} $f$ is reticular ${\cal P}$-${\cal K}$-versal.\\
{\rm (3)} $f$ is reticular ${\cal P}$-${\cal K}$-infinitesimally versal. \\
{\rm (4)} $f$ is reticular ${\cal P}$-${\cal K}$-infinitesimally stable. \\
{\rm (5)} $f$ is reticular ${\cal P}$-${\cal K}$-homotopically stable. 
\end{tth}

We say that $F,G\in{\cal E}(r;k+n+m)$ are {\it reticular $t$-${\cal P}$-${\cal K}$-equivalent} if
there exist $\Phi\in{\cal B}(r;k+n+m)$ and a unit $\alpha\in {\cal E}(r;k+n+m)$
such that \\
(1) $\Phi$ can be written in the form: $\Phi(x,y,u,t)=(x\phi_1(x,y,u,t),\phi_2(x,y,u,t),\phi_3(u,t),\phi_4(t))$, \\
(2) $G=\alpha\cdot F\circ \Phi$.\\

We say that $F(x,y,u,t)\in {\mathfrak M}(r;k+n+m)$ is {\it reticular $t$-${\cal P}$-${\cal K}$-infinitesimally stable} if
\begin{equation}
{\cal E}(r;k+n+m) =  \langle F,x\frac{\partial F}{\partial x},
\frac{\partial F}{\partial y} \rangle_{ {\cal E}(r;k+n+m) }+
\langle \frac{\partial F}{\partial u}\rangle_{{\cal E}(n+m)}+\langle 
\frac{\partial F}{\partial t}\rangle_{{\cal E}(m)}.\label{inftsa}
\end{equation}

We define the several stabilities of unfolding of function germ under the reticular $t$-${\cal P}$-${\cal K}$-equivalence
in ${\mathfrak M}(r;k+n+m)$ in \cite{tPKfunct}.
We have the following theorem:

\begin{tth}\label{mthft:th}{\rm (see \cite[Theorem 3.14]{tPKfunct})}
Let $F(x,y,u,t)\in {\mathfrak M}(r;k+n+m)$ be an unfolding of 
$f(x,y,u)\in {\mathfrak M}(r;k+n)$ 
and let $f$ is an unfolding of $f_0(x,y)\in {\mathfrak M}(r;k)$.  
Then 
following are equivalent. \\ 
{\rm (1)} There exists a non-negative number $l$ such that 
$f_0$ is reticular 
${\cal K}$-$l$-determined and  $F$ is reticular
$t$-${\cal P}$-${\cal K}$-$q$-transversal for $q\geq lm+l+m+1$.\\
{\rm (2)} $F$ is reticular $t$-${\cal P}$-${\cal K}$-stable.\\
{\rm (3)} $F$ is reticular $t$-${\cal P}$-${\cal K}$-versal.\\ {\rm (4)} $F$ is reticular
$t$-${\cal P}$-${\cal K}$-infinitesimally versal.\\ 
{\rm (5)} $F$ is reticular
$t$-${\cal P}$-${\cal K}$-infinitesimally stable.\\ 
{\rm (6)} $F$ is reticular
$t$-${\cal P}$-${\cal K}$-homotopically stable.
\end{tth}
This theorem is used in the proof of Theorem \ref{staLeg:th}.

\section{Reticular Legendrian unfoldings}\label{RLeg:unfo}
\quad 
We consider the $1$-jet bundle $J^1({\mathbb R}^n,{\mathbb R})$
with the canonical $1$-form $\theta$ and the canonical coordinate system 
$(q,z,p)=(q_1,\ldots,q_n,z,p_1,\ldots,p_n)$, the natural projection 
$\pi:J^1({\mathbb R}^n,{\mathbb R})\rightarrow {\mathbb R}^n\times {\mathbb R} ((q,z,p)\mapsto (q,z))$.
We also consider the big $1$-jet bundle  $J^1({\mathbb R}^m\times{\mathbb R}^n,{\mathbb R})$.
and the canonical $1$-form $\Theta$ on that space.
Let $(t,q)=(t_1,\ldots,t_m,q_1,\ldots,q_n)$ be the canonical coordinate system on 
${\mathbb R}^m\times{\mathbb R}^n$ and $(t,q,z,s,p)=
(t_1,\ldots,t_m,q_1,\ldots,q_n,z,s_1,\ldots,s_m,p_1,\ldots,p_n)$ 
be the corresponding coordinate system on 
$J^1({\mathbb R}^m\times{\mathbb R}^n,{\mathbb R})$.
Then the canonical $1$-form $\Theta$ is given by 
\[ \Theta=dz-\sum_{i=1}^np_idq_i-\sum_{i=1}^ms_idt_i. \]
There exists the natural projection 
\[ \Pi:J^1({\mathbb R}^m\times{\mathbb R}^n,{\mathbb R})\rightarrow 
{\mathbb R}^m\times{\mathbb R}^n\times {\mathbb R}\ \ (t,q,z,s,p)\mapsto (t,q,z).\]

Then we consider the following contact diffeomorphism germ $C$ on $(J^1({\mathbb R}^m\times {\mathbb R}^n,
{\mathbb R}),0)$:
\begin{lem}\label{C:lem}{\rm (cf., \cite[Lemma 3.1]{bifsemi})}
For any multi-family of contact embedding germs $C_t:
(J^1({\mathbb R}^n,{\mathbb R}),0)
\rightarrow J^1({\mathbb R}^n,$ ${\mathbb R})\ (C_0(0)=0)$ depending 
smoothly on $t\in ({\mathbb R}^m,0)$,
there exist  unique function germs $h_1,\ldots,h_m$ on 
$(J^1({\mathbb R}^m\times {\mathbb R}^n,{\mathbb R}),0)$ such that 
$h_i$ depends only on $t,q,z,s_i,p$ for each $i$ and 
the map germ $C:(J^1({\mathbb R}^m\times {\mathbb R}^n,{\mathbb R}),0)
\rightarrow (J^1({\mathbb R}^m\times {\mathbb R}^n,{\mathbb R}),0)$ defined by 
\[ C(t,q,z,s,p)=(t,q\circ C_t(q,z,p),z\circ C_t(q,z,p),
h(t,q,z,s,p),p\circ C_t(q,z,p))\]
is a contact diffeomorphism.
\end{lem}
The function germ $h_i$ is uniquely determined by
\begin{equation}
h_i(t,q,z,s,p)=\frac{\partial z_t}{\partial t}(q,z,p)-p_t(q,z,p)
\frac{\partial q_t}{\partial t}(q,z,p)+\alpha(t,q,z,p)s_i. \label{hs:eqn}
\end{equation}
We define that  
$\tilde{L}^0_\sigma=\{(t,q,z,s,p) \in J^1({\mathbb R}^m\times 
{\mathbb R}^n,{\mathbb R})|
q_\sigma=p_{I_r-\sigma}=q_{r+1}=\cdots=q_n=s=z=0,q_{I_r-\sigma}\geq 0 \}$
for $\sigma\subset I_r$
and ${\mathbb L}=\{(t,q,z,s,p) \in J^1({\mathbb R}^m\times 
{\mathbb R}^n,{\mathbb R})|
q_1p_1=\cdots=q_rp_r=q_{r+1}=\cdots=q_n=s=z=0,q_{I_r}\geq 0 \}$ be a 
representative as a germ 
of the union of $\tilde{L}^0_\sigma$ for all $\sigma\subset I_r$.
\begin{dfn}\label{C:dfn}{\rm
Let $C$ be a contact diffeomorphism germ on
$(J^1({\mathbb R}^m\times {\mathbb R}^n,{\mathbb R}),0)$.
We say that $C$ is {\em a ${\cal P}$-contact diffeomorphism}
if $C$ has the form:
\begin{equation}
C(t,q,z,s,p)=(t,q_C(t,q,z,p),z_C(t,q,z,p),
h_C(t,q,z,s,p),p_C(t,q,z,p))\label{tconta:eqn}
\end{equation}
and the function germ $h_C^i$ depends only on $t,q,z,s_i,p$ for each
$i=1,\ldots,m$.
}\end{dfn}
\begin{dfn}{\rm
We say that a map germ
${\cal L}:({\mathbb L},0)\rightarrow 
(J^1({\mathbb R}^m\times {\mathbb R}^n,{\mathbb R}),0)$ is
{\em a reticular Legendrian unfolding} if 
${\cal L}$ is the restriction of a ${\cal P}$-contact diffeomorphism.
We call $\{ {\cal L}(\tilde{L}^{0}_\sigma) \}_{\sigma\subset I_r}$ 
{\em the unfolded contact regular $r$-cubic configuration of } ${\cal L}$.
}\end{dfn}

We note that:
Let $\{ \tilde{L}_{\sigma} \}_{\sigma\subset I_r}$ be
an unfolded contact regular $r$-cubic configuration
 associated with an $m$-parameter family of contact regular $r$-cubic configurations $\{ L_{\sigma,t} \}_{\sigma\subset I_r,t\in ({\mathbb R}^m,0)}$.
Then there is the following relation between the wavefront 
$W_\sigma=\Pi(\tilde{L}_{\sigma})$ and 
the family of wavefronts $W_{\sigma,t}=\pi(L_{\sigma,t})$:
\[ W_\sigma=\bigcup_{t\in ({\mathbb R}^m,0)}
\{t\}\times W_{\sigma,t} \ \ \ \ \mbox{ for all } \sigma\subset I_r.\]

Let $K,\Psi$ be contact diffeomorphism germs on 
$(J^1({\mathbb R}^m\times {\mathbb R}^n,{\mathbb R}),0)$.
We say that $K$ is {\em a ${\cal P}$-Legendrian equivalence} if $K$ has 
the form:
\begin{equation}
K(t,q,z,s,p)= (\phi_1(t),\phi_2(t,q,z),\phi_3(t,q,z),\phi_4(t,q,z,s,p),
\phi_5(t,q,z,s,p))\label{PLequi}.
\end{equation}
We say that $\Psi$ is {\em a reticular ${\cal P}$-diffeomorphism} if $\pi_t\circ \Psi$ depends only on $t$ and $\Psi$ preserves 
$\tilde{L}^0_\sigma$ for all $\sigma\subset I_r$. \\

Let $\{ \tilde{L}^i_{\sigma} \}_{\sigma\subset I_r}(i=1,2)$
be unfolded contact regular $r$-cubic configurations on 
$(J^1({\mathbb R}^m\times {\mathbb R}^n,{\mathbb R}),0)$.
We say that they are 
{\em ${\cal P}$-Legendrian equivalent} 
if there exist a ${\cal P}$-Legendrian equivalence $K$
 such that 
$\tilde{L}^2_{\sigma}=K(\tilde{L}^1_{\sigma})$
for all $\sigma\subset I_r$.\vspace{2mm}

In order to understand the meaning of ${\cal P}$-Legendrian equivalence,
we observe the following: 
Let  $\{ \tilde{L}^i_{\sigma} \}_{\sigma\subset I_r}(i=1,2)$  be 
unfolded contact regular $r$-cubic configurations on 
$(J^1({\mathbb R}^m\times {\mathbb R}^n,{\mathbb R}),0)$
and 
$\{ L^i_{\sigma,t} \}_{\sigma\subset I_r,t\in ({\mathbb R}^m,0)}$
be the corresponding $m$-parameter families of 
contact regular $r$-cubic configurations on 
$J^1({\mathbb R}^n,{\mathbb R})$
 respectively.
We take the smooth $m$-parameter path germs $w_i:({\mathbb R}^m,0)\rightarrow 
(J^1({\mathbb R}^n,{\mathbb R}),0)$ such that 
$\{ L^i_{\sigma,t} \}_{\sigma\subset I_r}$ are defined at $w_i(t)$ 
for $i=1,2$.
Suppose that there exists a ${\cal P}$-Legendrian equivalence $K$ from 
$\{ \tilde{L}^1_{\sigma} \}_{\sigma\subset I_r}$
to $\{ \tilde{L}^2_{\sigma} \}_{\sigma\subset I_r}$ of the form (\ref{PLequi}).
We set  $W_{\sigma,t}^i$ be the wavefront of $L^i_{\sigma,t}$ 
for $\sigma\subset I_r,\ t\in ({\mathbb R}^m,0)$ and $i=1,2$.
We define the family of diffeomorphism $g_t:({\mathbb R}^n\times {\mathbb R},
\pi(w_1(t)))\rightarrow ({\mathbb R}^n\times {\mathbb R},
\pi(w_2(t)))$ by 
$g_t(q,z)=(\phi_2(t,q,z),\phi_3(t,q,z))$.
Then we have that $g_t(W^1_{\sigma,t})=W^1_{\sigma,\phi_1(t)}$ for all $\sigma\subset I_r,\ t\in ({\mathbb R}^m,0)$.
\vspace{2mm}

We also define the equivalence relation among reticular Legendrian unfoldings.
Let ${\cal L}_i:({\mathbb L},0)\rightarrow 
(J^1({\mathbb R}^m\times {\mathbb R}^n,{\mathbb R}),0), (i=1,2)$ 
be reticular Legendrian unfoldings.
We say that ${\cal L}_1$ and ${\cal L}_2$ are 
{\em ${\cal P}$-Legendrian equivalent} 
if there exist a ${\cal P}$-Legendrian equivalence $K$ 
and a reticular ${\cal P}$-diffeomorphism 
$\Psi$ 
such that 
$K\circ {\cal L}_1={\cal L}_2\circ \Psi$.

\begin{lem}\label{exleg:lm}{\rm (cf., \cite[Lemma 3.4]{bifsemi})}
Let $\{\tilde{L}_\sigma \}_{\sigma \subset I_r}$ be an 
unfolded contact regular $r$-cubic configuration on 
$(J^1({\mathbb R}^m\times {\mathbb R}^n,{\mathbb R}),0)$.
Then there exists 
a ${\cal P}$-contact diffeomorphism germ $C$ on $(J^1({\mathbb R}^m\times {\mathbb R}^n,{\mathbb R}),0)$
such that $C$ defines
$\{\tilde{L}_\sigma \}_{\sigma \subset I_r}$  and preserves the canonical 
$1$-form.
\end{lem}

We can construct generating families of reticular Legendrian unfoldings.
A function germ $F(x,y,t,q,z)\in{\mathfrak M}(r;k+m+n+1)$ is said to be 
{\em ${\cal P}$-$C$-non-degenerate} if 
$\frac{\partial F}{\partial x}(0)=\frac{\partial F}{\partial y}(0)=0$ and 
$x,t,F,\frac{\partial F}{\partial x},
\frac{\partial F}{\partial y}$
are independent on $({\mathbb H}^k\times {\mathbb R}^{k+m+n+1},0)$.\\

A ${\cal P}$-$C$-non-degenerate function germ 
$F(x,y,t,q,z)\in {\mathfrak M}(r;k+m+n+1)$
is called {\em a generating family} of 
a reticular Legendrian unfoldings ${\cal L}$ if
\begin{eqnarray*}
{\cal L}(\tilde{L}^0_{\sigma}) =
\{ (t,q,z,\frac{\partial F}{\partial t}/(-\frac{\partial F}{\partial z}),
\frac{\partial F}{\partial q}/(-\frac{\partial F}{\partial z}))\in
(J^1({\mathbb R}^m\times {\mathbb R}^n,{\mathbb R}),0)|\hspace{1cm}\\
\hspace{2cm}
x_\sigma=F=\frac{\partial F}{\partial x_{I_r-\sigma}}=
\frac{\partial F}{\partial y}=0,x_{I_r-\sigma}\geq 0\} 
\mbox{ for all }\sigma\subset I_r.
\end{eqnarray*}

By Lemma \ref{exleg:lm} we may assume that an extension of reticular Legendrian
unfolding preserves the canonical $1$-form.
\begin{lem} {\rm (cf., \cite[Lemma 3.5]{bifsemi})}
Let $C$ be a ${\cal P}$-contact diffeomorphism germ  
$(J^1({\mathbb R}^m\times {\mathbb R}^n,{\mathbb R}),0)\rightarrow
(J^1({\mathbb R}^m\times {\mathbb R}^n,{\mathbb R}),0)$
which preserves the canonical $1$-form.
If 
the map germ
\[ (T,Q,Z,S,P)\rightarrow (T,Q,Z,s_C(T,Q,Z,S,P),p_C(T,Q,Z,S,P))\]
is a  diffeomorphism,
there exists a function germ $H(T,Q,p)\in 
{\mathfrak M}(m+n+n)^2$ such that 
the canonical relation $P_C$ associated with $C$ has the form:
\begin{eqnarray}
P_C=\{ (T,Q,Z,-\frac{\partial H}{\partial T}(T,Q,p)+s,
-\frac{\partial H}{\partial Q},
T,-\frac{\partial H}{\partial p},
H-\langle \frac{\partial H}{\partial p},p\rangle +Z,
s,p)\},\label{cano:eqn}
\end{eqnarray}
and the function germ $F\in {\mathfrak M}(r;n+m+n+1)$ defined by
$F(x,y,t,q,z)=-z+H(t,x,0,y)+\langle y,q\rangle$
is a generating family of the reticular  Legendrian unfolding $C|_{{\mathbb L}}$.
\end{lem}

We have the following theorem which gives the relations 
between reticular Legendrian unfoldings and their generating families.
\begin{tth}\label{UCgf:th}{\rm (cf., \cite[Theorem 3.6]{bifsemi})}
{\rm (1)} For any reticular Legendrian unfolding ${\cal L}:({\mathbb L},0)
\rightarrow (J^1({\mathbb R}^m\times {\mathbb R}^n,{\mathbb R}),0)$, 
there exists a function germ $F(x,y,t,q,z)\in{\mathfrak M}(r;k+m+n+1)$ 
which is a generating family of ${\cal L}$.\\
{\rm (2)} For any ${\cal P}$-$C$-non-degenerate function germ 
$F(x,y,t,q,z)\in{\mathfrak M}(r;k+m+n+1)$ with
$\frac{\partial F}{\partial t}(0)=\frac{\partial F}{\partial q}(0)=0$, there exists a reticular
 Legendrian unfolding ${\cal L}:({\mathbb L},0)\rightarrow 
(J^1({\mathbb R}^m\times {\mathbb R}^n,{\mathbb R}),0)$ 
of which $F$
is a generating family.\\
{\rm (3)} Two reticular Legendrian unfolding are ${\cal P}$-Legendrian 
equivalent if and only if their generating families are stably reticular 
$t$-${\cal P}$-${\cal K}$-equivalent.
\end{tth}

\section{Stabilities of reticular Legendrian unfoldings}
\quad 
Let $U$ be an open set in $J^1({\mathbb R}^m\times {\mathbb R}^n,{\mathbb R})$.
We consider contact embedding germs 
$(J^1({\mathbb R}^m\times {\mathbb R}^n,{\mathbb R}),0)
\rightarrow J^1({\mathbb R}^m\times {\mathbb R}^n,{\mathbb R})$ 
and contact embeddings
$U\rightarrow J^1({\mathbb R}^m\times {\mathbb R}^n,{\mathbb R})$. 
Let $(T,Q,S,Z,P)$ and $(t,q,z,s,p)$ be canonical coordinates of the source
space and the target space respectively. 
We define the following notations:\\
$\imath:(J^1({\mathbb R}^m\times {\mathbb R}^n,{\mathbb R})\cap \{Z=0 \},0)\rightarrow 
(J^1({\mathbb R}^m\times {\mathbb R}^n,{\mathbb R}),0)$ be the inclusion map on the source space,
\begin{eqnarray*}
C_T(J^1({\mathbb R}^m\times {\mathbb R}^n,{\mathbb R}),0) & = &
\{ C| C \mbox{ is a ${\cal P}$-contact embedding germ}\\
& &  \hspace{1cm} (J^1({\mathbb R}^m\times {\mathbb R}^n,{\mathbb R}),0)\rightarrow J^1({\mathbb R}^m\times {\mathbb R}^n,{\mathbb R})\}, \\
C_T^\Theta (J^1({\mathbb R}^m\times {\mathbb R}^n,{\mathbb R}),0)
& = & \{ C\in C_T(J^1({{\mathbb R}^m\times \mathbb R}^n
,{\mathbb R}),0)|\ C^*\Theta=\Theta \}, \\
C_T^{Z} (J^1({\mathbb R}^m\times {\mathbb R}^n,{\mathbb R}),0)
 & = & \{ C\circ\imath\ |C \in 
C_T(J^1({\mathbb R}^m\times {\mathbb R}^n,{\mathbb R}),0) \},\\
C_T^{\Theta,Z} (J^1({\mathbb R}^m\times {\mathbb R}^n,{\mathbb R}),0)
 & = & \{ C\circ\imath\ |
C \in C_T^\Theta (J^1({\mathbb R}^m\times {\mathbb R}^n,{\mathbb R}),0) \}.
 \end{eqnarray*}
Let $V=U\cap \{Z=0\}$ and $\tilde{\imath}:V\rightarrow U$ be the 
inclusion map.
\begin{eqnarray*}
C_T(U,J^1({\mathbb R}^m\times {\mathbb R}^n,{\mathbb R})) & = & \{ 
\tilde{C}:U\rightarrow 
J^1({\mathbb R}^m\times {\mathbb R}^n,{\mathbb R})| \\
 & & \hspace{1cm} \tilde{C} \mbox{ is a contact embedding of the form (\ref{tconta:eqn})}\},\\
 C_T^\Theta (U,J^1({\mathbb R}^m\times {\mathbb R}^n,{\mathbb R}))  & = & \{ \tilde{C}\in 
C_T(U,J^1({\mathbb R}^m\times {\mathbb R}^n,{\mathbb R}))\ |
\tilde{C}^*\Theta=\Theta \},\\
 C_T^Z (V,J^1({\mathbb R}^m\times {\mathbb R}^n,{\mathbb R})) & = & \{ \tilde{C}\circ
\tilde{\imath}\ |\tilde{C}\in C_T(U,J^1({\mathbb R}^m\times {\mathbb R}^n,{\mathbb R}) )\},\\
 C_T^{\Theta,Z} (V,J^1({\mathbb R}^m\times {\mathbb R}^n,{\mathbb R})) & = & 
 \{ \tilde{C}\circ
\tilde{\imath}\ |\tilde{C}\in C_T^\Theta (U,J^1({\mathbb R}^m\times {\mathbb R}^n,{\mathbb R})) \}.
\end{eqnarray*}
\begin{dfn}{\rm
We define stabilities of reticular Legendrian unfoldings.
Let ${\cal L}$ be
a reticular Legendrian unfolding.\\
{\bf Stability}: We say that ${\cal L}$ is {\it stable}  if the 
following condition holds:
Let $C^{0}\in 
C_T(J^1({\mathbb R}^m\times {\mathbb R}^n,{\mathbb R}),0)$ 
be ${\cal P}$-contact embedding germs 
such that $C^{0}|_{{\mathbb L}}={\cal L}$ and 
$\tilde{C}^{0}\in C_T(U,J^1({\mathbb R}^m\times {\mathbb R}^n,{\mathbb R}))$ be
representatives of $C^{0}$.
Then there exist open neighborhoods $N_{\tilde{C}^{0}}$ of $\tilde{C}^{0}$ 
in $C^\infty$-topology
such that for any $\tilde{C}\in N_{\tilde{C}^{0}}$,
there exist points $x_0=(T,0,\ldots,0,P^0_{r+1},\ldots,P^0_n)\in U$ such that 
the reticular Legendrian unfolding 
${\cal L}_{x_0}$ 
and 
${\cal L}$ are 
${\cal P}_{(m)}$-Legendrian equivalent,
where the reticular Legendrian unfolding ${\cal L}_{x_0}$ is defined by
\[ x=(T,Q,Z,S,P)\mapsto \tilde{C}(x_0+x)-\tilde{C}(x_0)+(0,0,P^0_{r+1}Q_{r+1}+\cdots +P^0_nQ_n,0,0). \] 
{\bf Homotopical stability}: 
A one-parameter family
of ${\cal P}$-contact embedding germs  $\bar{C}:(J^1({\mathbb R}^m\times {\mathbb R}^n,{\mathbb R})
\times {\mathbb R},(0,0))\rightarrow
J^1({\mathbb R}^m\times {\mathbb R}^n,{\mathbb R})\ ((T,Q,Z,S,P,\tau)
\mapsto C_\tau(T,Q,Z,S,P))$
is called a {\em  ${\cal P}$-contact 
  deformation} of ${\cal L}$ if $C_0|_{{\mathbb L}}={\cal L}$.
A map germ $\bar{\Psi}:(J^1({\mathbb R}^m\times {\mathbb R}^n,{\mathbb R}) \times 
{\mathbb R},(0,0)) \rightarrow 
(J^1({\mathbb R}^m\times {\mathbb R}^n,{\mathbb R}),0)((T,Q,Z,S,P,\tau)
\mapsto \Psi_\tau(T,Q,Z,S,P))$
is called a {\em one-parameter
  deformation of reticular diffeomorphisms} if 
$\Psi_0=id_{J^1({\mathbb R}^m\times {\mathbb R}^n,{\mathbb R})}$ and 
$\Psi_t$ 
 is a ${\cal P}$-diffeomorphism for all $t$ around $0$. 
We say that ${\cal L}$ is {\it homotopically stable} if 
for any reticular ${\cal P}$-contact deformations 
$\bar{C}=\{ C_\tau\}$ of ${\cal L}$,
there exist 
 one-parameter families of ${\cal P}$-Legendrian 
  equivalences $\bar{K}=\{ K_\tau \}$ on 
$(J^1({\mathbb R}^m\times {\mathbb R}^n,{\mathbb R}),0)$
with $K_0=id$ of the form 
\begin{equation}K_\tau(t,q,z,s,p)=(\phi^1_\tau (t),\phi^2_\tau (t,q,z),\phi^3_\tau (t,q,z),
\phi^{4}_\tau (t,q,z,s,p),\phi^{5}_\tau (t,q,z,s,p))\label{Khomo:eqn}
\end{equation}
and 
one-parameter deformations of
reticular ${\cal P}$-diffeomorphisms $\bar{\Psi}=\{ \Psi_\tau \}$
such that $C_\tau=K_\tau\circ C_0\circ \Psi_\tau$ for $t$ around
$0$.\\
{\bf Infinitesimal stability}:
Let $C\in C_T(J^1({\mathbb R}^m\times {\mathbb R}^n,{\mathbb R}),0)$
be a ${\cal P}$-contact diffeomorphism germ.
We say that a vector field $v$ on $(J^1({\mathbb R}^m\times {\mathbb R}^n,{\mathbb R}),0)$ 
along $C$ 
is {\em an infinitesimal ${\cal P}$-contact transformation} of $C$ if 
there exists a ${\cal P}$-contact deformation 
$\bar{C}=\{C_\tau\}$ on $(J^1({\mathbb R}^m\times {\mathbb R}^n,{\mathbb R}),0)$ 
such that $C_0=C$  and
$\frac{dC_\tau}{d\tau}|_{\tau =0}=v$.
We say that a vector field $\xi$ on $(J^1({\mathbb R}^m\times {\mathbb R}^n,{\mathbb R}),0)$ is 
{\em an infinitesimal reticular ${\cal P}$-diffeomorphism} if there exists 
a one-parameter deformation of
reticular ${\cal P}$-diffeomorphisms  $\bar{\Psi}=\{ \Psi_\tau \}$ 
such that  $\frac{d\Psi_\tau}{d\tau }|_{\tau =0}=\xi$.
We say that a vector field $\eta$ on  $(J^1({\mathbb R}^m\times {\mathbb R}^n,{\mathbb R}),w)$ is
{\em an infinitesimal ${\cal P}$-Legendrian equivalence} if 
there exists a
 one-parameter family of ${\cal P}$-Legendrian equivalences $\bar{K}=\{K_\tau\}$ such that 
$K_0=id_{J^1({\mathbb R}^m\times {\mathbb R}^n,{\mathbb R})}$
and $\frac{dK_\tau}{d\tau}|_{\tau =0}=\eta$.
We say that ${\cal L}$ is  
 {\it infinitesimally stable} if for any
extension $C$ of ${\cal L}$ and any infinitesimal ${\cal P}$-contact 
transformation $v$ of $C$,
there exist infinitesimal reticular ${\cal P}$-diffeomorphisms $\xi$ and 
infinitesimal ${\cal P}$-Legendrian equivalences $\eta$
of 
the form
\begin{eqnarray}
\eta(t,q,z,s,p)=a_1(t)\frac{\partial}{\partial t}+
a_2(t,q,z)\frac{\partial}{\partial q}+a_3(t,q,z)\frac{\partial}{\partial z}
\hspace{2cm}\nonumber \\
+
a_4(t,q,z,s,p)\frac{\partial}{\partial s}+a_5(t,q,z,s,p)
\frac{\partial}{\partial p}  \label{etaform}
\end{eqnarray}
such  that $v=C_*\xi+\eta\circ C$.
}\end{dfn}

We may take an extension of a reticular Legendrian unfolding ${\cal L}$ by an element of $C^\Theta_T (J^1({\mathbb R}^m\times {\mathbb R}^n,{\mathbb R}),0)$ by
Lemma \ref{exleg:lm}.
Then as the remark after the definition of the stability of reticular 
Legendrian maps in \cite[p.121]{retLeg}, 
we may consider the following other definitions of stabilities of multi-reticular Legendrian unfoldings: 
(1) The definition given by replacing 
$C_T(J^1({\mathbb R}^m\times {\mathbb R}^n,{\mathbb R}),0)$ and
$C_T(U,J^1({\mathbb R}^m\times {\mathbb R}^n,{\mathbb R}))$ to
$C_T^\Theta(J^1({\mathbb R}^m\times {\mathbb R}^n,{\mathbb R}),0)$ and 
$C_T^\Theta(U,J^1({\mathbb R}^m\times {\mathbb R}^n,{\mathbb R}))$
of original definition respectively.
(2) The definition given by replacing to $C_T^Z(J^1({\mathbb R}^m\times {\mathbb R}^n,{\mathbb R}),0)$ and
$C_T^Z(V,J^1({\mathbb R}^m\times {\mathbb R}^n,{\mathbb R}))$ respectively. 
(3) The definition given by replacing to 
$C_T^{\Theta,Z}(J^1({\mathbb R}^m\times {\mathbb R}^n,{\mathbb R}),0)$ and
$C_T^{\Theta,Z}(V,J^1({\mathbb R}^m\times {\mathbb R}^n,{\mathbb R}))$ 
respectively, where $V=U\cap \{Z=0\}$.\\

Then we have the following lemma which is proved by the same method of the 
proof of \cite[Lemma 7.2]{retLeg}
\begin{lem}\label{sta:lm}{\rm (cf., \cite[Lemma 4.3]{bifsemi})}
The original definition and other three definitions of stabilities of reticular Legendrian unfoldings are all equivalent.
\end{lem}
By this lemma, we may choose an extension of a reticular Legendrian unfolding 
from among all of 
 $C_T(J^1({\mathbb R}^m\times {\mathbb R}^n,{\mathbb R}),0))$, 
$C_T^\Theta(J^1({\mathbb R}^m\times {\mathbb R}^n,{\mathbb R}),0))$,
 $C_T^Z(J^1({\mathbb R}^m\times {\mathbb R}^n,{\mathbb R}),0))$, and 
$C_T^{\Theta,Z}(J^1({\mathbb R}^m\times {\mathbb R}^n,{\mathbb R}),0))$.\\

We say that a function germ $H$ on $(J^1({\mathbb R}^m\times {\mathbb R}^n,{\mathbb R}),0)$ is {\em ${\cal P}$-fiber preserving} if 
$H$ has the form $H(t,q,z,s,p)=\sum_{i=1}^nh_j(t,q,z)p_j+h_0(t,q,z)+\sum_{i=1}^ma_i(t)s_i$.
\vspace{2mm}
\begin{lem}\label{infsta:t-Leglem}{\rm (cf., \cite[Lemma 4.4]{bifsemi})}
Let $C\in C_T(J^1({\mathbb R}^m\times {\mathbb R}^n,{\mathbb R}),0)$.
Then the following hold:
{\rm (1)} A vector field germ $v$ on $(J^1({\mathbb R}^m\times {\mathbb R}^n,{\mathbb R}),0)$ along $C$
 is an infinitesimal ${\cal P}$-contact transformation of $C$ if and only if 
there exists a function germ $f$ on 
$(J^1({\mathbb R}^m\times {\mathbb R}^n,{\mathbb R}),0)$ 
such that $f$ does not depend on $s$ and $v=X_f\circ C$.\\
{\rm (2)} A vector field germ $\eta$ on 
$(J^1({\mathbb R}^m\times {\mathbb R}^n,{\mathbb R}),0)$ is
an infinitesimal ${\cal P}$-Legendrian equivalence
 if and only if there exists a ${\cal P}$-fiber preserving function germ
$H$ on $(J^1({\mathbb R}^m\times {\mathbb R}^n,{\mathbb R}),0)$ 
such that $\eta=X_H$.\\
{\rm (3)} A vector field $\xi$ on 
$(J^1({\mathbb R}^m\times {\mathbb R}^n,{\mathbb R}),0)$ is
an infinitesimal reticular ${\cal P}$-diffeomorphism  if and only if 
there exists a function germ $g\in B$
such that $\xi=X_g$, where $B=\langle q_1p_1,\ldots,q_rp_r,$ $
q_{r+1},\ldots,q_n,z\rangle_{{\cal E}_{t,q,z,p}}
+\langle s \rangle_{{\cal E}_t}$.
\end{lem}

We define the several stabilities of reticular Legendrian unfoldings and we have the following theorem:
\begin{tth}\label{staLeg:th}{\rm (cf., \cite[Theorem 4.6]{bifsemi})}
Let ${\cal L}$ 
be a reticular Legendrian unfolding with a 
generating family $F(x,y,t,q,z)$.
Then the following are all equivalent.\\
{\rm (u)} $F$ is a reticular $t$-${\cal P}$-${\cal K}$-stable unfolding of 
$F|_{t=0}$.\\
{\rm (hs)} ${\cal L}$ is homotopically stable.\\
{\rm (is)} ${\cal L}$ is infinitesimally stable.\\ 
{\rm (a)} 
${\cal E}_{t,q,p}=
B_0+
\langle 1,p_1\circ C',\ldots,p_n\circ C'\rangle_{(\Pi\circ C')^*{\cal E}_{t,q,z}}+
\langle s\circ C'\rangle_{{\cal E}_t}$,
where $C'=C|_{z=s=0}$ and $B_0=\langle q_1p_1,\ldots,q_rp_r,
q_{r+1},\ldots,q_n\rangle_{{\cal E}_{t,q,p}}$.
\end{tth}
\section{Genericity of reticular Legendrian unfoldings}
\quad 
In order to give a generic classification of reticular Legendrian 
unfoldings,
we reduce our investigation to  finite dimensional jet spaces 
of ${\cal P}$-contact diffeomorphism germs.

\begin{dfn}\label{l,l+1detLeg}{\rm
Let ${\cal L}$ be a reticular Legendrian unfolding.
We say that ${\cal L}$ is $l$-determined if the following condition holds:
For any extension $C\in C_T(J^1({\mathbb R}^m\times{\mathbb R}^n,{\mathbb R}),0)$ of ${\cal L}$, 
the reticular Legendrian unfolding 
$C'|_{{\mathbb L}}$ 
and ${\cal L}$ are ${\cal P}$-Legendrian 
equivalent for all $C'\in C_T(J^1({\mathbb R}^m\times{\mathbb R}^n,{\mathbb R}),0)$ satisfying that  
$j^lC(0)=j^lC'(0)$ .
}\end{dfn}

As Lemma \ref{sta:lm}, we may consider the following other definition
of finitely determinacy of reticular Legendrian maps:\\
(1) The definition given by replacing $C_T(J^1({\mathbb R}^m\times{\mathbb R}^n,{\mathbb R}),0)$
to $C^\Theta_T(J^1({\mathbb R}^m\times{\mathbb R}^n,{\mathbb R}),0)$.\\
(2) The definition given by replacing $C_T(J^1({\mathbb R}^m\times{\mathbb R}^n,{\mathbb R}),0)$
to $C^Z_T(J^1({\mathbb R}^m\times{\mathbb R}^n,{\mathbb R}),0)$.\\
(3) The definition given by replacing $C_T(J^1({\mathbb R}^m\times{\mathbb R}^n,{\mathbb R}),0)$
to $C^{\Theta,Z}_T(J^1({\mathbb R}^m\times{\mathbb R}^n,{\mathbb R}),0)$.\\
Then the following holds by \cite[p.341 Proposition 5.6]{generic}:
\begin{prop}{\rm (cf., \cite[Proposition 5.2]{bifsemi})}
Let ${\cal L}$ be a reticular Legendrian unfolding.
Then \\
{\rm (A)} If ${\cal L}$ is $l$-determined of the original definition, then ${\cal L}$ is $l$-determined of the 
definition {\rm (1)}.\\
{\rm (B)} If ${\cal L}$ is $l$-determined of the definition {\rm (1)}, then ${\cal L}$ is $l$-determined of the definition {\rm (3)}.\\
{\rm (C)} If ${\cal L}$ is $(l+1)$-determined of the definition {\rm (3)}, then ${\cal L}$ is $l$-determined of the definition {\rm (2)}.\\
{\rm (D)} If ${\cal L}$ is $l$-determined of the definition {\rm (2)}, then 
${\cal L}$ is $l$-determined of the original definition.
\end{prop}

\begin{tth}\label{n+1det:Leg}{\rm (cf., \cite[Lemma 5.3]{bifsemi})}
Let ${\cal L}:({\mathbb L},0) \rightarrow 
(J^1({\mathbb R}^m\times{\mathbb R}^n,{\mathbb R}),0)$ be a 
reticular Legendrian unfolding.
If ${\cal L}$ is infinitesimally stable then ${\cal L}$ is $(n+m+3)$-determined.
\end{tth}
{\em Proof.}
It is enough to prove ${\cal L}$ is $(n+m+2)$-determined of Definition 
\ref{l,l+1detLeg} (3).
Let $C\in C^{\Theta,Z}_T(J^1({\mathbb R}^m\times{\mathbb R}^n,{\mathbb R}),0)$ be an 
extension of ${\cal L}$.
We may assume that $P_C$ has the form 
\[P_C=\{(T,Q,0,-\frac{\partial H}{\partial T}(T,Q,p)+s,
-\frac{\partial H}{\partial Q},T,-\frac{\partial H}{\partial p},
H-\langle \frac{\partial H}{\partial p},p\rangle ,s,p)\} \]
for some function germ $H(T,Q,p)\in {\mathfrak M}(2n+m)^2$.
Then the function germ $F(x,y,t,q,z)=-z+H_0(x,y,t)+\langle y,q\rangle  \in {\mathfrak M}(r;n+m+n+1)$
is a generating family of ${\cal L}$,
where $H_0(x,y,t)=H(t,x,0,y)\in {\mathfrak M}(r;n+m)^2$.
We have that 
 $F$ is a
reticular $t$-${\cal P}$-${\cal K}$-stable unfolding of 
$f(x,y,q,z):=-z+H_0(x,y,0)+\langle y,q\rangle \in {\mathfrak M}(r;n+n+1)$.
This means that 
\[{\cal E}(r;n+1+n+m) =\langle F,x\frac{\partial
 F}{\partial x},\frac{\partial
 F}{\partial y}\rangle_{ {\cal
 E}(r;n+1+n+m) }
+\langle 1,\frac{\partial F}{\partial q}\rangle_{{\cal E}(1+n+m)}+
\langle 
\frac{\partial F}{\partial t}\rangle_{{\cal E}(m)}. 
\]
By the restriction of this to $q=z=0$, we have that
\begin{equation}
{\cal E}(r;n+m)=\langle H_0,x\frac{\partial
 H_0}{\partial x},\frac{\partial
 H_0}{\partial y}\rangle_{ {\cal
 E}(r;n+m) }
+\langle 1,y_1,\ldots,y_n,\frac{\partial
 H_0}{\partial t}\rangle_{{\cal E}(m)}.\label{rest:eqn}
\end{equation}
This means  that 
\begin{equation}
{\mathfrak M}(r;n+m)^{n+m+1}\subset \langle H_0,x\frac{\partial H_0}{\partial x},
\frac{\partial H_0}{\partial y}
\rangle_{{\cal E}(r;n+m)}+{\mathfrak M}(m){\cal E}(r;n+m).
\label{n+3det:leg}
\end{equation}
Let $C'\in C^{\Theta,Z}_T(J^1({\mathbb R}^m\times{\mathbb R}^n,{\mathbb R}),0)$ 
satisfying $j^{n+m+2}C(0)=j^{n+m+2}C'(0)$ be given.
There exists a function germ $H'(T,Q,p)\in {\mathfrak M}(2n+1)$ such that 
\[ P_{C'}=\{(T,Q,0,-\frac{\partial H'}{\partial T}(T,Q,p)+s,
-\frac{\partial H'}{\partial Q},T,-\frac{\partial H'}{\partial p},
H'-\langle \frac{\partial H'}{\partial p},p\rangle ,s,p)\}.\]
Since $H=z-qp$ on $P_C$ and $H'=z-qp$ on $P_{C'}$,
we have that $j^{n+m+2}H_0(0)=j^{n+m+2}H_0'(0)$, where 
$H_0'(x,y,t)=H'(t,x,0,y)\in {\mathfrak M}(r;n+m)^2$.
By (\ref{n+3det:leg}) we have that 
\[ {\mathfrak M}(r;n)^{n+m+1}\subset \langle H_0,x\frac{\partial H_0}{\partial x}
(x,y,0),
\frac{\partial H_0}{\partial y}(x,y,0)
\rangle_{{\cal E}(r;n)} \]
and this means that $H_0(x,y,0)$ is reticular ${\cal K}$-$(n+m+2)$-determined by Lemma \ref{findetc:lm}.
Therefore we may assume that 
$H_0|_{t=0}=H'_0|_{t=0}$. It follows that
$H_0-H_0'\in {\mathfrak M}(m){\mathfrak M}(r;n+m)^{n+m+2}$.
Then the function germ $G(x,y,t,q,z)=-z+H_0'(x,y,t)+
\langle y,q\rangle \in {\mathfrak M}(r;n+1+n+m)$
is a generating family of $C'|_{{\mathbb L}}$.

We define the function germ $E_{\tau_0}(x,y,t,\tau)\in {\cal E}(r;n+m+1)$ by
$E_{\tau_0}(x,y,t,\tau)=(1-\tau-\tau_0)H_0(x,y,t)+(\tau+\tau_0)H_0'(x,y,t)$
for $\tau_0\in [0,1]$.
By (\ref{rest:eqn}) and (\ref{n+3det:leg}), we have that
\begin{equation}
{\mathfrak M}(r;n+m)^{n+m+2}
\subset \langle H_0,x\frac{\partial H_0}{\partial x}
\rangle_{{\cal E}(r;n+m)}+{\mathfrak M}(r;n+m)\langle 
\frac{\partial H_0}{\partial y}\rangle 
+{\mathfrak M}(m)\langle 1,y,\frac{\partial H_0}{\partial t}\rangle.\label{n+3eqn:eqn}
\end{equation}
Then we have that 
\begin{eqnarray*}
& & {\mathfrak M}_t{\mathfrak M}_{x,y,t}^{n+m+2}{\cal E}_{x,y,t,\tau} \\
& \subset & 
{\mathfrak M}_{t,\tau}
\langle E_{\tau_0},x\frac{\partial E_{\tau_0}}{\partial x}\rangle_{{\cal E}_{x,y,t,\tau}}+{\mathfrak M}_{t,\tau}{\mathfrak M}_{x,y,t,\tau}\langle\frac{\partial E_{\tau_0}}{\partial y}
\rangle \\
& & \hspace{5cm}+{\mathfrak M}_{t,\tau}^2 \langle 1,y,\frac{\partial
 E_{\tau_0}}{\partial t}\rangle +{\mathfrak M}_{t,\tau}{\mathfrak M}_t
{\mathfrak M}_{x,y,t}^{n+m+2}{\cal E}_{x,y,t,\tau}.
\end{eqnarray*}
By Malgrange preparation theorem we have that 
\begin{eqnarray*}
\frac{\partial E_{\tau_0}}{\partial \tau}\in {\mathfrak M}_t{\mathfrak M}_{x,y,t}^{n+3}\subset 
{\mathfrak M}_t{\mathfrak M}_{x,y,t}^{n+3}{\cal E}_{x,y,t,\tau}\hspace{5cm}
\\
\subset  {\mathfrak M}_{t,\tau}(
\langle E_{\tau_0},x\frac{\partial E_{\tau_0}}{\partial x}\rangle_{{\cal E}_{x,y,t,\tau}}+
{\mathfrak M}_{x,y,t,\tau}\langle\frac{\partial E_{\tau_0}}{\partial y}
\rangle) 
+{\mathfrak M}_{t,\tau}^2 \langle 1,y,\frac{\partial
 E_{\tau_0}}{\partial t}\rangle.
\end{eqnarray*}
for $\tau_0\in [0,1]$.
Then there exist $\Phi(x,y,t)\in {\cal B}_m(r;n+m)$ and a unit $a\in
{\cal E}(r;n+m)$ and $b_1(t),\ldots,b_n(t),c(t)\in {\mathfrak M}(m)$ such that\\
(1) $\Phi$ has the form:
$
\Phi(x,y,t)=(x\phi_1(x,y,t),\phi_2(x,y,t),\phi_3(t))
$,\\
(2) $H_0(x,y,t)=a(x,y,t)\cdot H_0'\circ \Phi(x,y,t)+\sum_{i=1}^ny_ib_i(t)+c(t)$ 
for $(x,y,t)\in ({\mathbb H}^r\times {\mathbb R}^{n+m},0)$

We define the reticular $t$-${\cal P}$-${\cal K}$-isomorphism $(\Psi,d)$
by 
\[
\Psi(x,y,t,q,z)=(x\phi_1(x,y,t),\phi_2(x,y,t),\phi_3(t),q(1-b(t)),z),
d(x,y,t,q,z)=a(x,y,t).
\]
We set $G':=d\cdot G\circ\Psi\in {\mathfrak M}(r;n+n+m)$.
Since $\frac{\partial E_{\tau_0}}{\partial \tau}|_{t=0}=0$, we have that 
$a(x,y,0)=1$ and $\Phi(x,y,0)=(x,y,0)$.
Therefore we have that $G'|_{t=0}=f$.
Then $F$ and $G'$ are reticular $t$-${\cal P}$-${\cal K}$-infinitesimal versal
unfoldings of $F|_{t=0}$.
Since $G$ and $G'$ are reticular $t$-${\cal P}$-${\cal K}$-equivalent,
it follows that $F$ and $G$ are reticular $t$-${\cal P}$-${\cal K}$-equivalent.
Therefore ${\cal L}$ and $C'|_{{\mathbb L}}$ are
${\cal P}$-Legendrian equivalent.\hfill $\blacksquare$\\

Let ${\cal L}$ be a stable reticular Legendrian unfolding.
We say that ${\cal L}$ is {\it simple} if 
there exists a representative $\tilde{C}\in C_T(
U,J^1({\mathbb R}^m\times{\mathbb R}^n,{\mathbb R}))$ of
 a extension of ${\cal L}$ such that 
$\{ \tilde{C}_x| x\in  U\}$ is covered by 
finite orbits $[C_1],\ldots,[C_l]$ for  
 some ${\cal P}$-contact embedding germs $C_1,\ldots,C_l\in 
C_T(U,J^1({\mathbb R}^m\times{\mathbb R}^n,{\mathbb R}))$. 
 \begin{lem}\label{simpleLeg:lem}{\rm (cf., \cite[Proposition 5.5]{bifsemi})}
A stable reticular Legendrian  unfolding ${\cal L}$ is simple if and only if 
for a generating family 
$F(x,y,t,q,z)\in {\mathfrak M}(r;k+m+n+1)$ of ${\cal L}$,
$f(x,y)=F(x,y,0,0)\in {\mathfrak M}(r;k)^2$ is a  ${\cal K}$-simple singularity.
\end{lem}
Let  $J^l(2n+2m+1,2n+2m+1)$ be the set of $l$-jets of map germs  from  $(J^1({\mathbb R}^m
\times{\mathbb R}^n,{\mathbb R}),0)$ 
to $(J^1({\mathbb R}^m\times{\mathbb R}^n,{\mathbb R}),0)$ and 
$tC^l(n)$ be the immersed manifold in $J^l(2n+2m+1,2n+2m+1)$ 
which consists of $l$-jets of ${\cal P}$-contact embedding  germs.
Let $L^l(2n+2m+1)$ be the Lie group which consists of $l$-jets of diffeomorphism germs on 
$(J^1({\mathbb R}\times{\mathbb R}^n,{\mathbb R}),0)$.

We  consider the Lie subgroup $rtLe^l(2n+2m+1)$ of $L^l(2n+2n+1)\times L^l(2n+2m+1)$ which consists of 
 $l$-jets of reticular ${\cal P}$-diffeomorphisms on the  source space and $l$-jets of ${\cal P}$-Legendrian 
equivalences of $\Pi$ at $0$: 
\begin{eqnarray*}
rtLe^l(n,m)=\{ (j^l\Psi(0),j^lK(0))\in 
L^l(2n+2m+1)\times L^l(2n+2m+1)\ | \\
 \Psi_i
 \mbox{ is a reticular} \mbox{ ${\cal P}$-diffeomorphism on }
(J^1({\mathbb R}\times{\mathbb R}^n,{\mathbb R}),0), \\
K \mbox{ is a ${\cal P}$-Legendrian equivalence of } 
\Pi \}.
\end{eqnarray*}

The group $rtLe^l(2n+2m+1)$ acts on $J^l(2n+2m+1,2n+2m+1)$ and $tC^l(2n+2m+1)$ 
is invariant 
under this action.

Let $C$ be  a ${\cal P}$-contact diffeomorphism germ on 
$(J^1({\mathbb R}\times{\mathbb R}^n,{\mathbb R}),0)$
and set  
$z=j^lC(0)$, ${\cal L}=C|_{{\mathbb L}}$.
We denote  the orbit $rtLe^l(2n+2m+1)\cdot z$ by $[z]$.
Then 
\[
[z]=\{ j^lC'(0)\in tC^l(2n+2m+1)\ | \ {\cal L}  \mbox{ and } 
C'|_{{\mathbb L}}
\mbox{ are ${\cal P}$-Legendrian equivalent} \}. \]

For $\tilde{C}=\in C_T(U,J^1({\mathbb R}^m\times{\mathbb R}^n,{\mathbb R}))$,
we define
 the continuous map 
$j^l_0\tilde{C}:U\rightarrow tC^l(n)$ by $x$ to the $l$-jet of 
$\tilde{C}_{x}$.
For $C\in C_T(J^1({\mathbb R}^m\times{\mathbb R}^n,{\mathbb R}),0)$,
we define $j^l_0C:(J^1({\mathbb R}^m\times{\mathbb R}^n,{\mathbb R}),0)\rightarrow 
 tC^l(n)$ by the analogous method.
\begin{tth}\label{stabletrans_tleg2:th}{\rm (cf., \cite[Theorem 5.4]{bifsemi})}
Let ${\cal L}$ be a reticular Legendrian unfolding. 
Let $C$ be an extension of ${\cal L}$ and $l\geq (n+m+1)^2$.
Then the followings are equivalent:\\
{\rm (s)} ${\cal L}$  is stable.\\
{\rm (t)} $j^l_0C$ is transversal to $[j^l_0C(0)]$.\\
{\rm (a')} ${\cal E}_{t,q,p}=
B_0+
\langle 1,p_1\circ C',\ldots,p_n\circ C'\rangle_{(\Pi\circ C')^*{\cal E}_{t,q,z}}+
\langle s\circ C'\rangle_{{\cal E}_t}+{\mathfrak M}_{t,q,p}^l$,
where $C'=C|_{z=s=0}$ and $B_0=\langle q_1p_1,\ldots,q_{r}p_{r},
q_{r+1},$ $\ldots,q_n\rangle_{{\cal E}_{t,q,p}}$,\\
{\rm (a)} ${\cal E}_{t,q,p}=
B_0+
\langle 1,p_1\circ C',\ldots,p_n\circ C'\rangle_{(\Pi\circ C')^*{\cal E}_{t,q,z}}+
\langle s\circ C'\rangle_{{\cal E}_t}$,\\
{\rm (is)} ${\cal L}$ is infinitesimally stable,\\
{\rm (hs)} ${\cal L}$ is homotopically stable,\\
{\rm (u)} A generating family $F$ of ${\cal L}$ is reticular 
$t$-${\cal P}$-${\cal K}$-stable unfolding of $F|_{t=0}$.
\end{tth}
{\it Proof}. We prove only  (a')$\Rightarrow$(a)
By the restriction of  (a') to $t=0$ we have that:
\[{\cal E}_{q,p}=
B_1+
\langle 1,p_1\circ C'',\ldots,p_n\circ C''\rangle_{(\Pi\circ C'')^*{\cal E}_{t,q,z}}
+
\langle s\circ C''\rangle_{{\mathbb R}}+
{\mathfrak M}_{q,p}^l,
\]
where $C''=C'|_{t=0}$ and $B_1=B_0|_{t=0}$.
Then we have that 
\[{\cal E}_{q,p}=
B_1+
(\Pi\circ C'')^*
{\mathfrak M}_{t,q,p}{\cal E}_{q,p(m)}
+
\langle 1,p_1\circ C'',\ldots,p_n\circ C'',
s\circ C''\rangle_{{\mathbb R}}+{\mathfrak M}_{q,p}^l.
\]
It follows that 
\[ {\mathfrak M}_{q,p}^{n+m+1}\subset 
B_1+(\Pi\circ C'')^*
{\mathfrak M}_{t,q,p}{\cal E}_{q,p}.\]
Therefore 
\[ {\mathfrak M}_{t,q,p}^{n+m+1}\subset B_0\times+
(\Pi\circ C')^*
{\mathfrak M}_{t,q,p}{\cal E}_{q,p}+
{\mathfrak M}_{t}{\cal E}_{t,q,p},\]
and we have that 
\[{\mathfrak M}_{t,q,p}^l=({\mathfrak M}_{t,q,p(m)}^{n+m+1})^{n+m+1}
\subset 
B_0+(\Pi\circ C')^*
{\mathfrak M}_{t,q,p}^{n+m+1}{\cal E}_{t,q,p}+
{\mathfrak M}_{t}{\cal E}_{t,q,p}.
\]
It follows that
\[
{\cal E}_{t,q,p}=
B_0+
\langle 1,p_1\circ C',\ldots,p_n\circ C'\rangle_{(\Pi\circ C')^*{\cal E}_{t,q,z}}+
\langle s\circ C'\rangle_{{\cal E}_t}+
(\Pi\circ C')^*
{\mathfrak M}_{t,q,p}^{n+m+1}{\cal E}_{t,q,p}+
{\mathfrak M}_{t}{\cal E}_{t,q,p}.
\]
This means (a) by Lemma \ref{gw1.8:cor}.\\

\begin{tth}{\rm (cf., \cite[Theorem 5.6]{bifsemi})}
Let $r=0,n\leq 5,m=2,3$ or $r=1,n\leq 3,m=2,3$. 
Let  $U$ be a neighborhood of $0$ in 
$J^1({\mathbb R}^m\times{\mathbb R}^n,{\mathbb R})$.
Then there exists a residual set $O\subset  C^\Theta_T(U,J^1({\mathbb R}^m\times{\mathbb R}^n,{\mathbb R}))$ 
such that for any $\tilde{C}\in O$ and $w\in U$,
the reticular Legendrian unfolding $\tilde{C}_w|_{{\mathbb L}}$ 
is stable and has a generating family which is stably reticular 
$t$-${\cal P}$-${\cal K}$-equivalent to one of the types in the 
classification list below.
\end{tth}

Let $F(x,y,u,t)\in {\mathfrak M}(r;k+n+2)$ be a reticular $t$-${\cal P}$-${\cal K}$-stable unfolding of $f(x,y)$ be given
($r=0,n\leq 6$ or $r=1,n\leq 4$).
Since $f$ is simple singularity we may assume that $f$ has the normal form of $A,D,E(r=0)$ or $B,C,F(r=1)$.
By analogous method of \cite[p.200]{tPKfunct}, we may assume that $F$ has the form
\[ F(x,y,u,t)=f(x,y)+a(u_l,\ldots ,u_n,t)\varphi_0(x,y)+u_1\varphi_1(x,y)+\cdots + u_{l-1}\varphi_{l-1}(x,y), \]
where the function germ $f(x,y)+t\varphi_0(x,y)+u_1\varphi_1(x,y)+\cdots + u_{l-1}\varphi_{l-1}(x,y)\in  {\mathfrak M}(r;k+n+1)$ is a reticular $t$-${\cal P}$-${\cal K}$-universal unfolding.
Since $F$ is also a  reticular $t$-${\cal P}$-${\cal K}$-universal unfolding of $f$, we have that 
\[ {\cal E}_{x,y,u,t}=\langle f,x\frac{\partial f}{\partial x},\frac{\partial f}{\partial y}\rangle_{{\cal E}_{x,y,u,t}} 
+\langle \varphi_1,\ldots,\varphi_{l-1}\rangle_{{\cal E}_{u,t}}+
\varphi_0(\langle \frac{\partial a}{\partial u_l},\ldots,\frac{\partial a}{\partial u_n}
\rangle_{{\cal E}_{u_l,\ldots,u_n,t}}
+\langle \frac{\partial a}{\partial t}\rangle_{{\cal E}_{t}}).
\]
This means that $a(u_l,\ldots,u_n,t)$ is a ${\cal P}$-${\cal R}$-versal unfolding of $a(u_l,\ldots,n,0)$ with 
codimension $\leq 3$.
Since the ${\cal P}$-${\cal R}$-equivalence of $a$ is allowed under the reticular 
$t$-${\cal P}$-${\cal K}$-equivalence, 
it follows the classification of  ${\cal P}$-${\cal R}$-versal unfolding of functions on $u_l,\ldots,u_n$
of type $A_2,A_3$.\\

We classify $F(x,y,q,t)\in {\mathfrak M}(r;k+n+m)$ with $r=0,n\leq 6,m=2$ and $r=1,n\leq 4,m=2$.\\
$({}^2A_1$)  $y^2+(t_1+t_2u_1+u_1^3\pm u_2^2\pm \ldots \pm u_l^2)$,\\
$({}^2A_2$)  $y^3+(t_1+t_2u_2+u_2^3\pm u_3^2\pm \ldots \pm u_l^2)y+u_1$,\\
$({}^2A_3$)  $y^4+(t_1+t_2u_3+u_3^3\pm u_4^2\pm \ldots \pm u_l^2)y^2+u_1y+u_2$,\\
$({}^2A_4$)  $y^5+(t_1+t_2u_4+u_4^3\pm u_5^2\pm \ldots \pm u_l^2)y^3+u_1y^2+u_2y+u_3$,\\
$({}^2A_5$)  $y^6+(t_1+t_2u_5+u_5^3)y^4+u_1y^3+u_2y^2+u_3y+u_4$,  
                       $y^6+(t_1+t_2u_5+u_5^3\pm u_6^2)y^4+u_1y^3+u_2y^2+u_3y+u_4$,\\
$({}^2A_6$)  $y^7+(t_1+t_2u_6+u_6^3)y^5+u_1y^4+u_2y^3+u_3y^2+u_4y+u_5$,\\
$({}^2D^\pm_4)$   $y_1^2y_2\pm y_2^3+(t_1+t_2u_4+u_4^3\pm u_5^2\pm \ldots \pm u_l^2)y_2^2+u_1y_2+u_2y_1+u_3$,\\
$({}^2D_5)$  $y_1^2y_2+ y_2^4+(t_1+t_2u_5+u_5^3)y_2^3+u_1y_2^2+u_2y_2+u_3y_1+u_4$,
                        $y_1^2y_2+ y_2^4+(t_1+t_2u_5+u_5^3\pm u_6^2)y_2^3+u_1y_2^2+u_2y_2+u_3y_1+u_4$,\\
$({}^2D^\pm_6)$ $y_1^2y_2\pm  y_2^5+(t_1+t_2u_6+u_6^3)y_2^6+u_1y_2^3+u_2y_2^2+u_3y_2+u_4y_1+u_5$,\\
$({}^2E_6)$ $y_1^3+ y_2^4+(t_1+t_2u_6+u_6^3)y_1y_2^2+u_1y_1y_2+u_2y_2^2+u_3y_1+u_4y_2+u_5$\\
,where $l\leq 6$.\\
$({}^2B_2)$  $x^2+(t_1+t_2u_2+u_2^3\pm u_2^2\pm \ldots \pm u_l^2)x+u_1$,\\
$({}^2B_3)$  $x^3+(t_1+t_2u_3+u_3^3)x+u_1x+u_2$,
                        $x^3+(t_1+t_2u_3+u_3^3\pm u_4^2)x+u_1x+u_2$,\\
$({}^2B_4)$  $x^4+(t_1+t_2u_4+u_4^3)x^2+u_1x^2+u_2x+u_3$,\\
$({}^2C_3^\pm$) $\pm xy+y^3+(t_1+t_2u_3+u_3^2)x+u_1y+u_2$,
                                $\pm xy+y^3+(t_1+t_2u_3+u_3^2\pm u_4^2)x+u_1y+u_2$,\\
$({}^2C_4)$ $xy+y^4+(t_1+t_2u_4+u_4^3)y^3+u_1y^2+u_2y+u_3$,\\
$({}^2F_4)$ $x^2+y^3+(t_1+t_2u_4+u_4^3)xy+u_1x+u_2y+u_3$\\
,where $l\leq 4$.
We give all figures of bifurcations of generic wavefronts with $n=2,3$ and $m=2$: $({}^2A_1$), $({}^2A_2$), $({}^2A_3$), 
$({}^2B_2)$, $({}^2B_3)$, $({}^2C_3^\pm$).\\
We give the positions of the figures as the following way:
\newpage
\begin{figure}[ht]
\begin{center}
 \begin{minipage}{0.60\hsize}
  \begin{center}
    \includegraphics*[width=10cm,height=10cm]{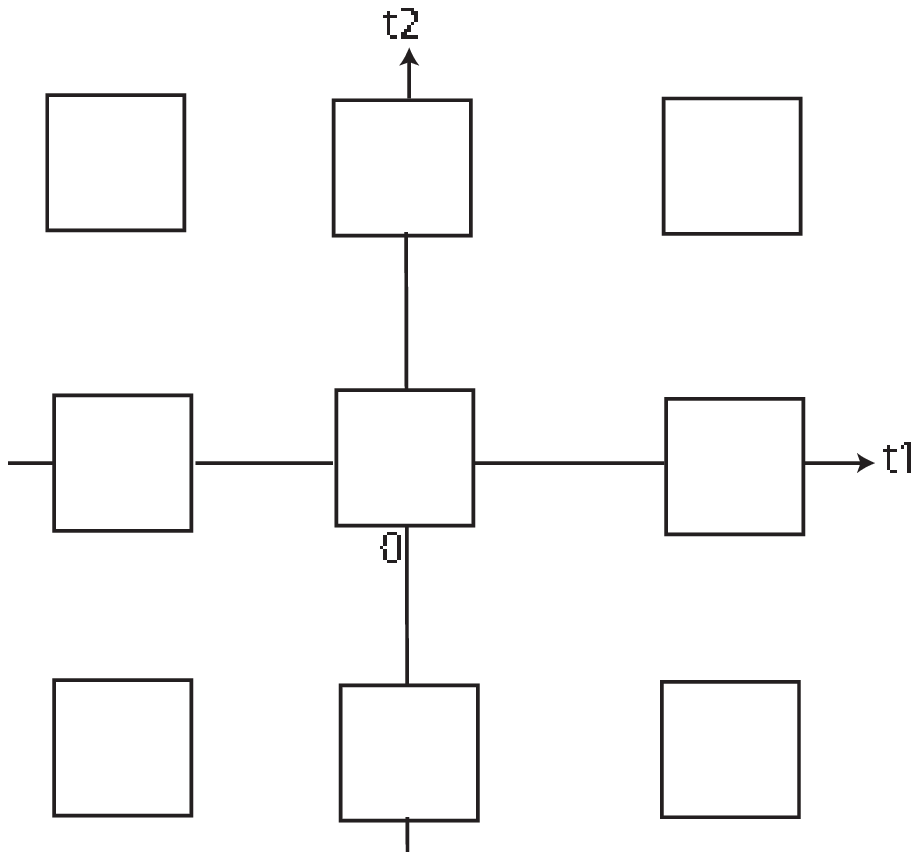} 
\end{center}
 \end{minipage}
\end{center}
\end{figure}

\newpage
\begin{figure}[ht]
\begin{center}
 \begin{minipage}{0.30\hsize}
  \begin{center}
    \includegraphics*[width=4cm,height=4cm]{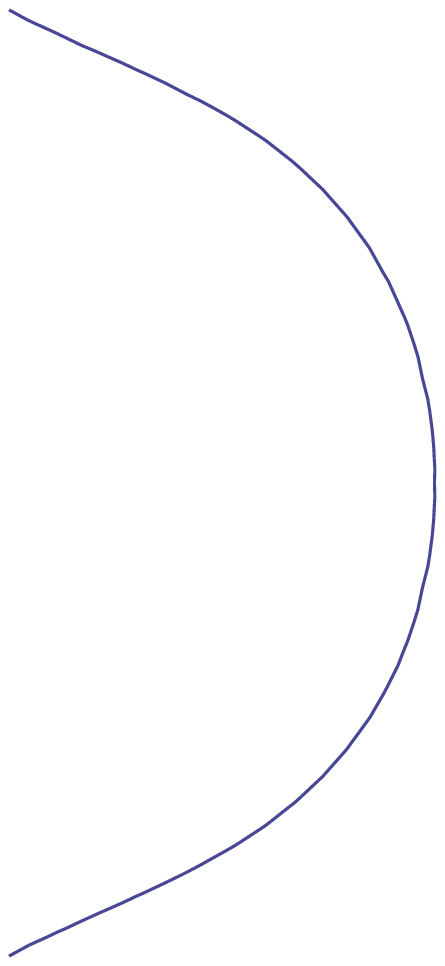} 
\end{center}
 \end{minipage}
 \begin{minipage}{0.30\hsize}
  \begin{center}
    \includegraphics*[width=4cm,height=4cm]{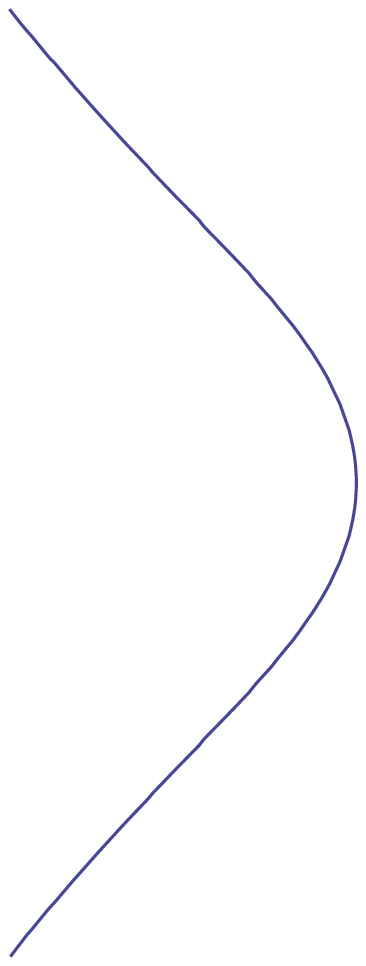}
\end{center}
 \end{minipage}
  \begin{minipage}{0.30\hsize}
  \begin{center}
   \includegraphics*[width=4cm,height=4cm]{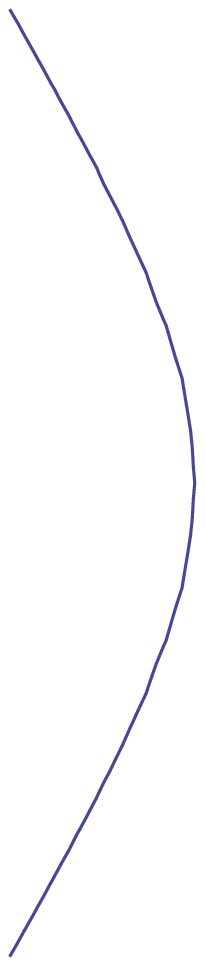}
\end{center}
 \end{minipage}
 \end{center}
 \end{figure}
 \begin{figure}[ht]
\begin{center}
 \begin{minipage}{0.30\hsize}
  \begin{center}
    \includegraphics*[width=4cm,height=4cm]{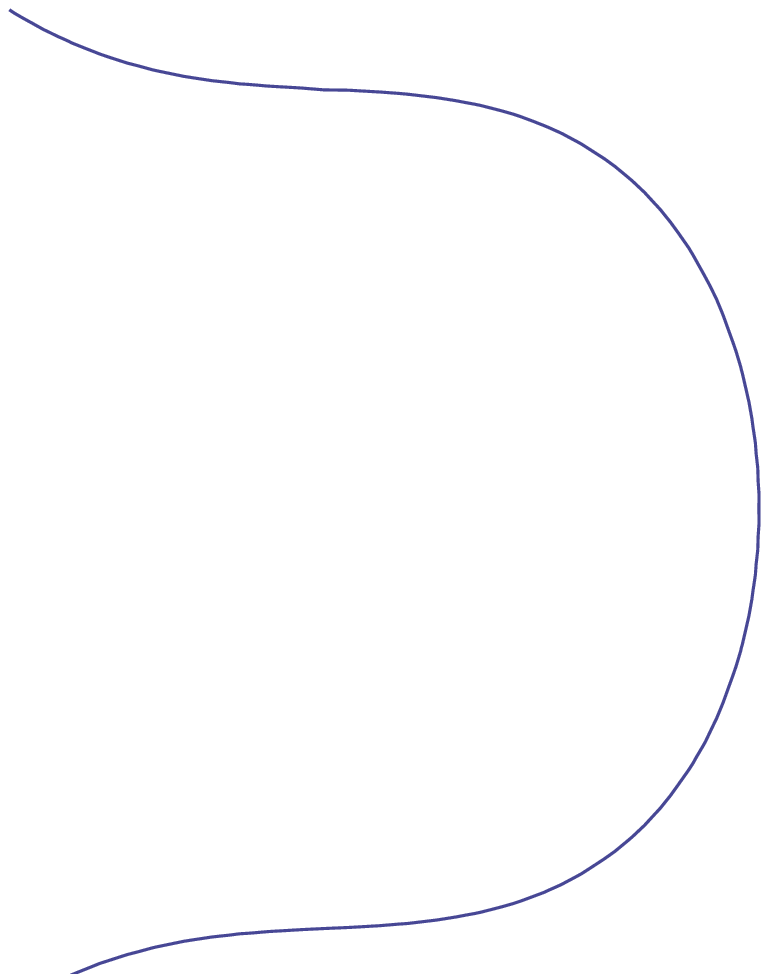} 
\end{center}
 \end{minipage}
 \begin{minipage}{0.30\hsize}
  \begin{center}
    \includegraphics*[width=4cm,height=4cm]{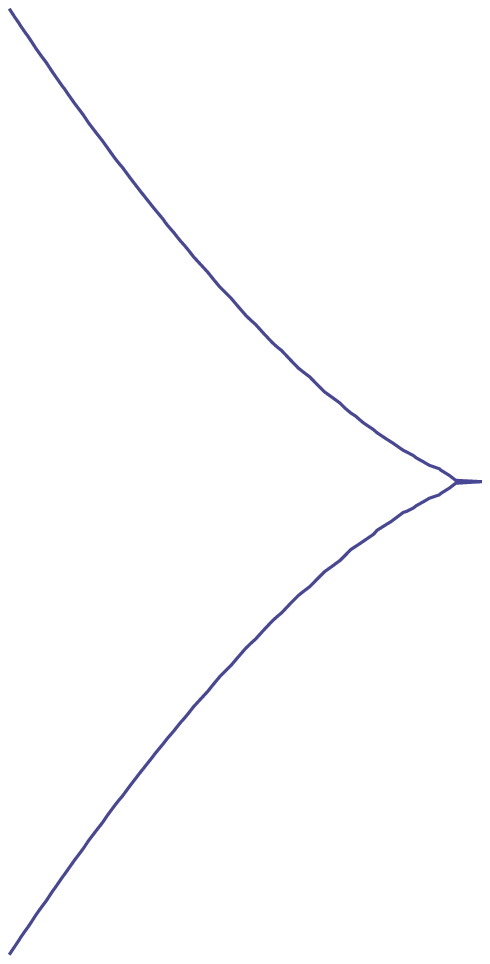}
\end{center}
 \end{minipage}
  \begin{minipage}{0.30\hsize}
  \begin{center}
   \includegraphics*[width=4cm,height=4cm]{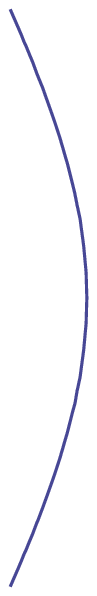}
\end{center}
 \end{minipage}
 \end{center}
 \end{figure}
 \begin{figure}[ht]
\begin{center}
 \begin{minipage}{0.30\hsize}
  \begin{center}
    \includegraphics*[width=4cm,height=4cm]{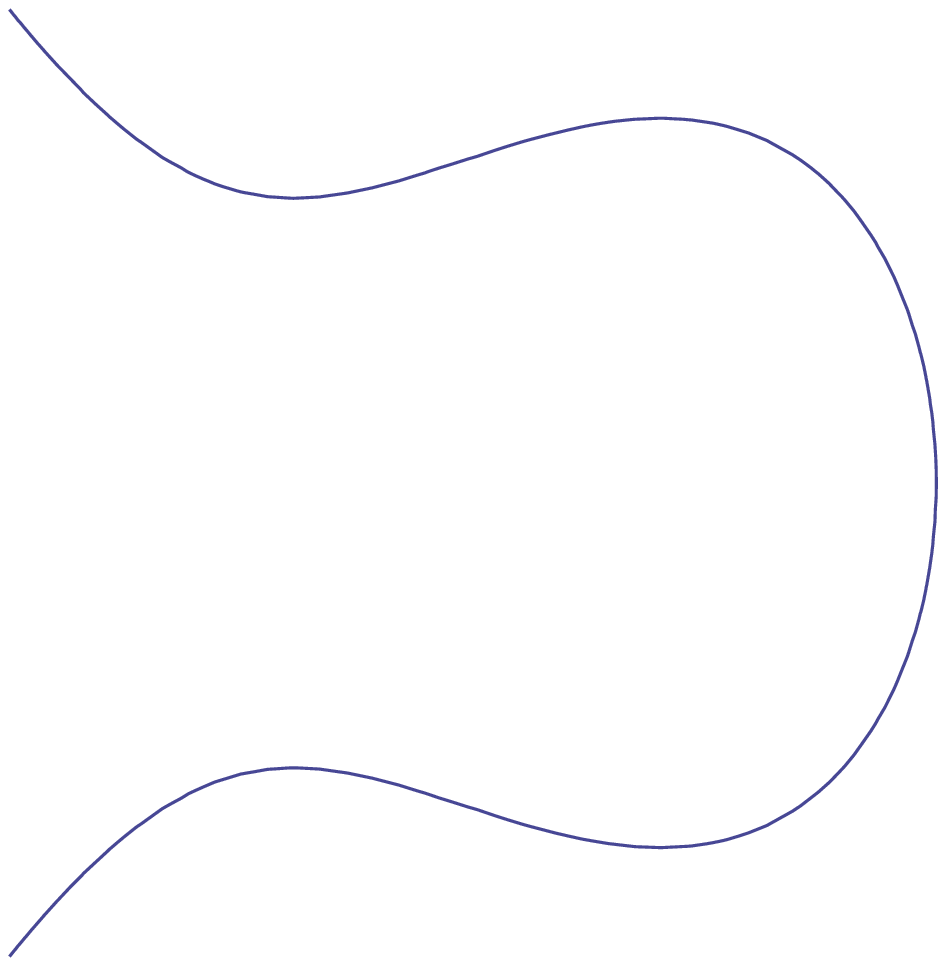} 
\end{center}
 \end{minipage}
 \begin{minipage}{0.30\hsize}
  \begin{center}
    \includegraphics*[width=4cm,height=4cm]{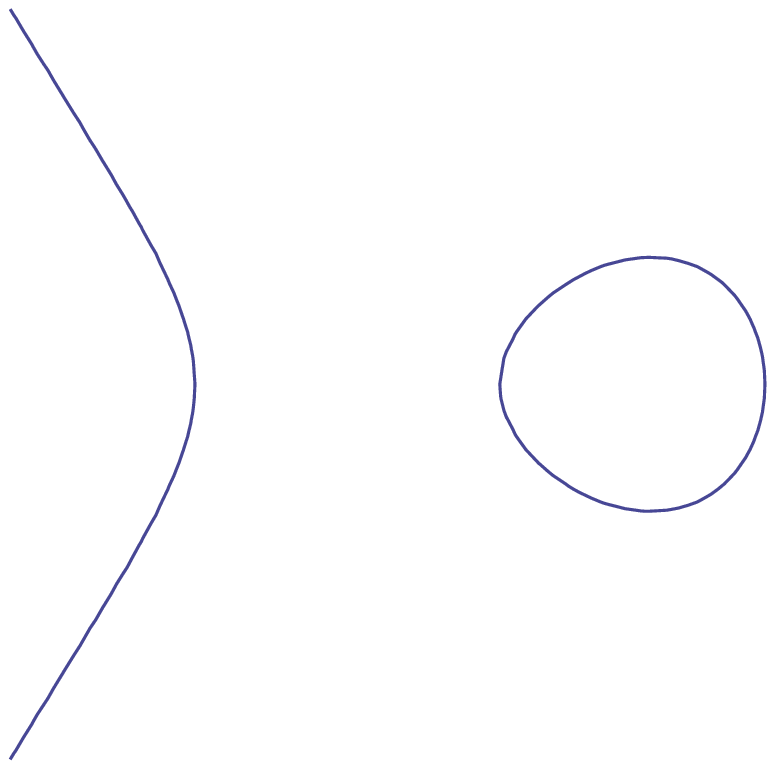}
\end{center}
 \end{minipage}
  \begin{minipage}{0.30\hsize}
  \begin{center}
   \includegraphics*[width=4cm,height=4cm]{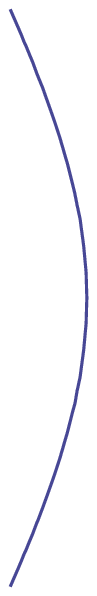}
\end{center}
 \end{minipage}
\end{center}
\caption{${}^2A_1$}
\end{figure}
 
 \newpage
\begin{figure}[ht]
\begin{center}
 \begin{minipage}{0.30\hsize}
  \begin{center}
    \includegraphics*[width=4cm,height=4cm]{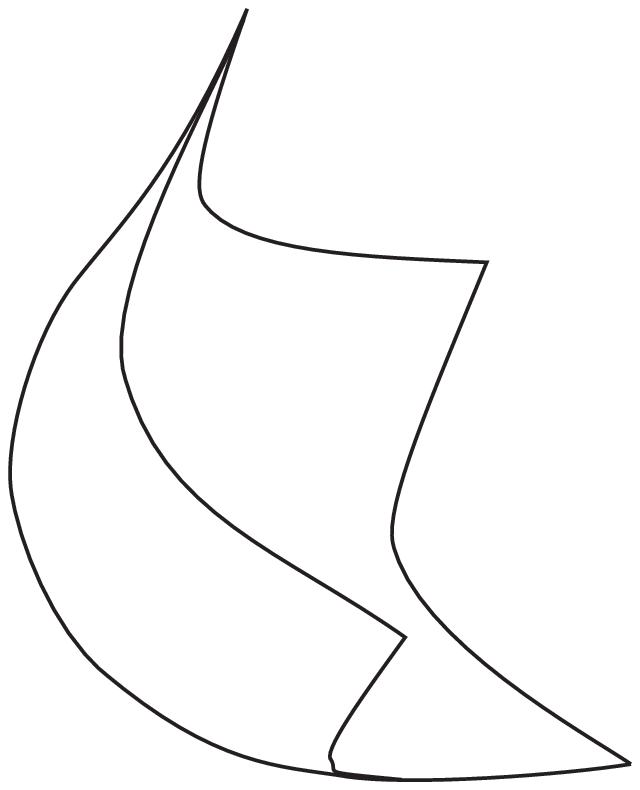} 
\end{center}
 \end{minipage}
 \begin{minipage}{0.30\hsize}
  \begin{center}
    \includegraphics*[width=4cm,height=4cm]{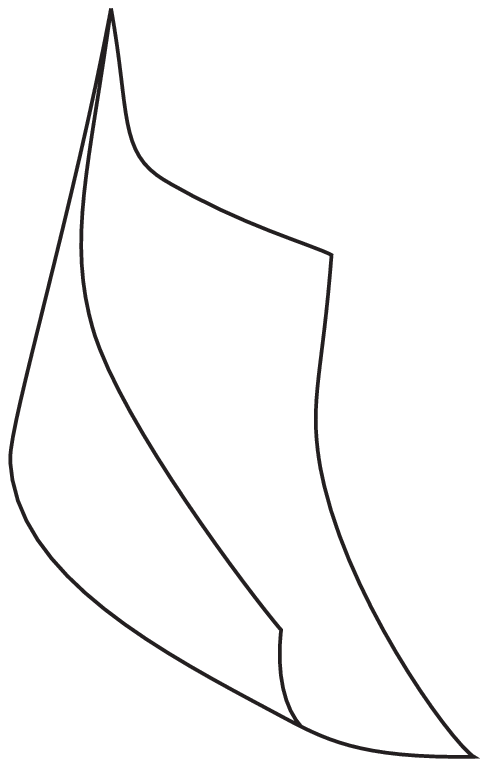}
\end{center}
 \end{minipage}
  \begin{minipage}{0.30\hsize}
  \begin{center}
   \includegraphics*[width=4cm,height=4cm]{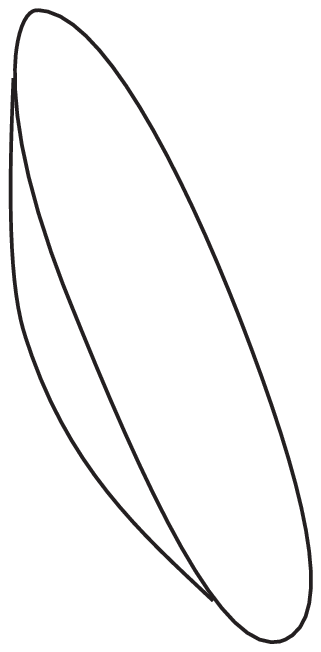}
\end{center}
 \end{minipage}
 \end{center}
 \end{figure}
 \begin{figure}[ht]
\begin{center}
 \begin{minipage}{0.30\hsize}
  \begin{center}
    \includegraphics*[width=4cm,height=4cm]{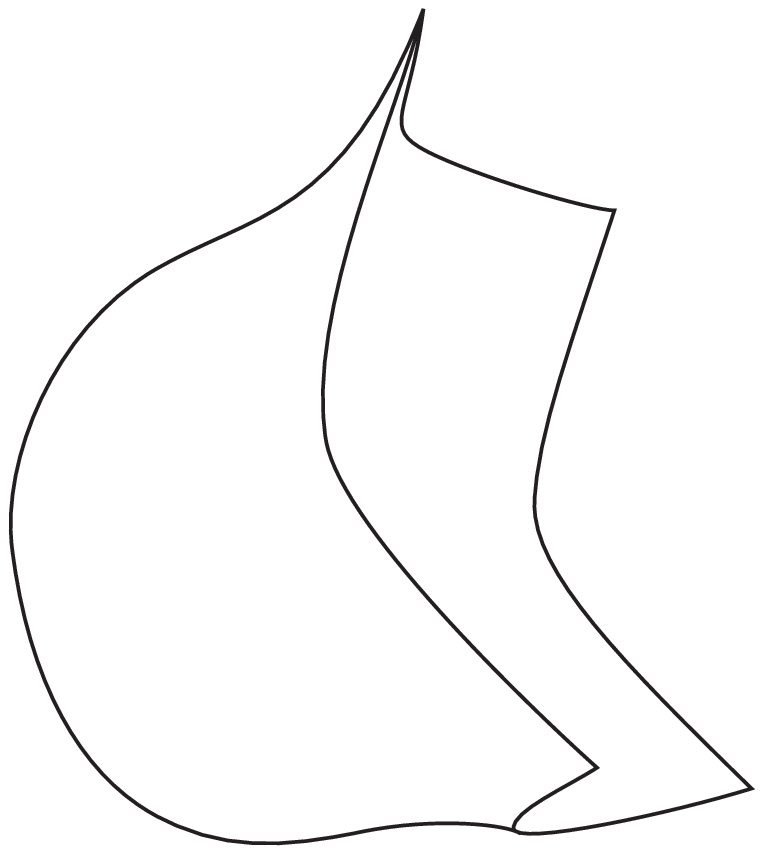} 
\end{center}
 \end{minipage}
 \begin{minipage}{0.30\hsize}
  \begin{center}
    \includegraphics*[width=4cm,height=4cm]{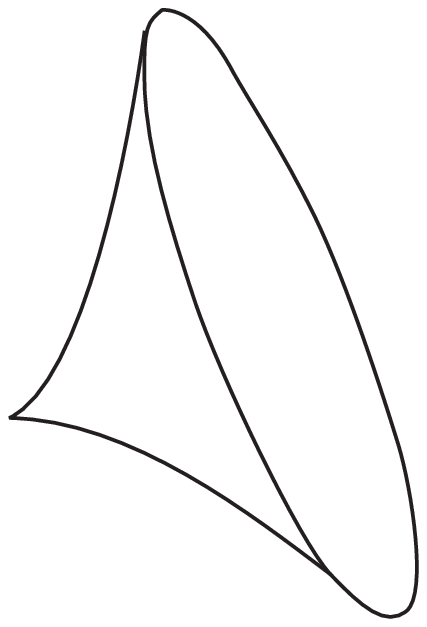}
\end{center}
 \end{minipage}
  \begin{minipage}{0.30\hsize}
  \begin{center}
   \includegraphics*[width=4cm,height=4cm]{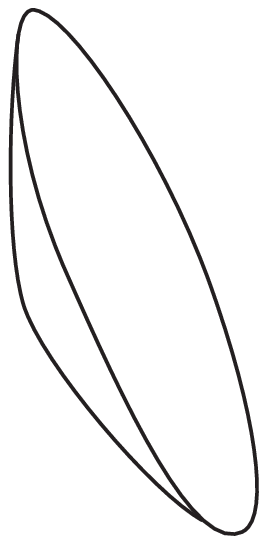}
\end{center}
 \end{minipage}
 \end{center}
 \end{figure}
 \begin{figure}[ht]
\begin{center}
 \begin{minipage}{0.30\hsize}
  \begin{center}
    \includegraphics*[width=4cm,height=4cm]{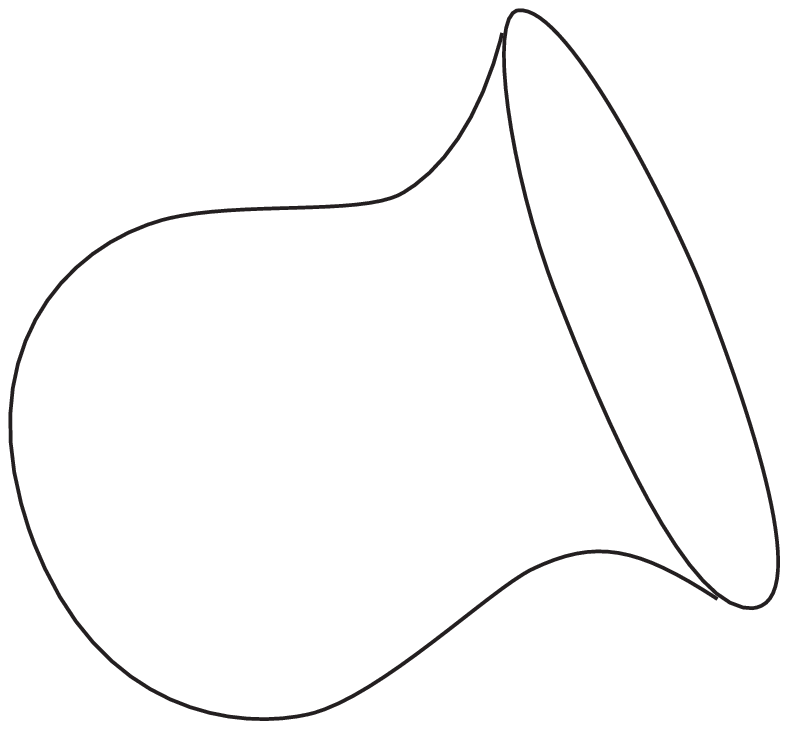} 
\end{center}
 \end{minipage}
 \begin{minipage}{0.30\hsize}
  \begin{center}
    \includegraphics*[width=4cm,height=4cm]{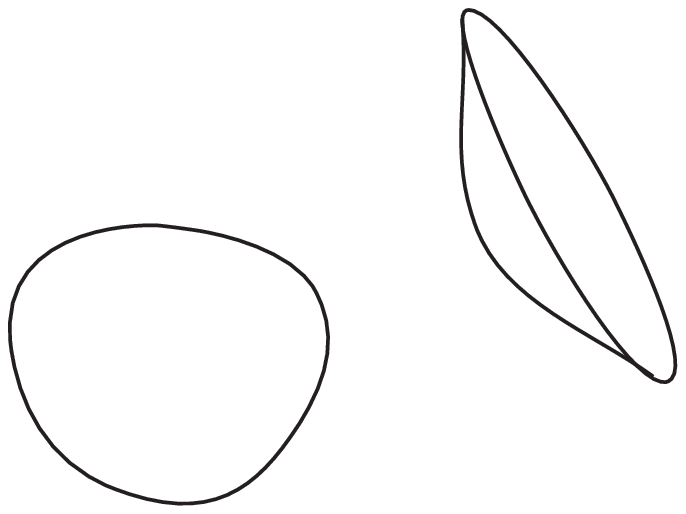}
\end{center}
 \end{minipage}
  \begin{minipage}{0.30\hsize}
  \begin{center}
   \includegraphics*[width=4cm,height=4cm]{2A1d+-.eps}
\end{center}
 \end{minipage}
\end{center}
\caption{${}^2A_1$}
\end{figure}
 
 \newpage
\begin{figure}[ht]
\begin{center}
 \begin{minipage}{0.30\hsize}
  \begin{center}
    \includegraphics*[width=4cm,height=4cm]{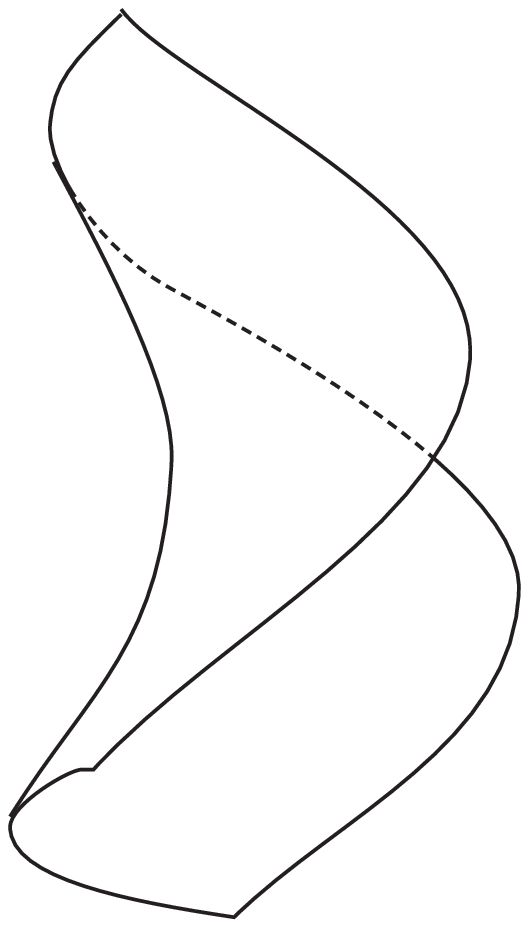} 
\end{center}
 \end{minipage}
 \begin{minipage}{0.30\hsize}
  \begin{center}
    \includegraphics*[width=4cm,height=4cm]{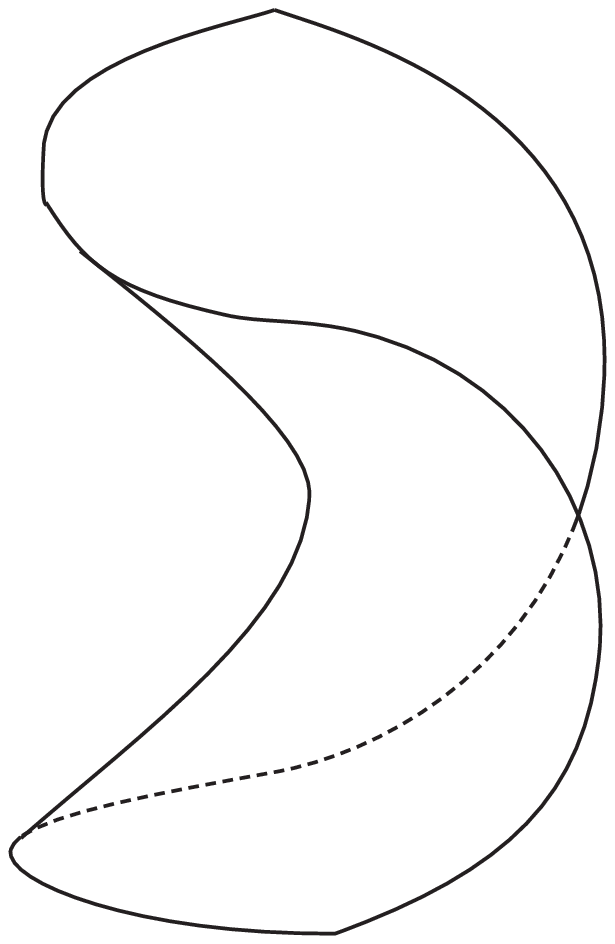}
\end{center}
 \end{minipage}
  \begin{minipage}{0.30\hsize}
  \begin{center}
   \includegraphics*[width=4cm,height=4cm]{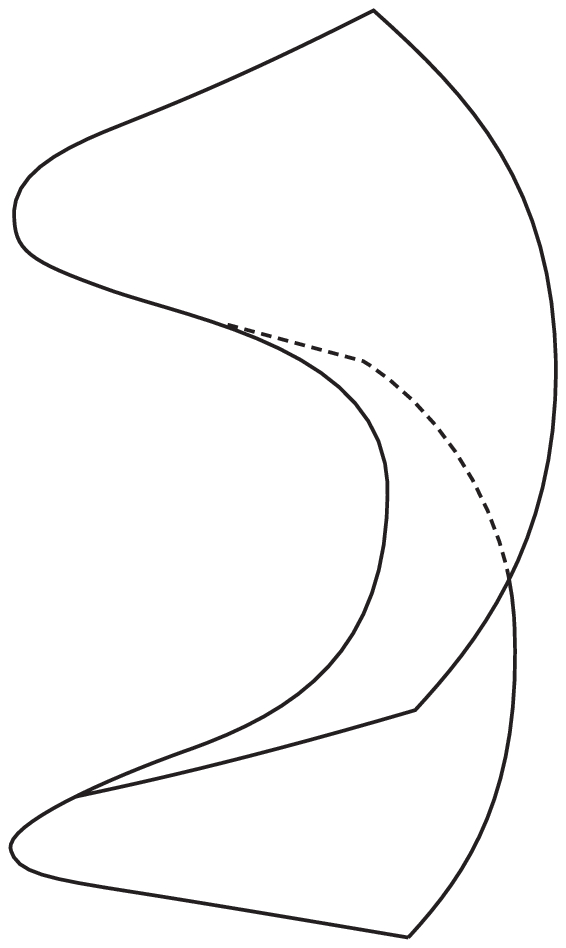}
\end{center}
 \end{minipage}
 \end{center}
 \end{figure}
 \begin{figure}[ht]
\begin{center}
 \begin{minipage}{0.30\hsize}
  \begin{center}
    \includegraphics*[width=4cm,height=4cm]{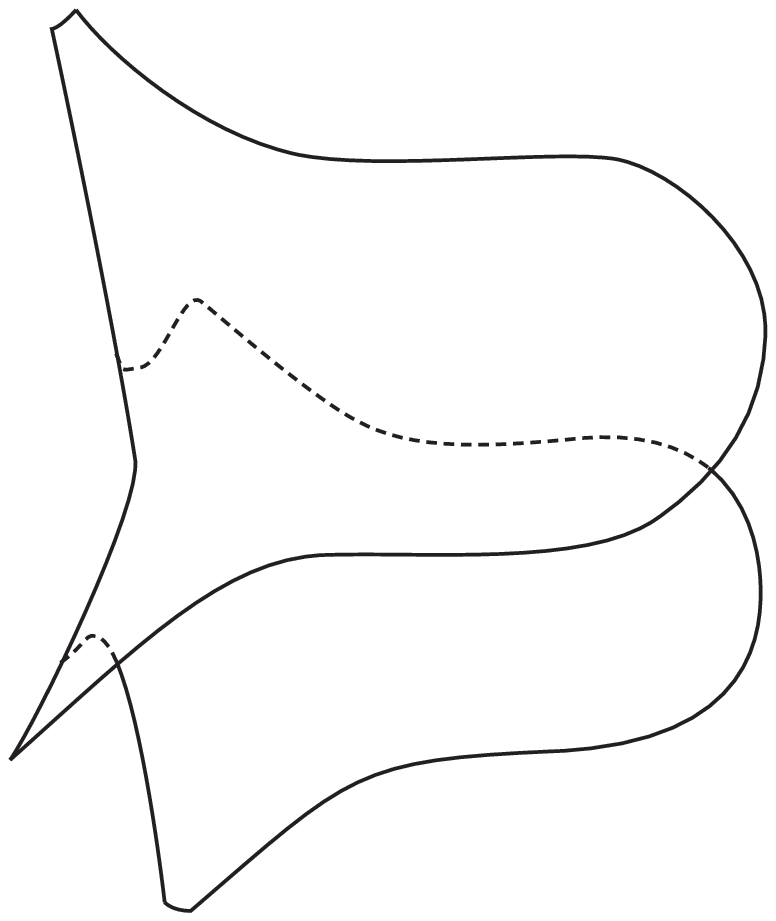} 
\end{center}
 \end{minipage}
 \begin{minipage}{0.30\hsize}
  \begin{center}
    \includegraphics*[width=4cm,height=4cm]{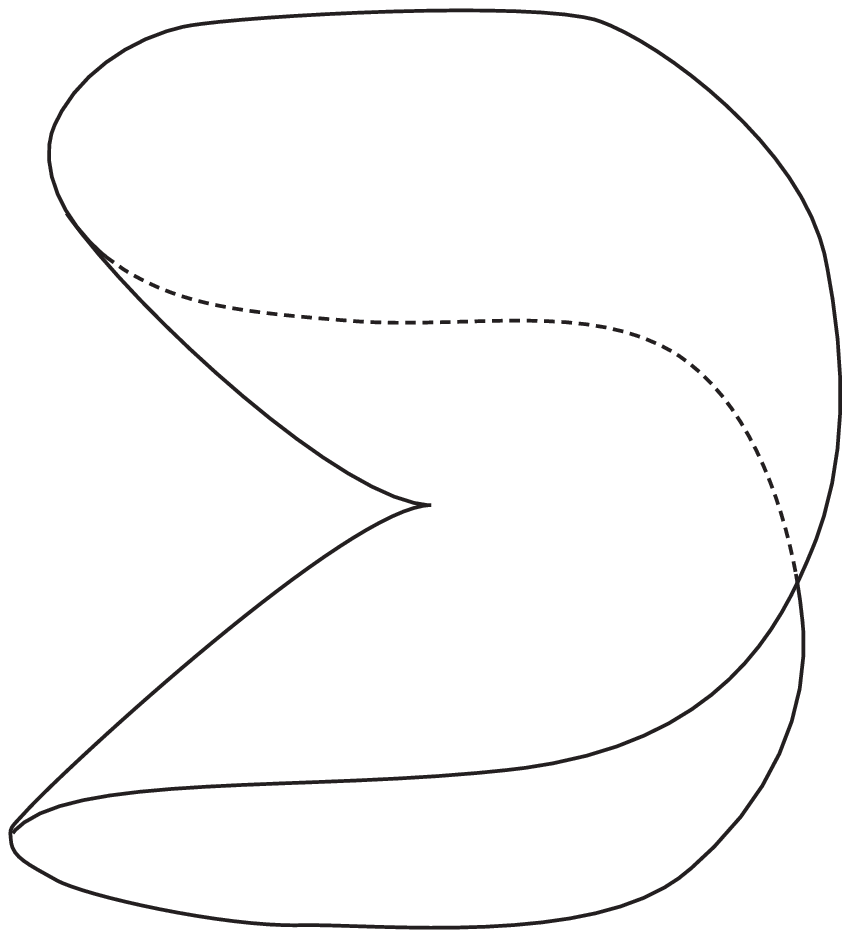}
\end{center}
 \end{minipage}
  \begin{minipage}{0.30\hsize}
  \begin{center}
   \includegraphics*[width=4cm,height=4cm]{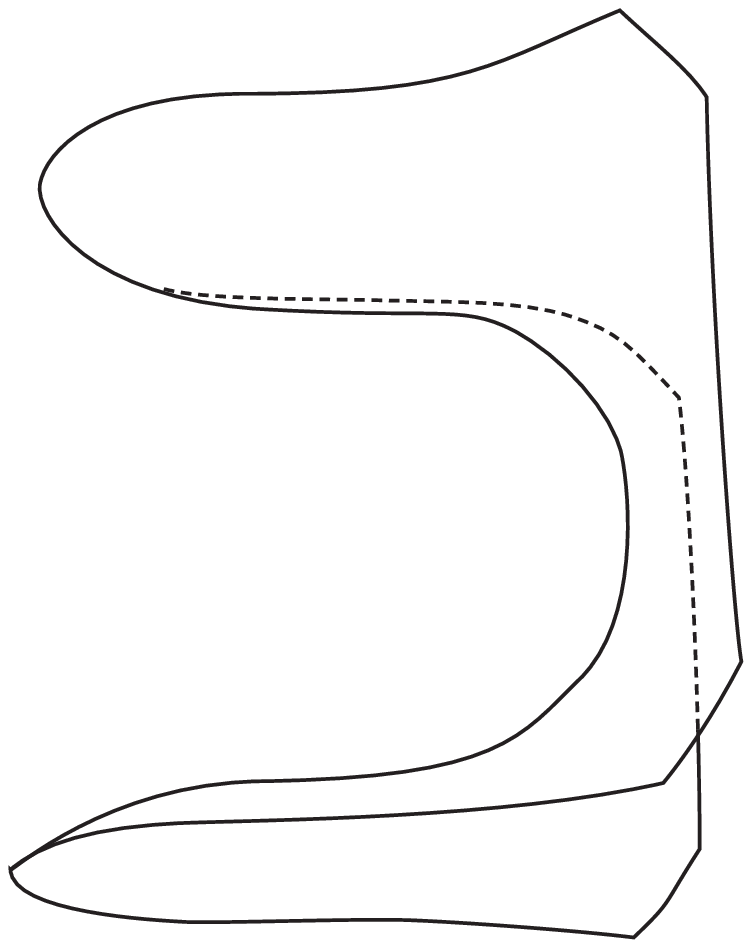}
\end{center}
 \end{minipage}
 \end{center}
 \end{figure}
 \begin{figure}[ht]
\begin{center}
 \begin{minipage}{0.30\hsize}
  \begin{center}
    \includegraphics*[width=4cm,height=4cm]{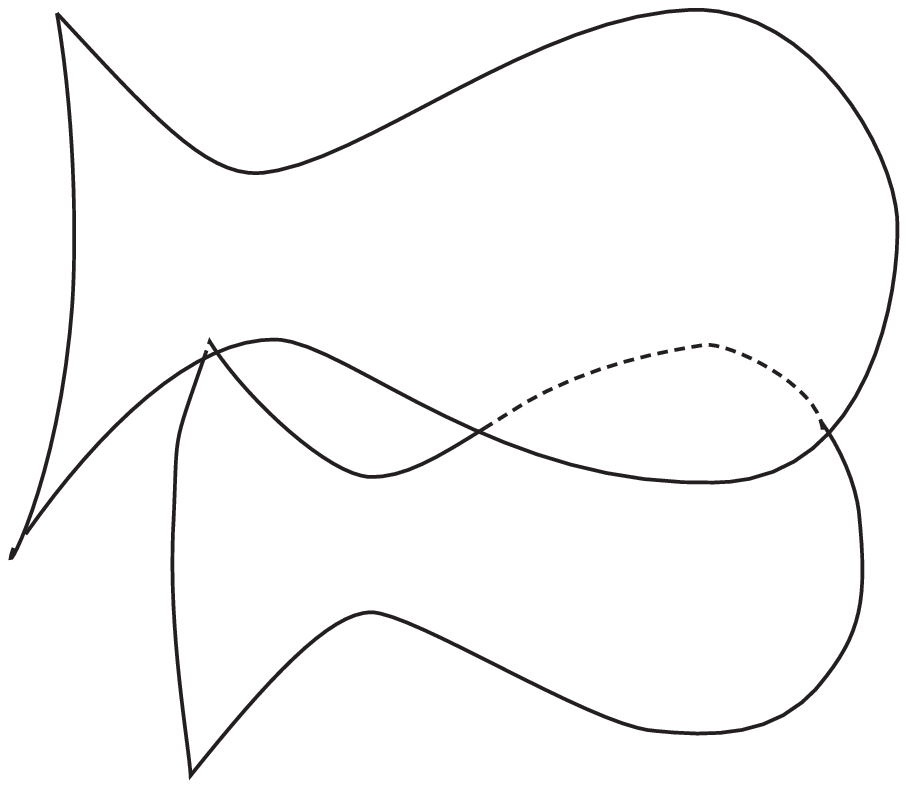} 
\end{center}
 \end{minipage}
 \begin{minipage}{0.30\hsize}
  \begin{center}
    \includegraphics*[width=4cm,height=4cm]{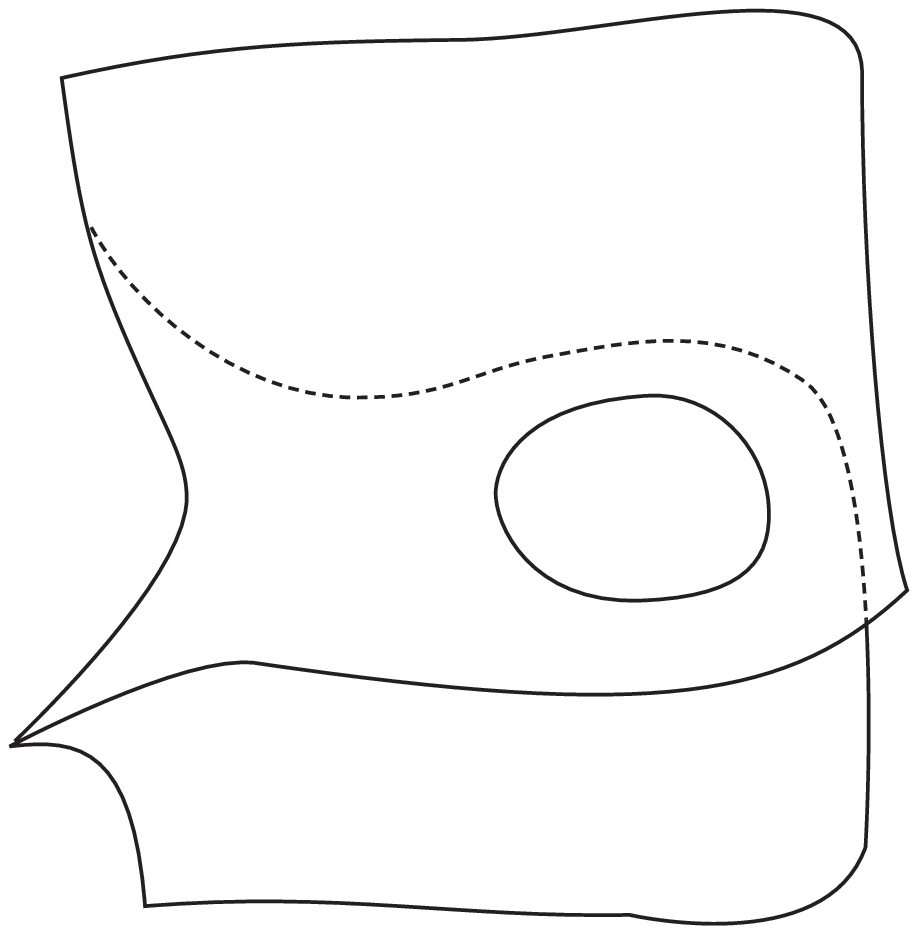}
\end{center}
 \end{minipage}
  \begin{minipage}{0.30\hsize}
  \begin{center}
   \includegraphics*[width=4cm,height=4cm]{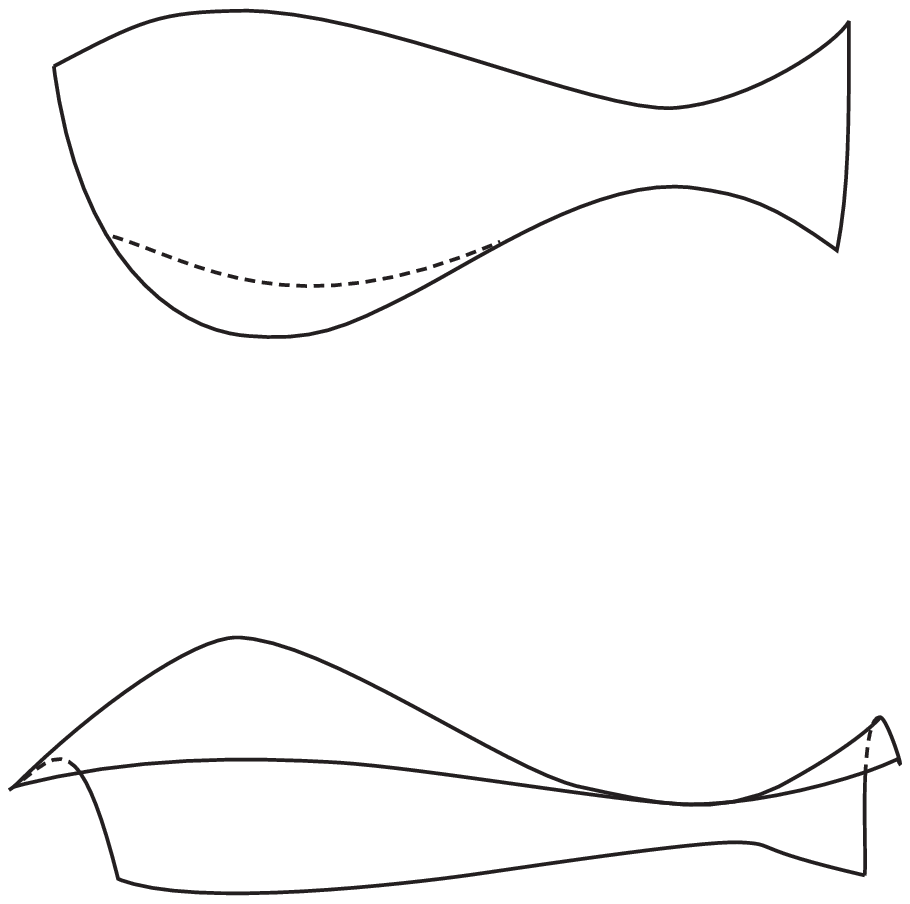}
\end{center}
 \end{minipage}
\end{center}
\caption{${}^2A_1$}
\end{figure}
 \newpage
\begin{figure}[ht]
\begin{center}
 \begin{minipage}{0.30\hsize}
  \begin{center}
    \includegraphics*[width=4cm,height=4cm]{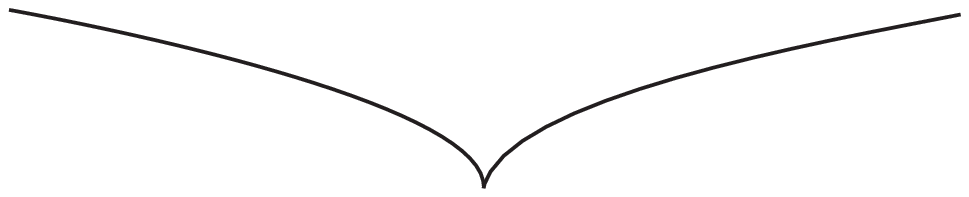} 
\end{center}
 \end{minipage}
 \begin{minipage}{0.30\hsize}
  \begin{center}
    \includegraphics*[width=4cm,height=4cm]{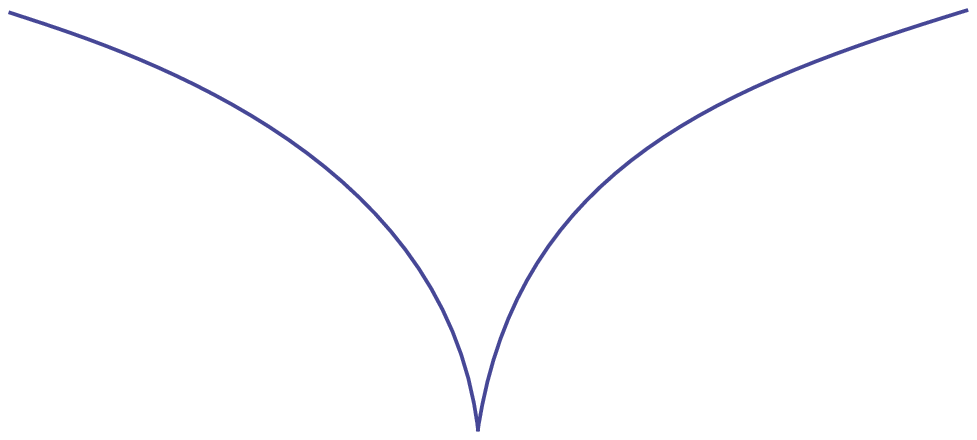}
\end{center}
 \end{minipage}
  \begin{minipage}{0.30\hsize}
  \begin{center}
   \includegraphics*[width=4cm,height=4cm]{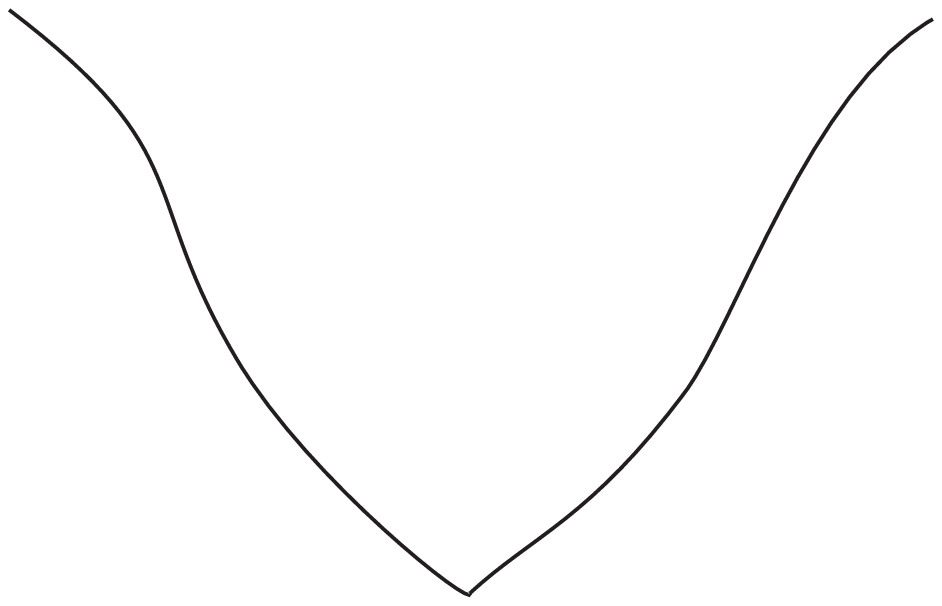}
\end{center}
 \end{minipage}
 \end{center}
 \end{figure}
 \begin{figure}[ht]
\begin{center}
 \begin{minipage}{0.30\hsize}
  \begin{center}
    \includegraphics*[width=4cm,height=4cm]{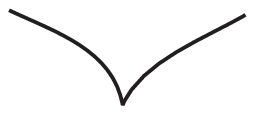} 
\end{center}
 \end{minipage}
 \begin{minipage}{0.30\hsize}
  \begin{center}
    \includegraphics*[width=4cm,height=4cm]{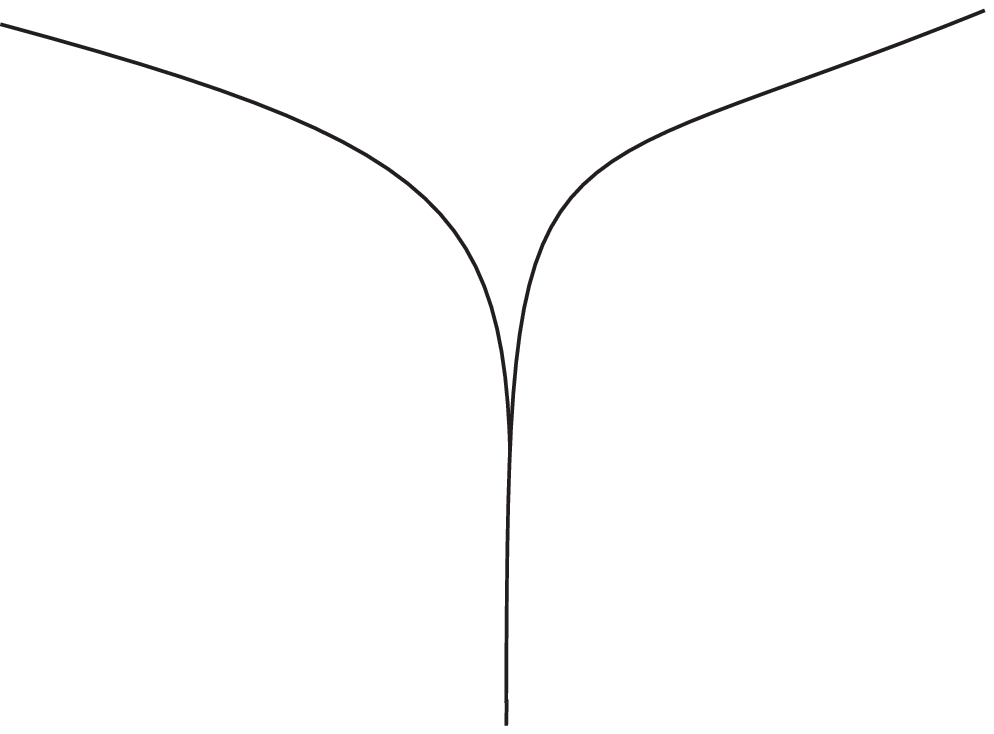}
\end{center}
 \end{minipage}
  \begin{minipage}{0.30\hsize}
  \begin{center}
   \includegraphics*[width=4cm,height=4cm]{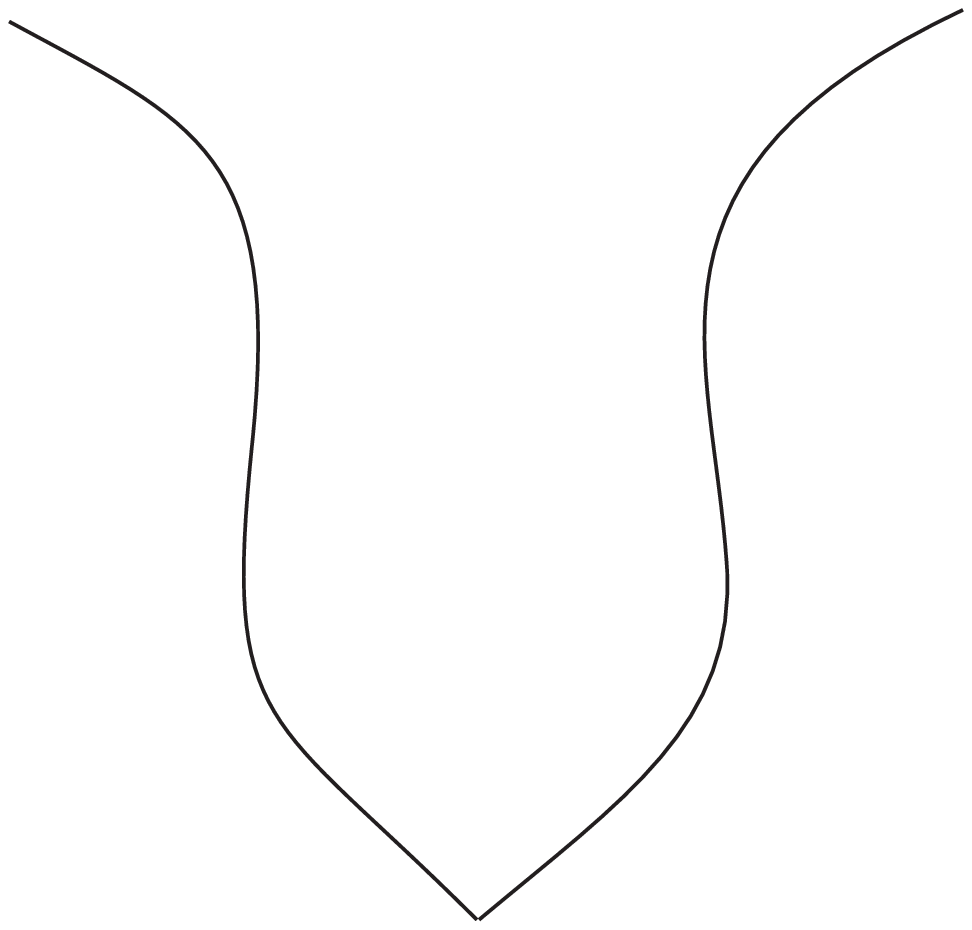}
\end{center}
 \end{minipage}
 \end{center}
 \end{figure}
 \begin{figure}[ht]
\begin{center}
 \begin{minipage}{0.30\hsize}
  \begin{center}
    \includegraphics*[width=4cm,height=4cm]{empty.eps} 
\end{center}
 \end{minipage}
 \begin{minipage}{0.30\hsize}
  \begin{center}
    \includegraphics*[width=4cm,height=4cm]{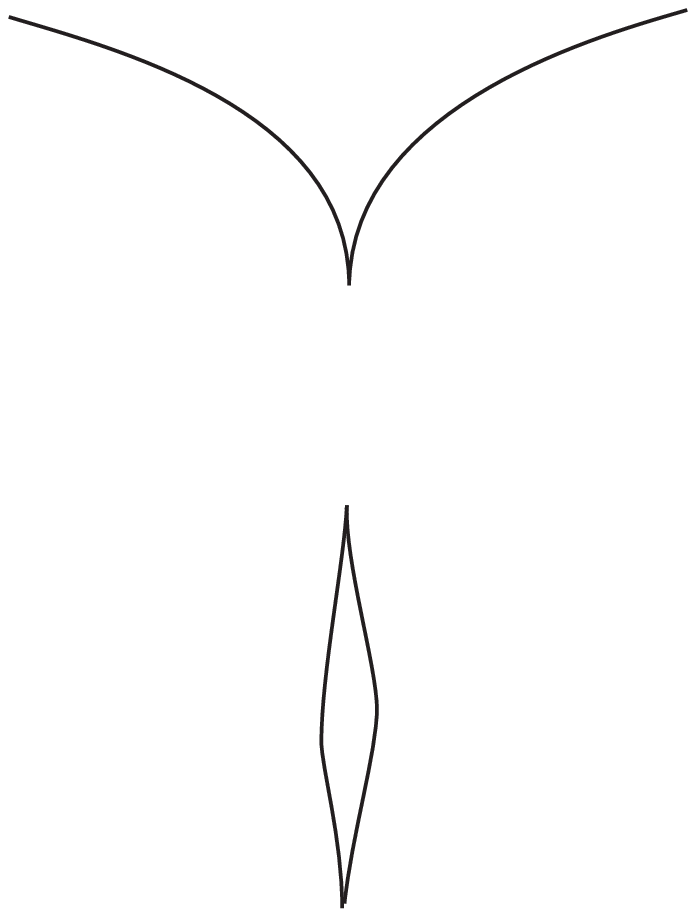}
\end{center}
 \end{minipage}
  \begin{minipage}{0.30\hsize}
  \begin{center}
   \includegraphics*[width=4cm,height=4cm]{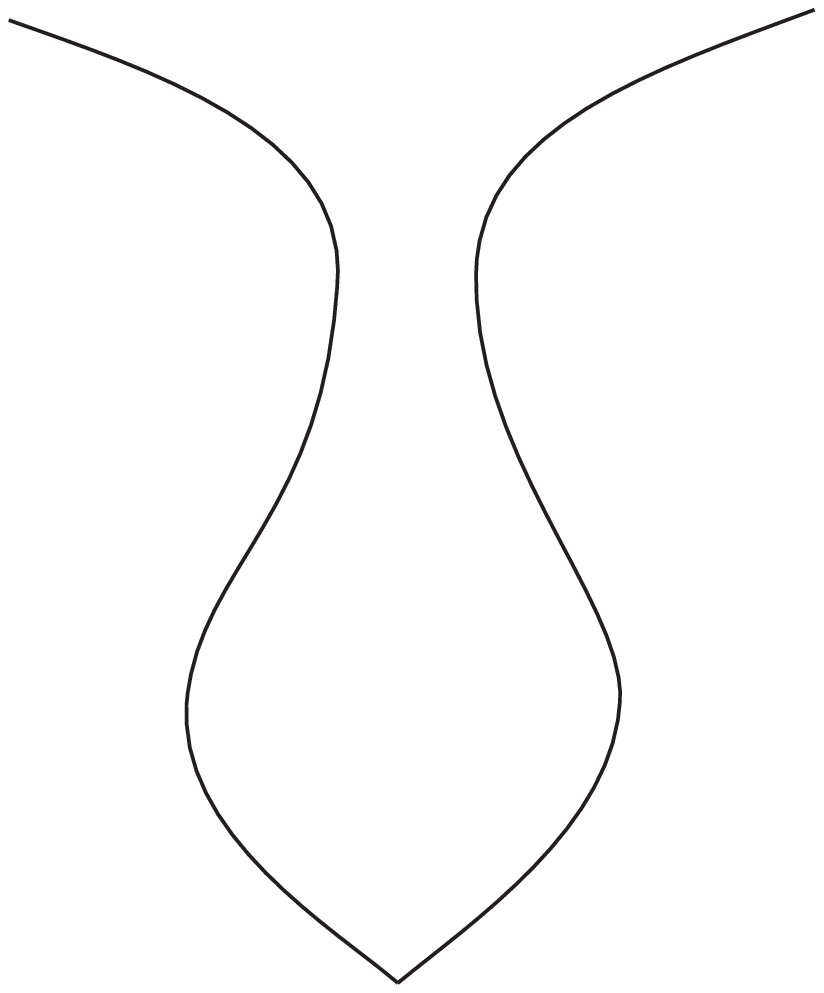}
\end{center}
 \end{minipage}
\end{center}
\caption{${}^2A_2$}
\end{figure}
\newpage
\begin{figure}[ht]
\begin{center}
 \begin{minipage}{0.30\hsize}
  \begin{center}
    \includegraphics*[width=4cm,height=4cm]{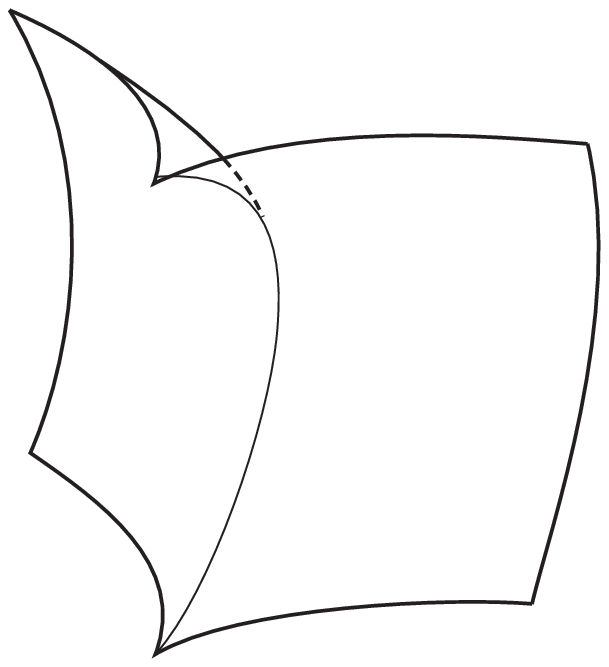} 
\end{center}
 \end{minipage}
 \begin{minipage}{0.30\hsize}
  \begin{center}
    \includegraphics*[width=4cm,height=4cm]{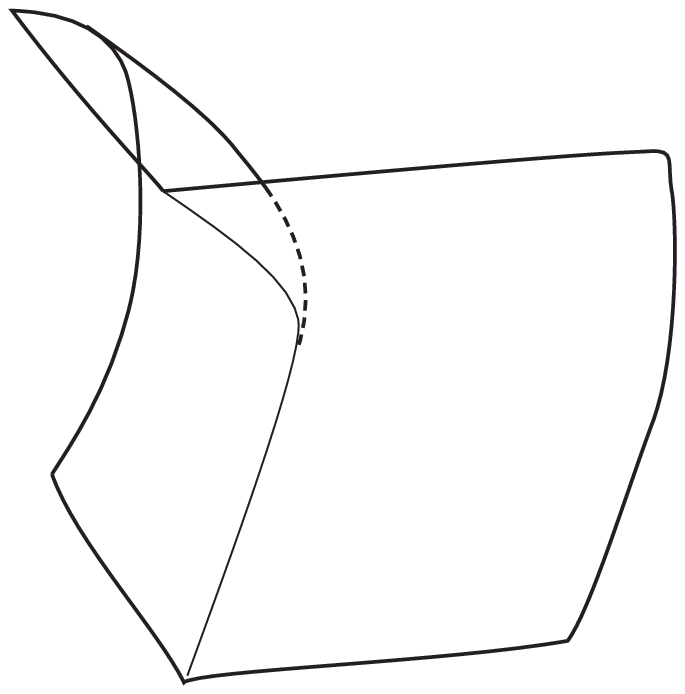}
\end{center}
 \end{minipage}
  \begin{minipage}{0.30\hsize}
  \begin{center}
   \includegraphics*[width=4cm,height=4cm]{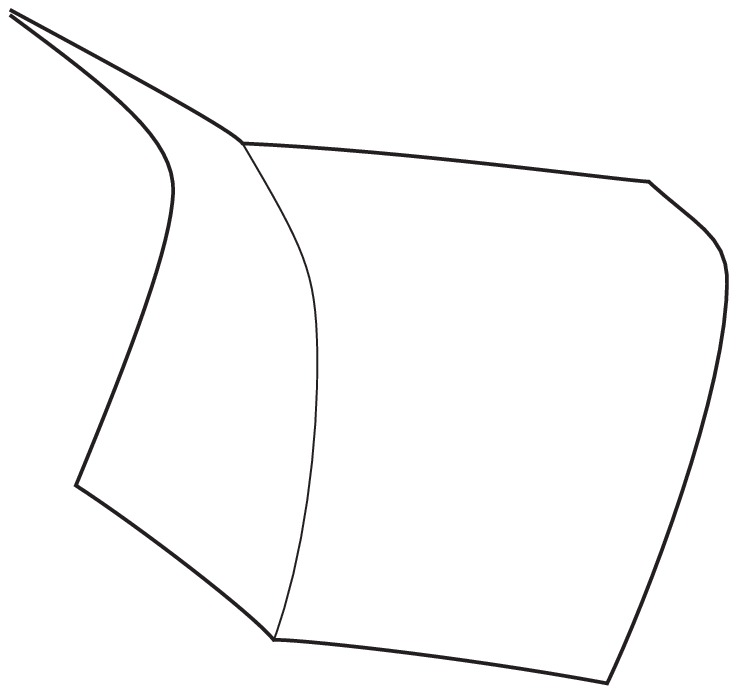}
\end{center}
 \end{minipage}
 \end{center}
 \end{figure}
 \begin{figure}[ht]
\begin{center}
 \begin{minipage}{0.30\hsize}
  \begin{center}
    \includegraphics*[width=4cm,height=4cm]{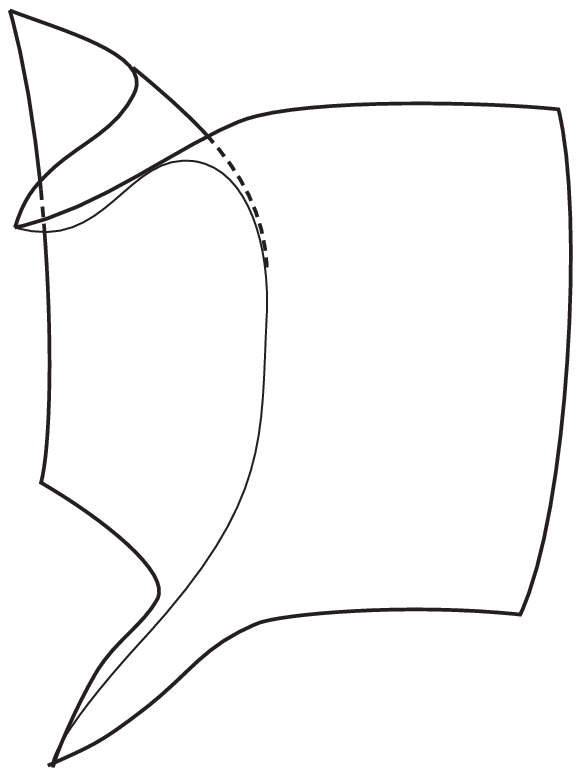} 
\end{center}
 \end{minipage}
 \begin{minipage}{0.30\hsize}
  \begin{center}
    \includegraphics*[width=4cm,height=4cm]{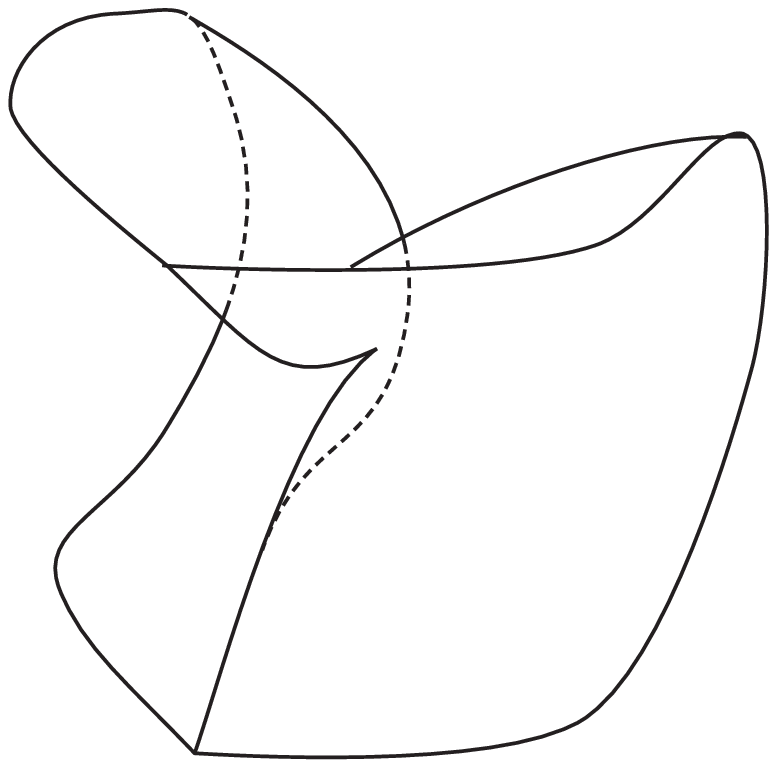}
\end{center}
 \end{minipage}
  \begin{minipage}{0.30\hsize}
  \begin{center}
   \includegraphics*[width=4cm,height=4cm]{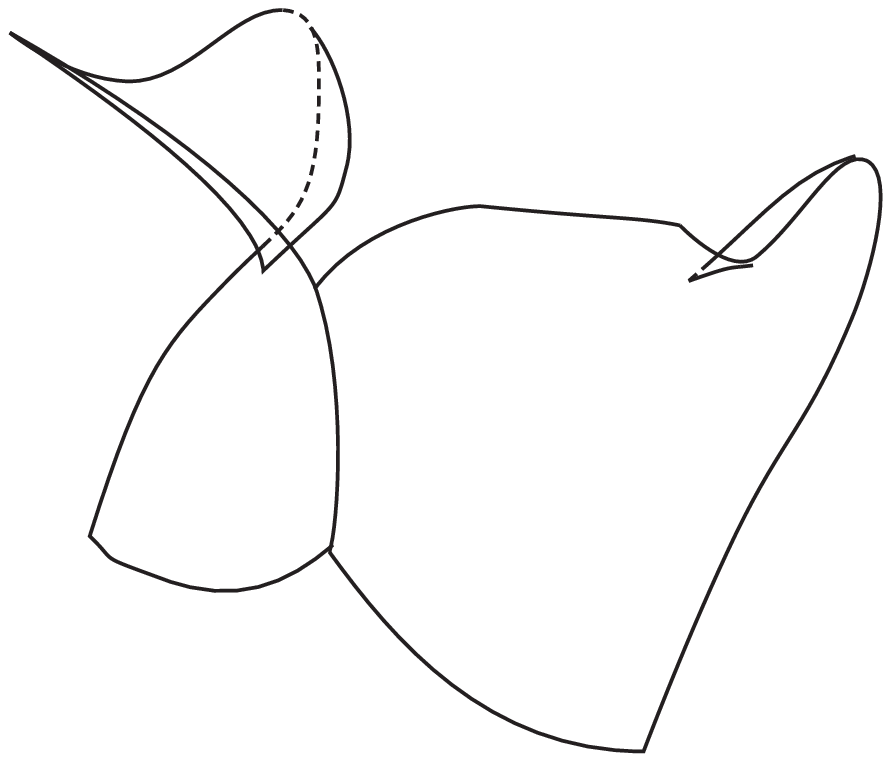}
\end{center}
 \end{minipage}
 \end{center}
 \end{figure}
 \begin{figure}[ht]
\begin{center}
 \begin{minipage}{0.30\hsize}
  \begin{center}
    \includegraphics*[width=4cm,height=4cm]{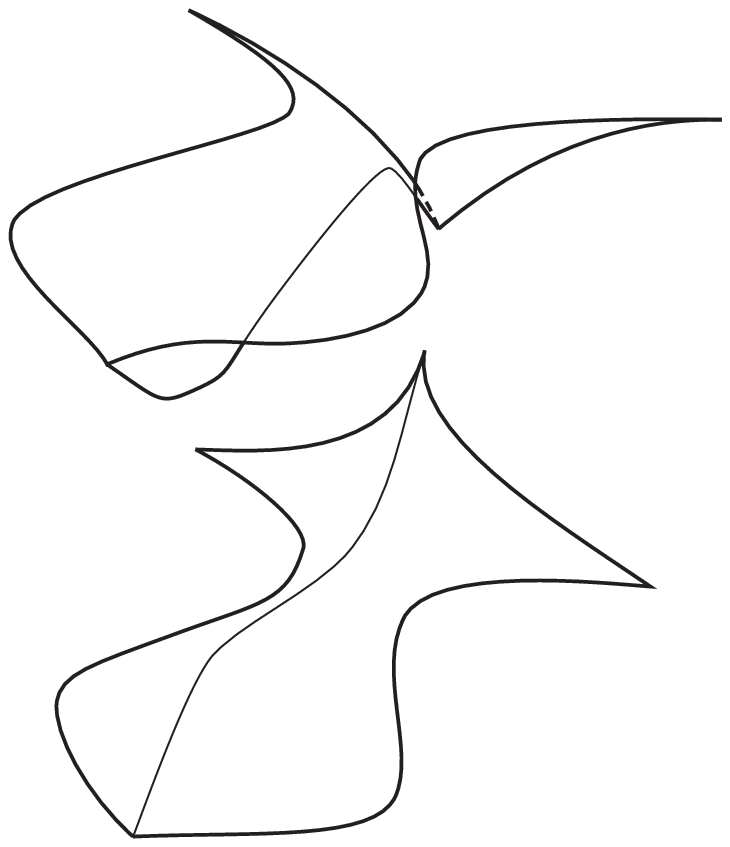} 
\end{center}
 \end{minipage}
 \begin{minipage}{0.30\hsize}
  \begin{center}
    \includegraphics*[width=4cm,height=4cm]{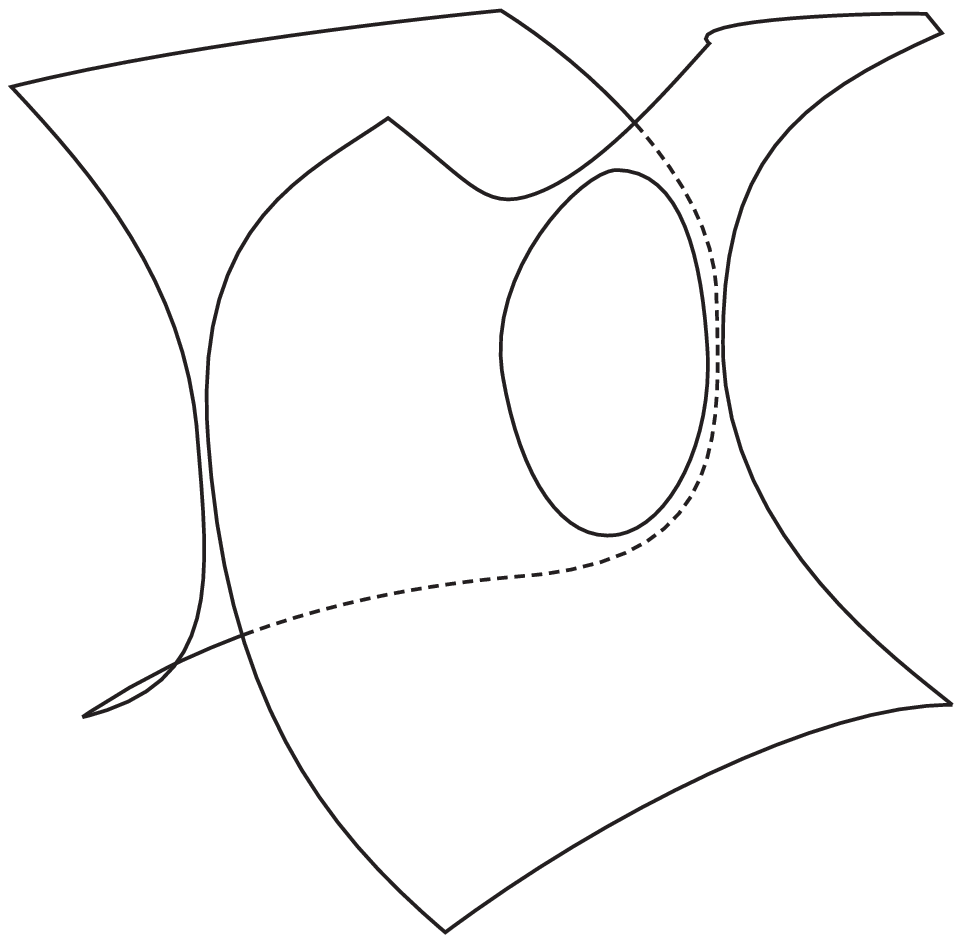}
\end{center}
 \end{minipage}
  \begin{minipage}{0.30\hsize}
  \begin{center}
   \includegraphics*[width=4cm,height=4cm]{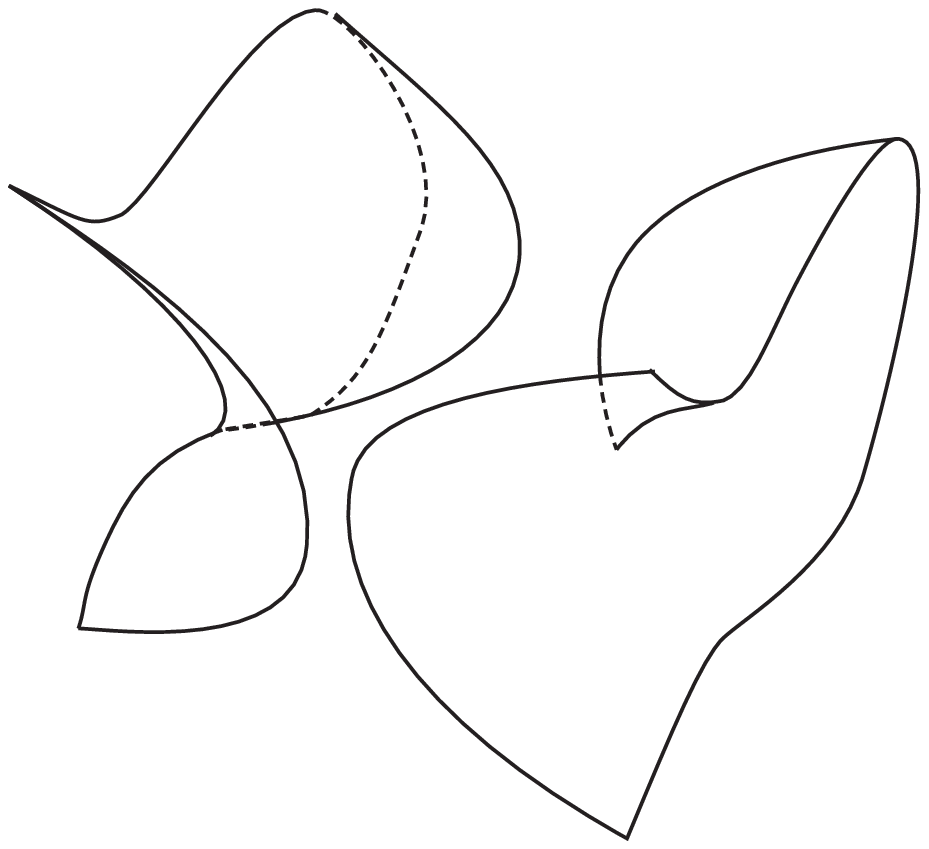}
\end{center}
 \end{minipage}
\end{center}
\caption{${}^2A_2$}
\end{figure}
\newpage
\begin{figure}[ht]
\begin{center}
 \begin{minipage}{0.30\hsize}
  \begin{center}
    \includegraphics*[width=4cm,height=4cm]{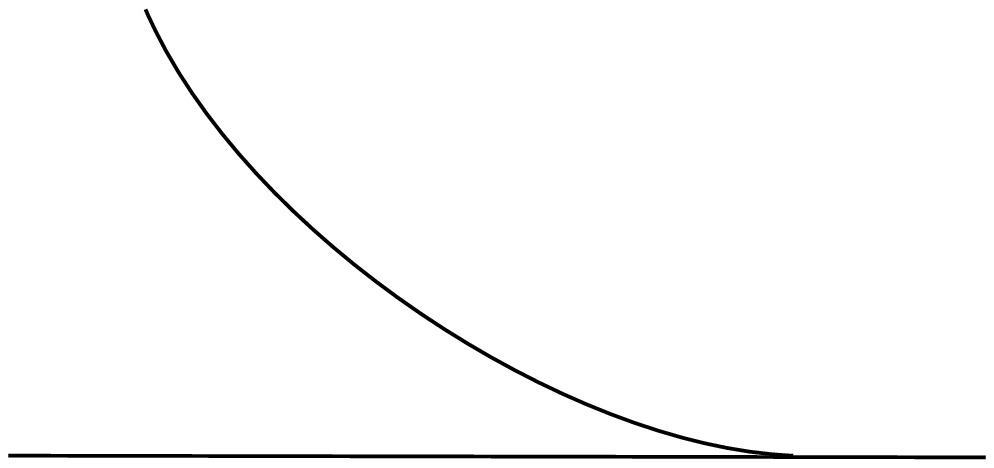} 
\end{center}
 \end{minipage}
 \begin{minipage}{0.30\hsize}
  \begin{center}
    \includegraphics*[width=4cm,height=4cm]{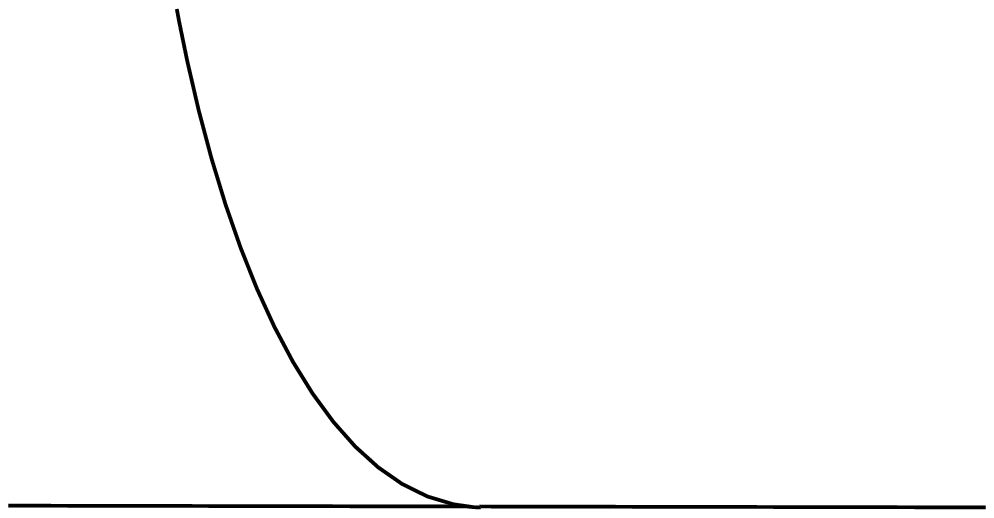}
\end{center}
 \end{minipage}
  \begin{minipage}{0.30\hsize}
  \begin{center}
   \includegraphics*[width=4cm,height=4cm]{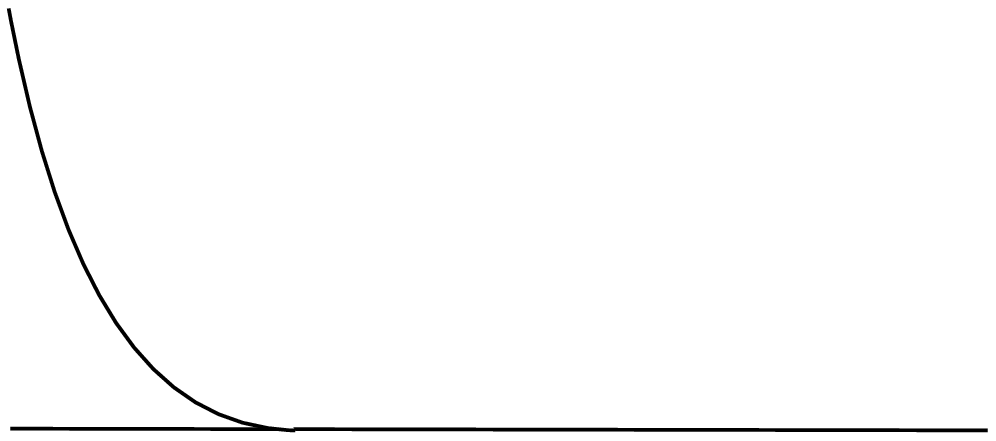}
\end{center}
 \end{minipage}
 \end{center}
 \end{figure}
 \begin{figure}[ht]
\begin{center}
 \begin{minipage}{0.30\hsize}
  \begin{center}
    \includegraphics*[width=4cm,height=4cm]{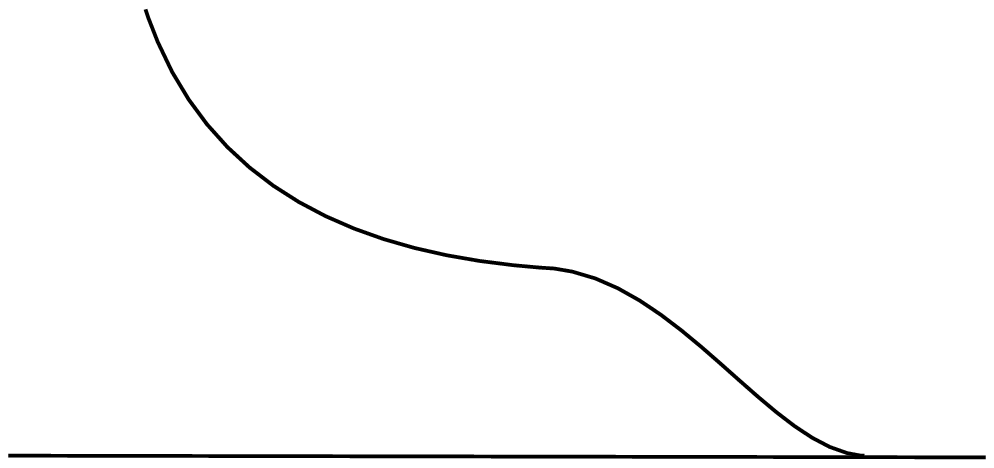} 
\end{center}
 \end{minipage}
 \begin{minipage}{0.30\hsize}
  \begin{center}
    \includegraphics*[width=4cm,height=4cm]{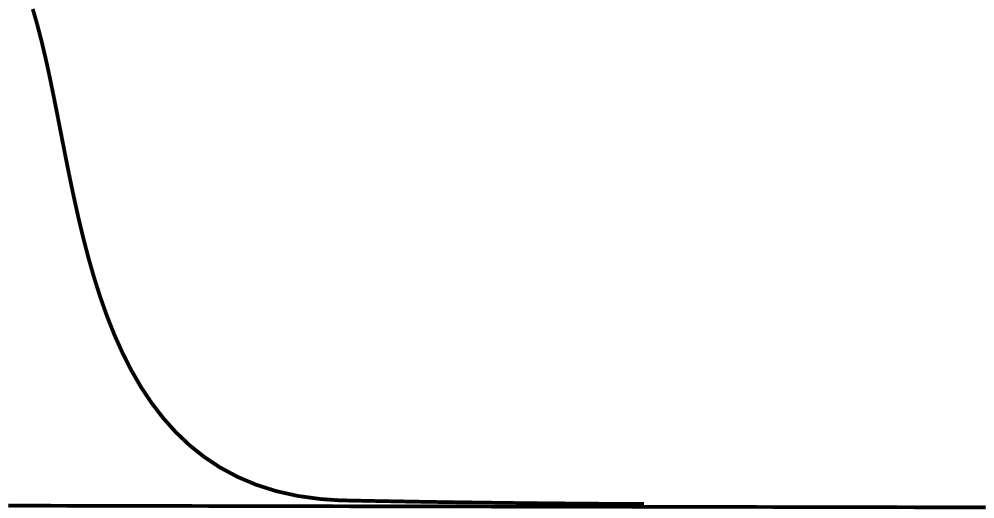}
\end{center}
 \end{minipage}
  \begin{minipage}{0.30\hsize}
  \begin{center}
   \includegraphics*[width=4cm,height=4cm]{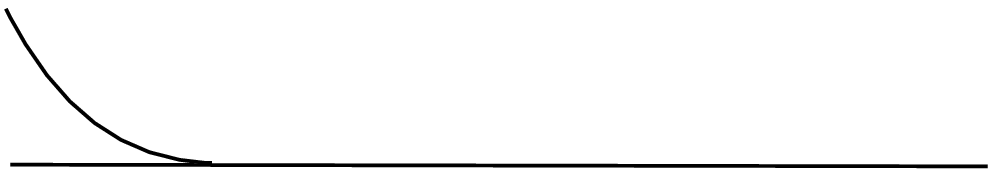}
\end{center}
 \end{minipage}
 \end{center}
 \end{figure}
 \begin{figure}[ht]
\begin{center}
 \begin{minipage}{0.30\hsize}
  \begin{center}
    \includegraphics*[width=4cm,height=4cm]{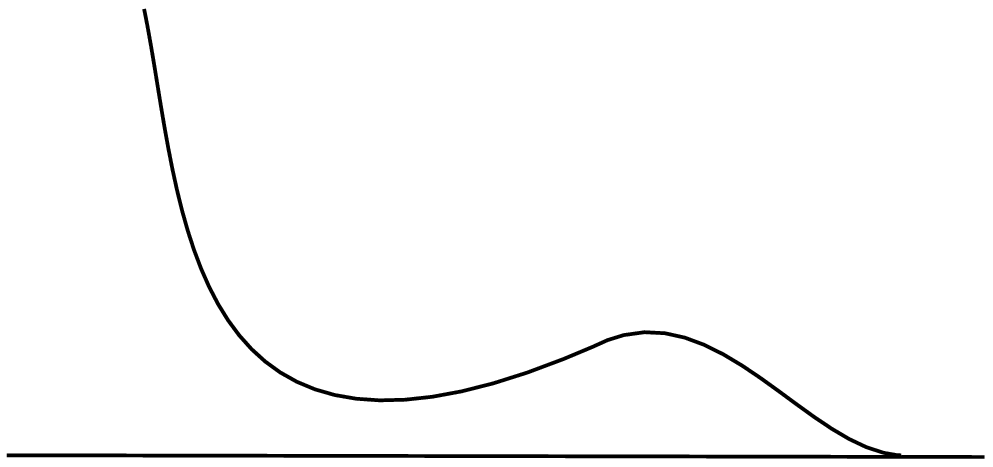} 
\end{center}
 \end{minipage}
 \begin{minipage}{0.30\hsize}
  \begin{center}
    \includegraphics*[width=4cm,height=4cm]{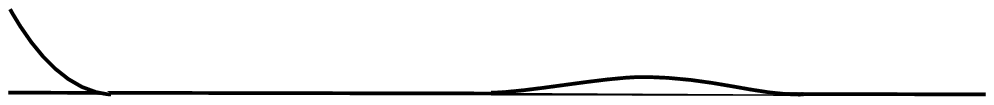}
\end{center}
 \end{minipage}
  \begin{minipage}{0.30\hsize}
  \begin{center}
   \includegraphics*[width=4cm,height=4cm]{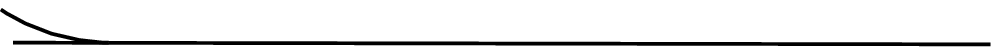}
\end{center}
 \end{minipage}
\end{center}
\caption{${}^2B_2$}
\end{figure}
\newpage
\begin{figure}[ht]
\begin{center}
 \begin{minipage}{0.30\hsize}
  \begin{center}
    \includegraphics*[width=4cm,height=4cm]{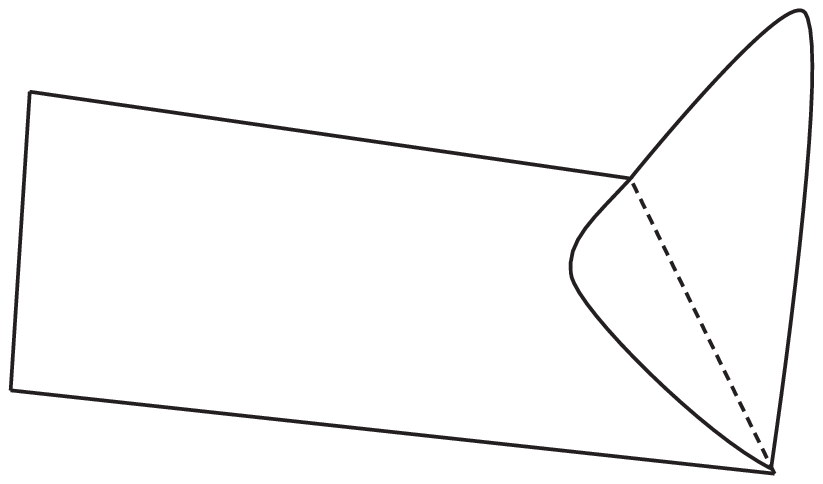} 
\end{center}
 \end{minipage}
 \begin{minipage}{0.30\hsize}
  \begin{center}
    \includegraphics*[width=4cm,height=4cm]{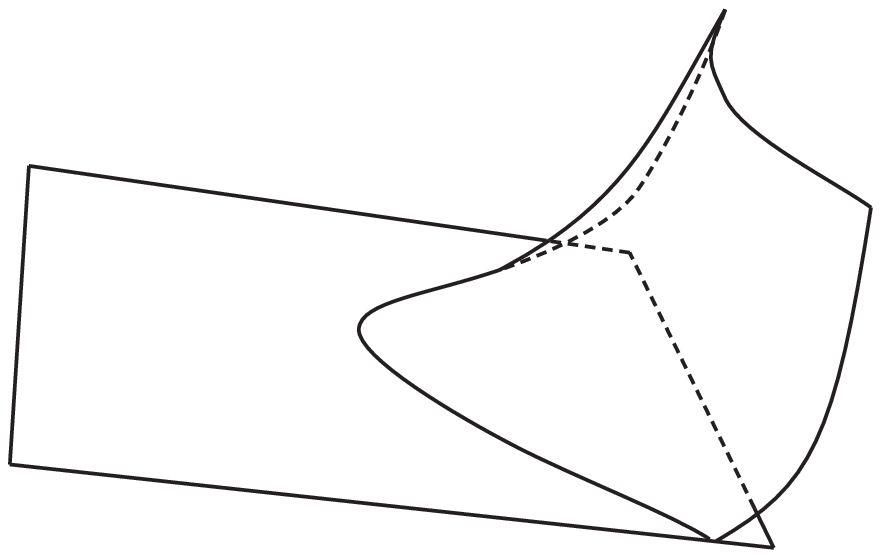}
\end{center}
 \end{minipage}
  \begin{minipage}{0.30\hsize}
  \begin{center}
   \includegraphics*[width=4cm,height=4cm]{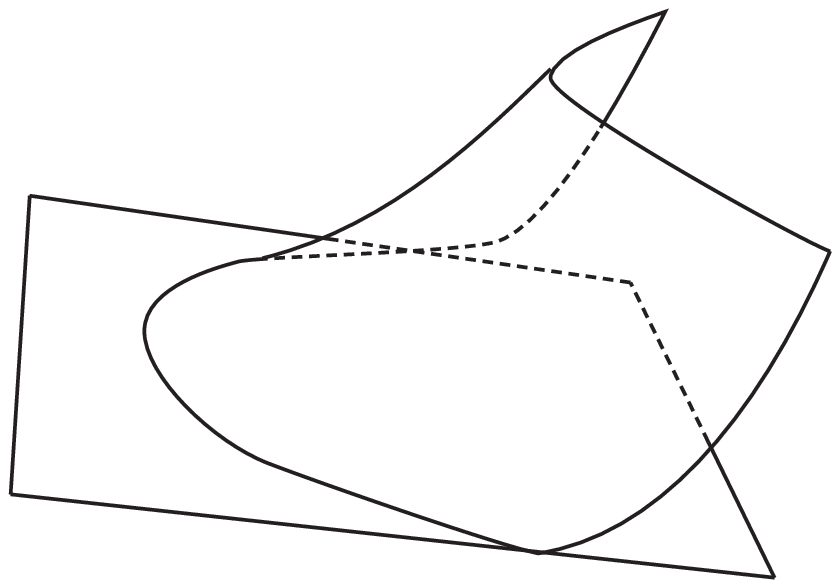}
\end{center}
 \end{minipage}
 \end{center}
 \end{figure}
 \begin{figure}[ht]
\begin{center}
 \begin{minipage}{0.30\hsize}
  \begin{center}
    \includegraphics*[width=4cm,height=4cm]{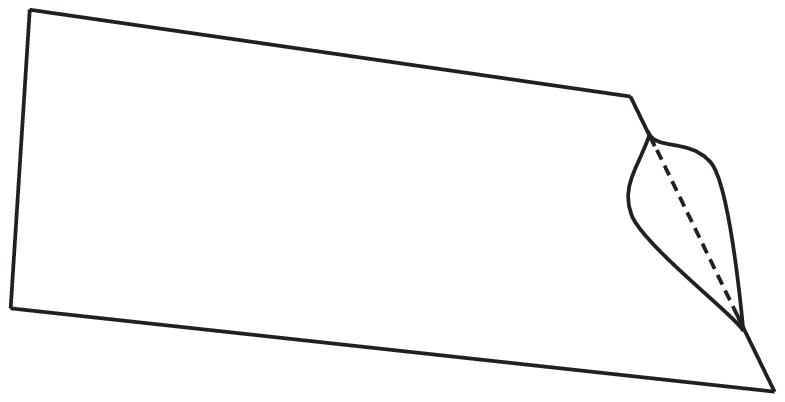} 
\end{center}
 \end{minipage}
 \begin{minipage}{0.30\hsize}
  \begin{center}
    \includegraphics*[width=4cm,height=4cm]{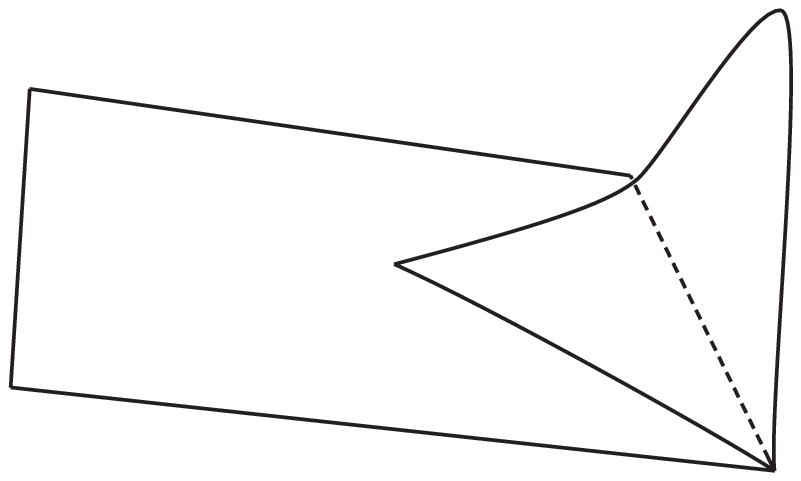}
\end{center}
 \end{minipage}
  \begin{minipage}{0.30\hsize}
  \begin{center}
   \includegraphics*[width=4cm,height=4cm]{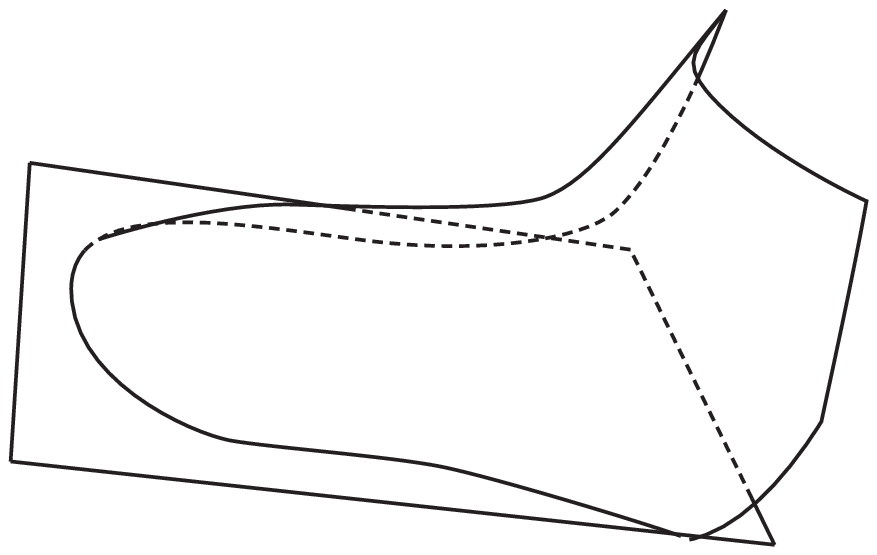}
\end{center}
 \end{minipage}
 \end{center}
 \end{figure}
 \begin{figure}[ht]
\begin{center}
 \begin{minipage}{0.30\hsize}
  \begin{center}
    \includegraphics*[width=4cm,height=4cm]{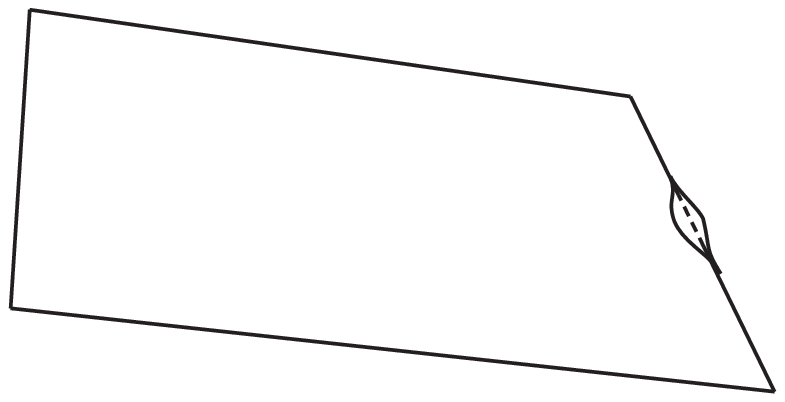} 
\end{center}
 \end{minipage}
 \begin{minipage}{0.30\hsize}
  \begin{center}
    \includegraphics*[width=4cm,height=4cm]{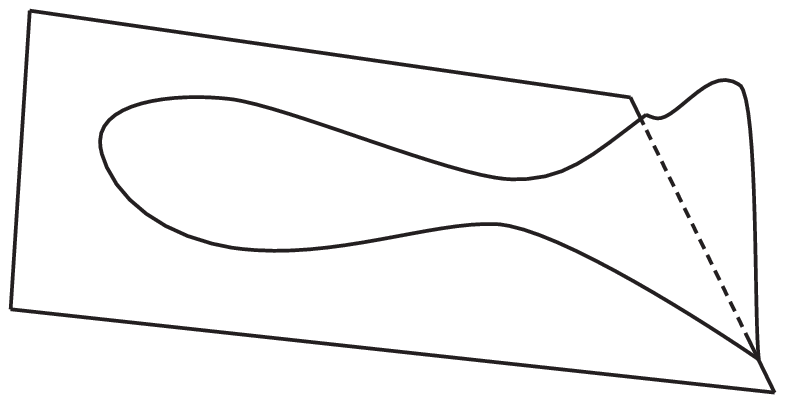}
\end{center}
 \end{minipage}
  \begin{minipage}{0.30\hsize}
  \begin{center}
   \includegraphics*[width=4cm,height=4cm]{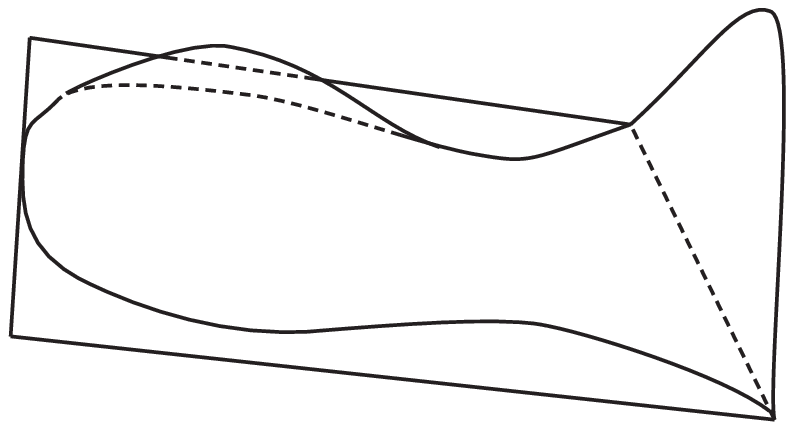}
\end{center}
 \end{minipage}
\end{center}
\caption{${}^2B_2$}
\end{figure}

\newpage
\begin{figure}[ht]
\begin{center}
 \begin{minipage}{0.30\hsize}
  \begin{center}
    \includegraphics*[width=4cm,height=4cm]{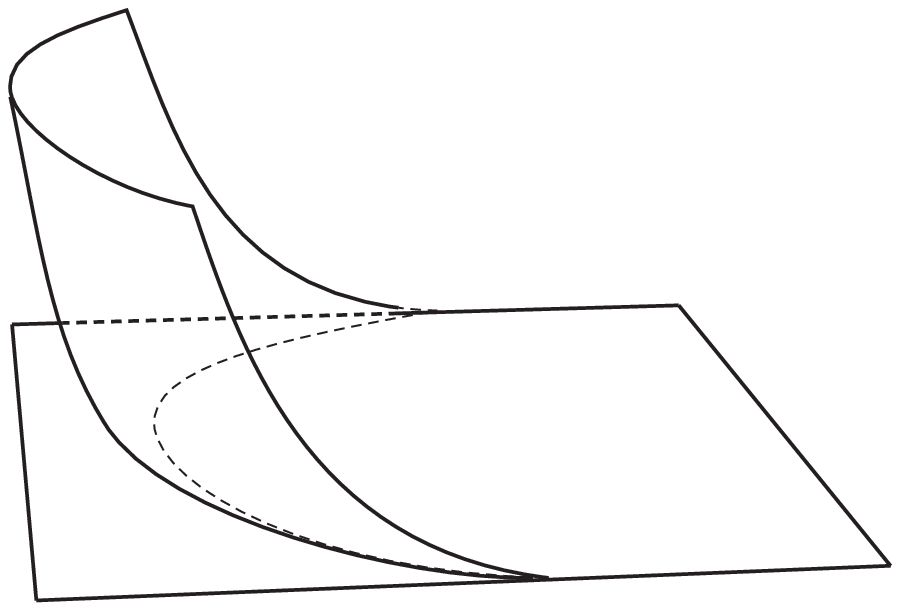} 
\end{center}
 \end{minipage}
 \begin{minipage}{0.30\hsize}
  \begin{center}
    \includegraphics*[width=4cm,height=4cm]{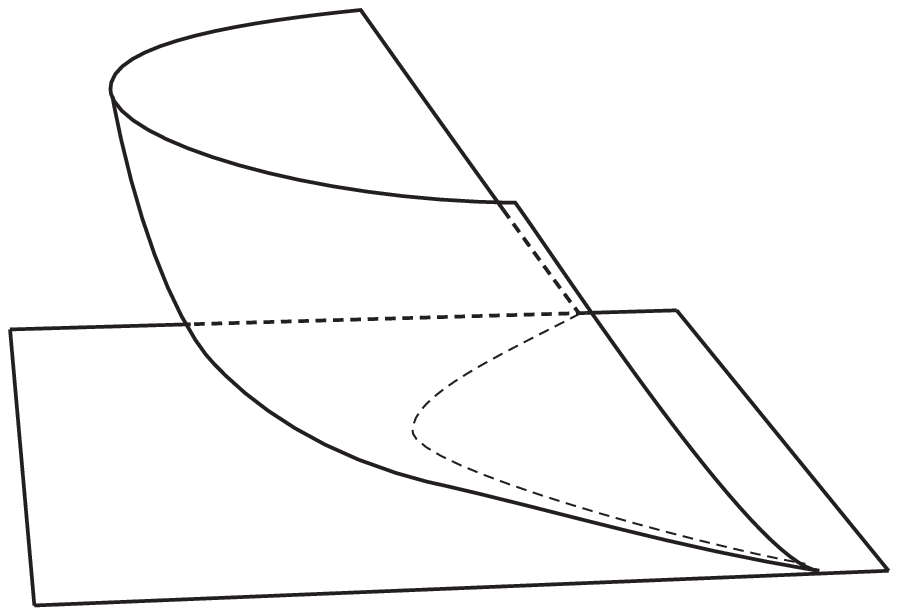}
\end{center}
 \end{minipage}
  \begin{minipage}{0.30\hsize}
  \begin{center}
   \includegraphics*[width=4cm,height=4cm]{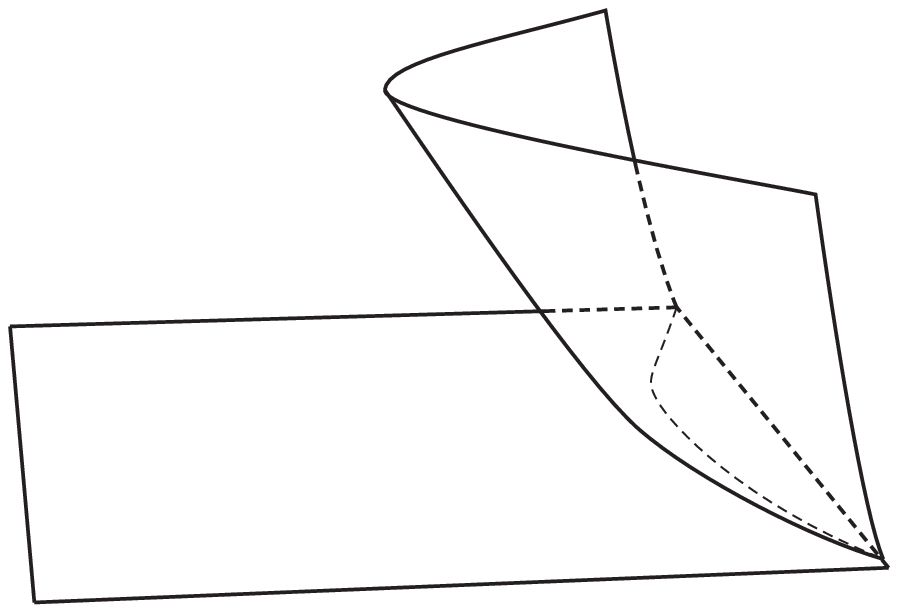}
\end{center}
 \end{minipage}
 \end{center}
 \end{figure}
 \begin{figure}[ht]
\begin{center}
 \begin{minipage}{0.30\hsize}
  \begin{center}
    \includegraphics*[width=4cm,height=4cm]{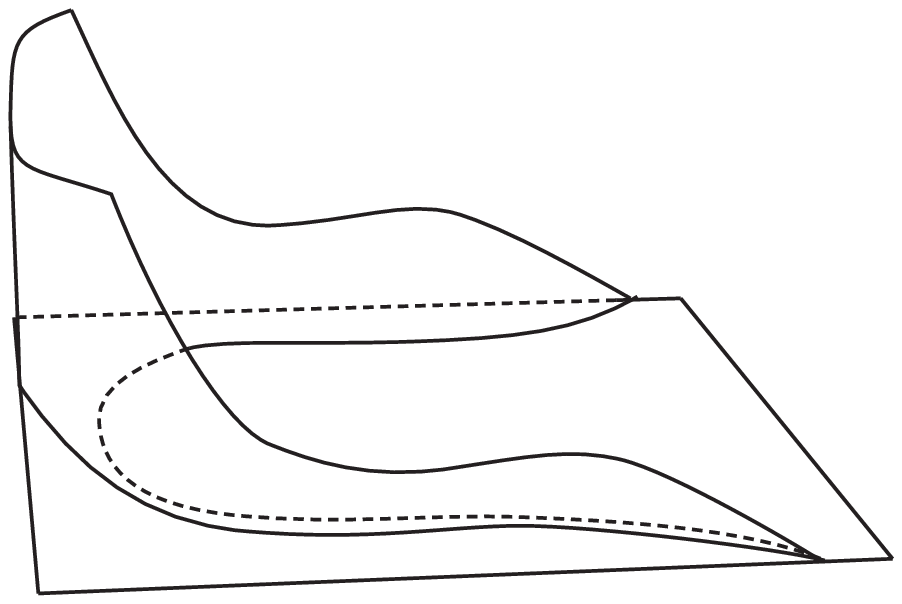} 
\end{center}
 \end{minipage}
 \begin{minipage}{0.30\hsize}
  \begin{center}
    \includegraphics*[width=4cm,height=4cm]{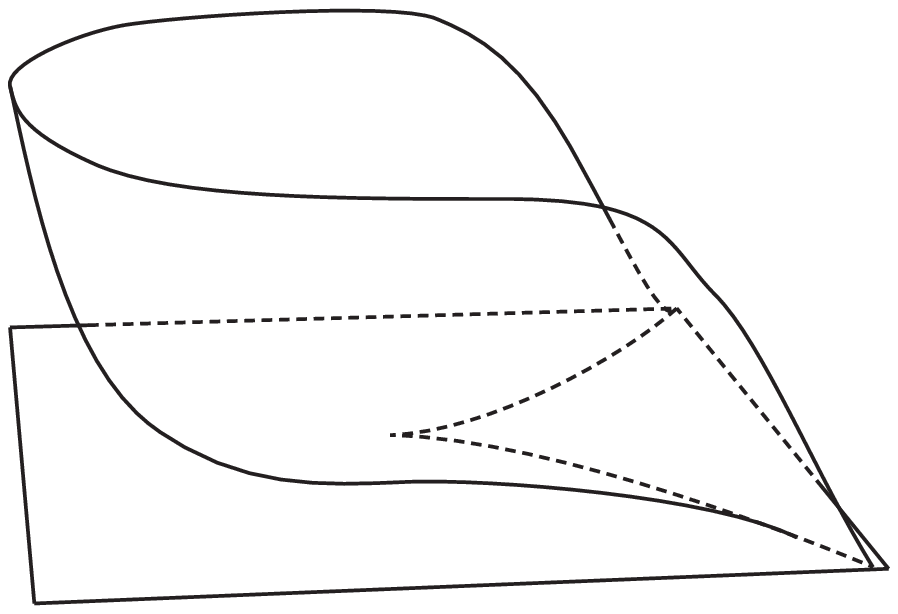}
\end{center}
 \end{minipage}
  \begin{minipage}{0.30\hsize}
  \begin{center}
   \includegraphics*[width=4cm,height=4cm]{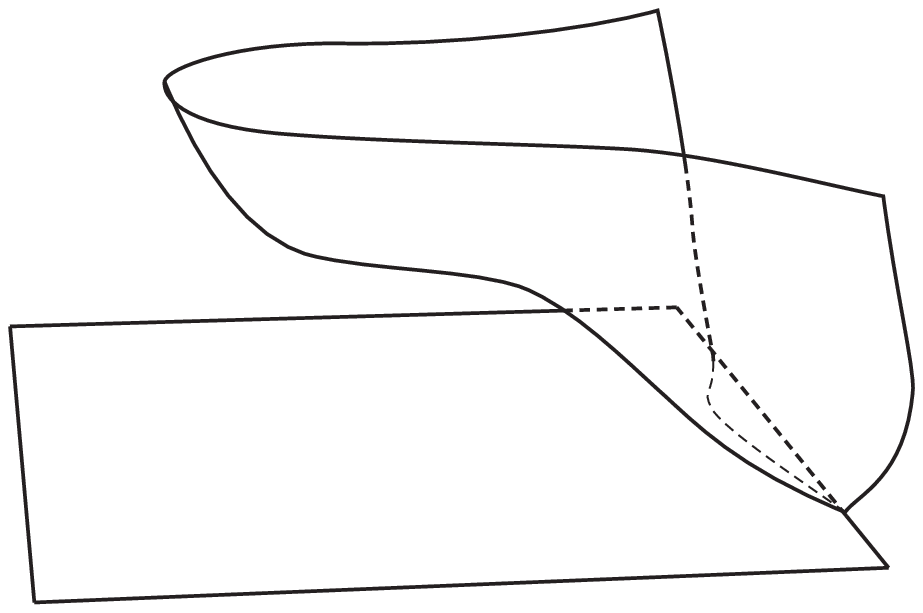}
\end{center}
 \end{minipage}
 \end{center}
 \end{figure}
 \begin{figure}[ht]
\begin{center}
 \begin{minipage}{0.30\hsize}
  \begin{center}
    \includegraphics*[width=4cm,height=4cm]{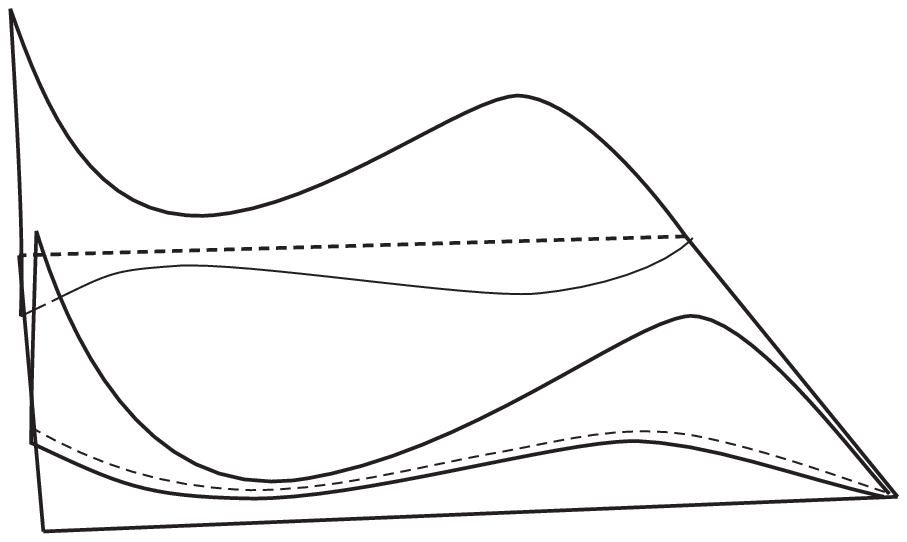} 
\end{center}
 \end{minipage}
 \begin{minipage}{0.30\hsize}
  \begin{center}
    \includegraphics*[width=4cm,height=4cm]{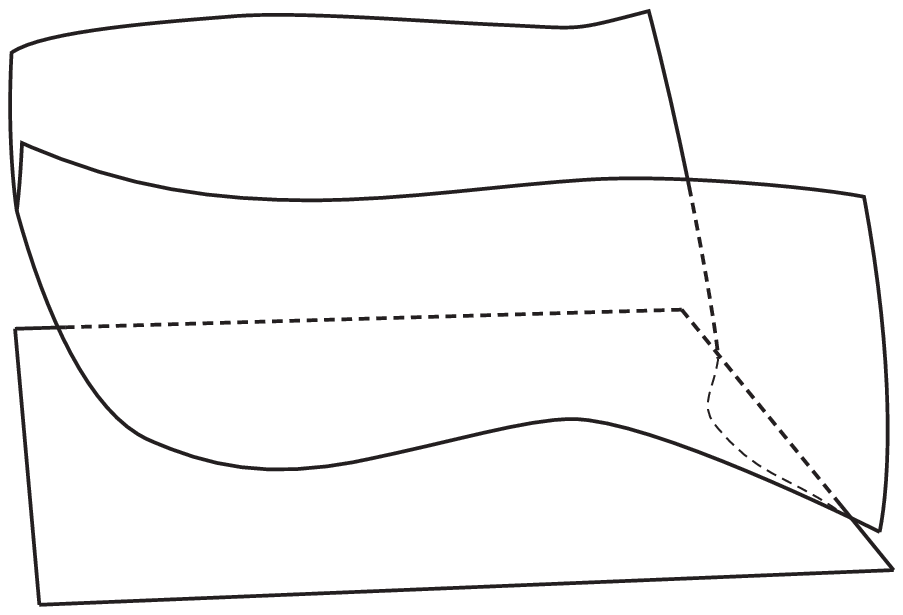}
\end{center}
 \end{minipage}
  \begin{minipage}{0.30\hsize}
  \begin{center}
   \includegraphics*[width=4cm,height=4cm]{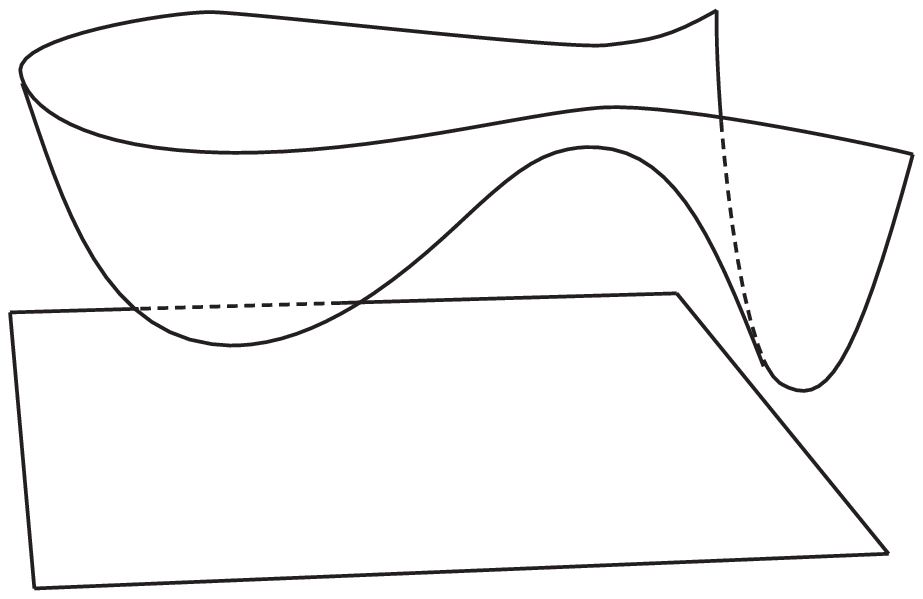}
\end{center}
 \end{minipage}
\end{center}
\caption{${}^2B_2$}
\end{figure}
\newpage
\begin{figure}[ht]
\begin{center}
 \begin{minipage}{0.30\hsize}
  \begin{center}
    \includegraphics*[width=4cm,height=4cm]{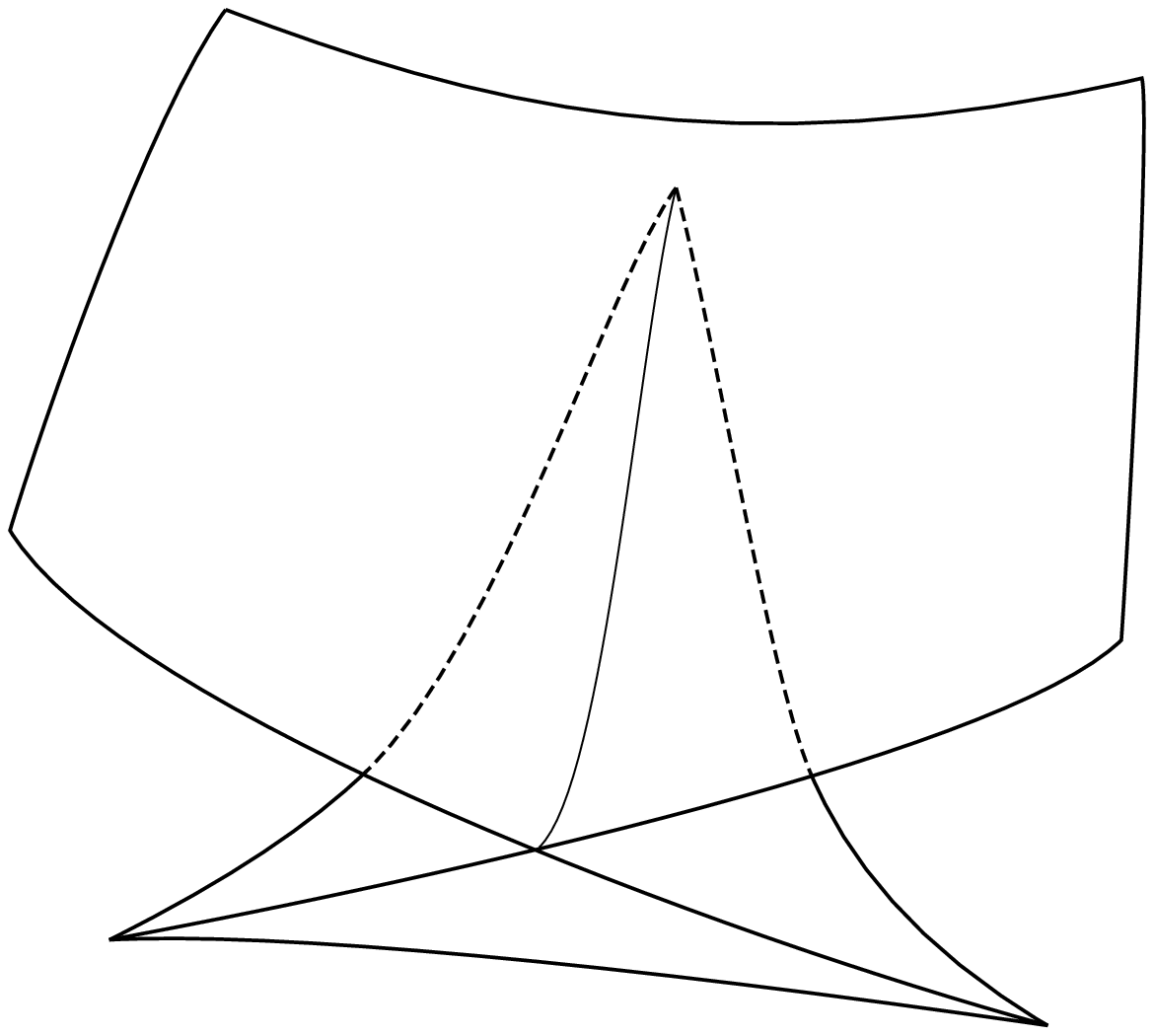} 
\end{center}
 \end{minipage}
 \begin{minipage}{0.30\hsize}
  \begin{center}
    \includegraphics*[width=4cm,height=4cm]{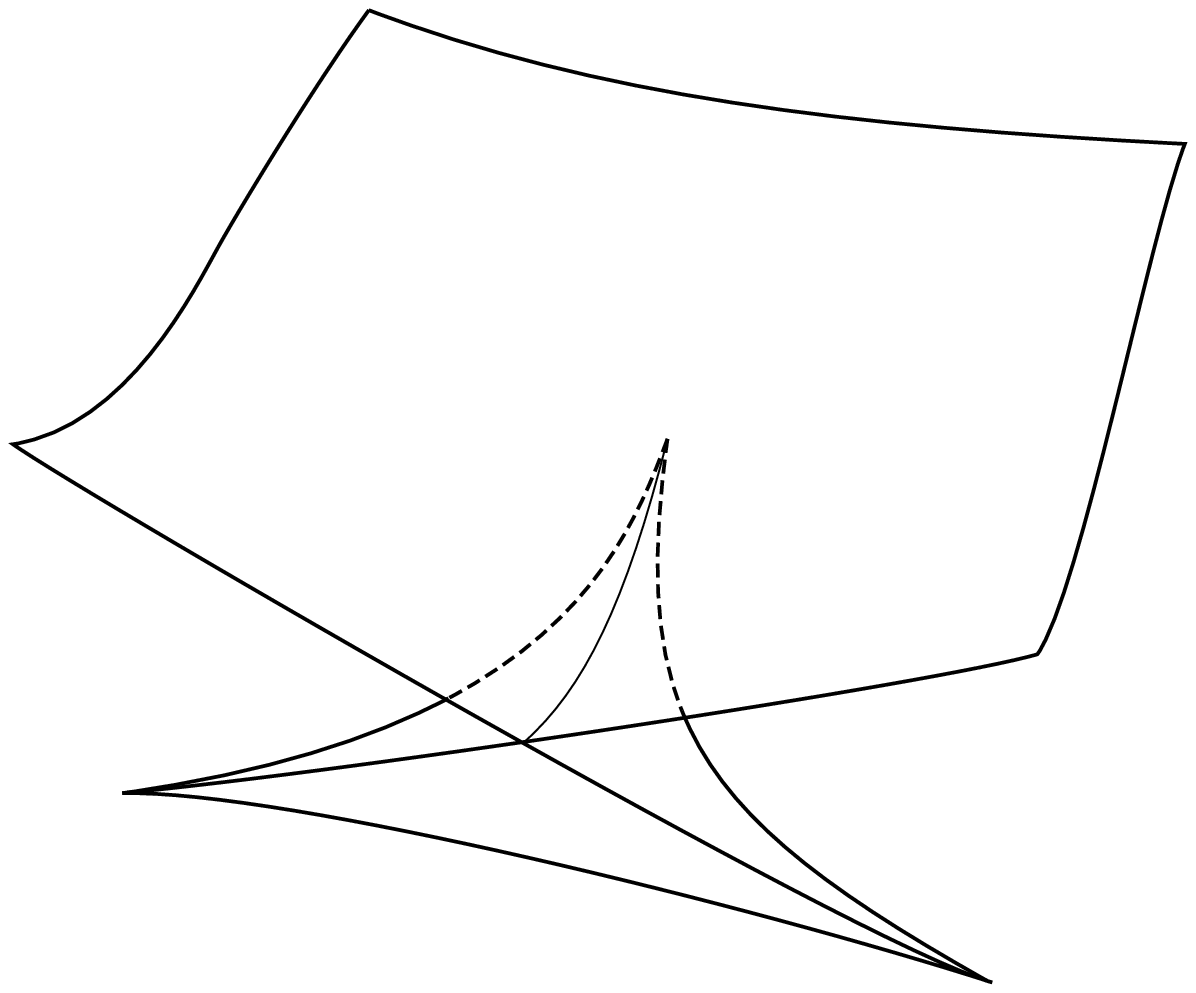}
\end{center}
 \end{minipage}
  \begin{minipage}{0.30\hsize}
  \begin{center}
   \includegraphics*[width=4cm,height=4cm]{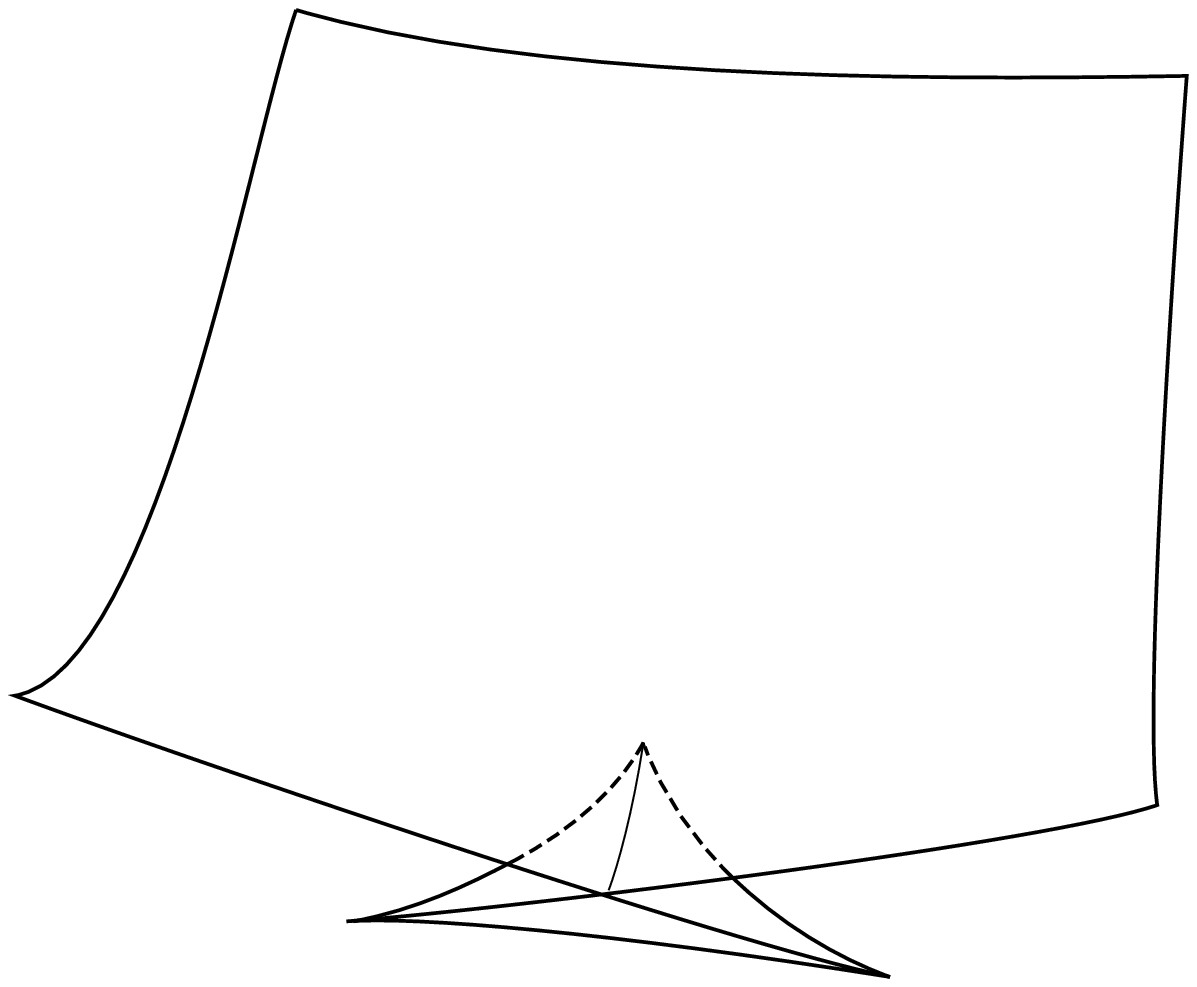}
\end{center}
 \end{minipage}
 \end{center}
 \end{figure}
 \begin{figure}[ht]
\begin{center}
 \begin{minipage}{0.30\hsize}
  \begin{center}
    \includegraphics*[width=4cm,height=4cm]{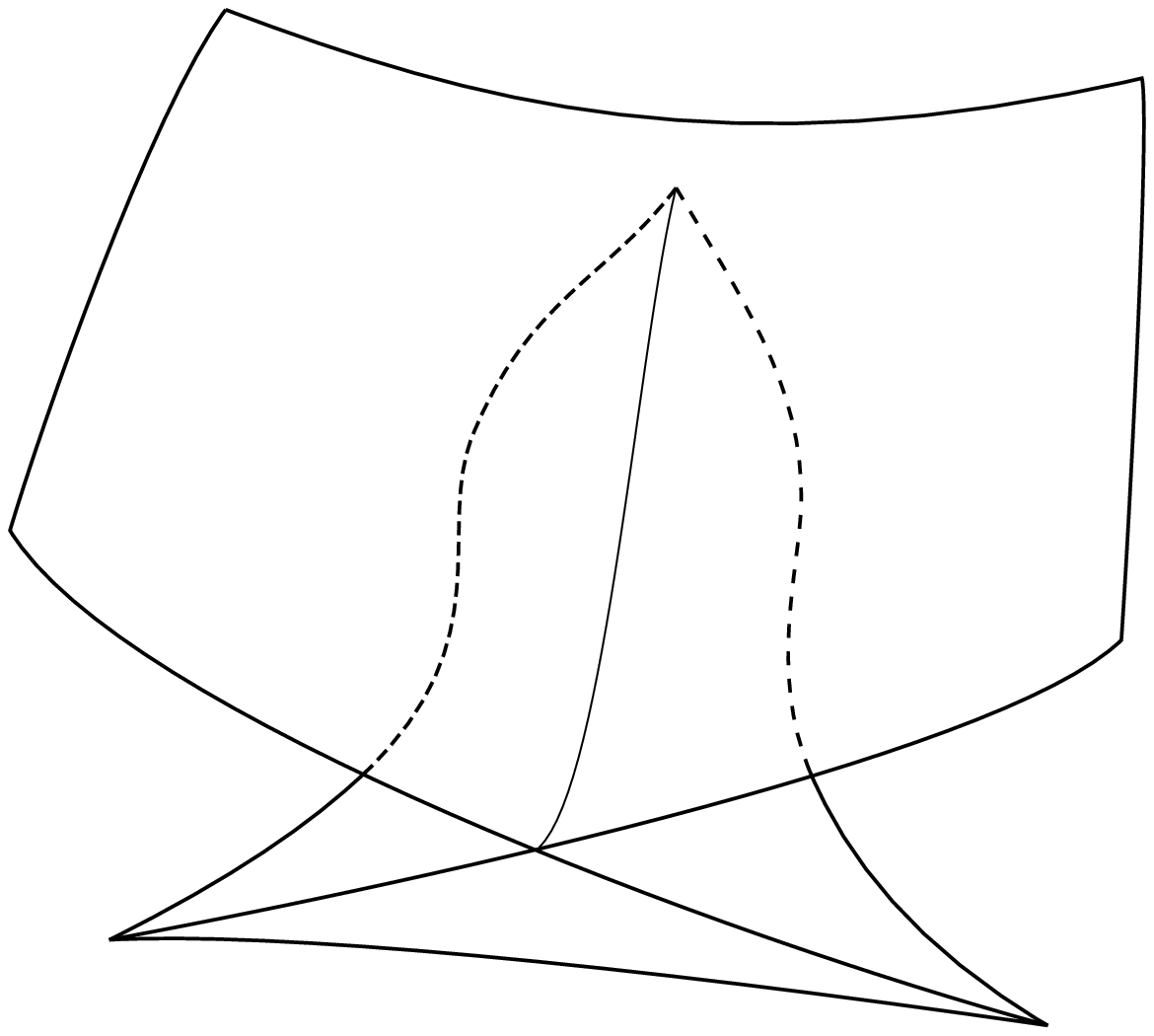} 
\end{center}
 \end{minipage}
 \begin{minipage}{0.30\hsize}
  \begin{center}
    \includegraphics*[width=4cm,height=4cm]{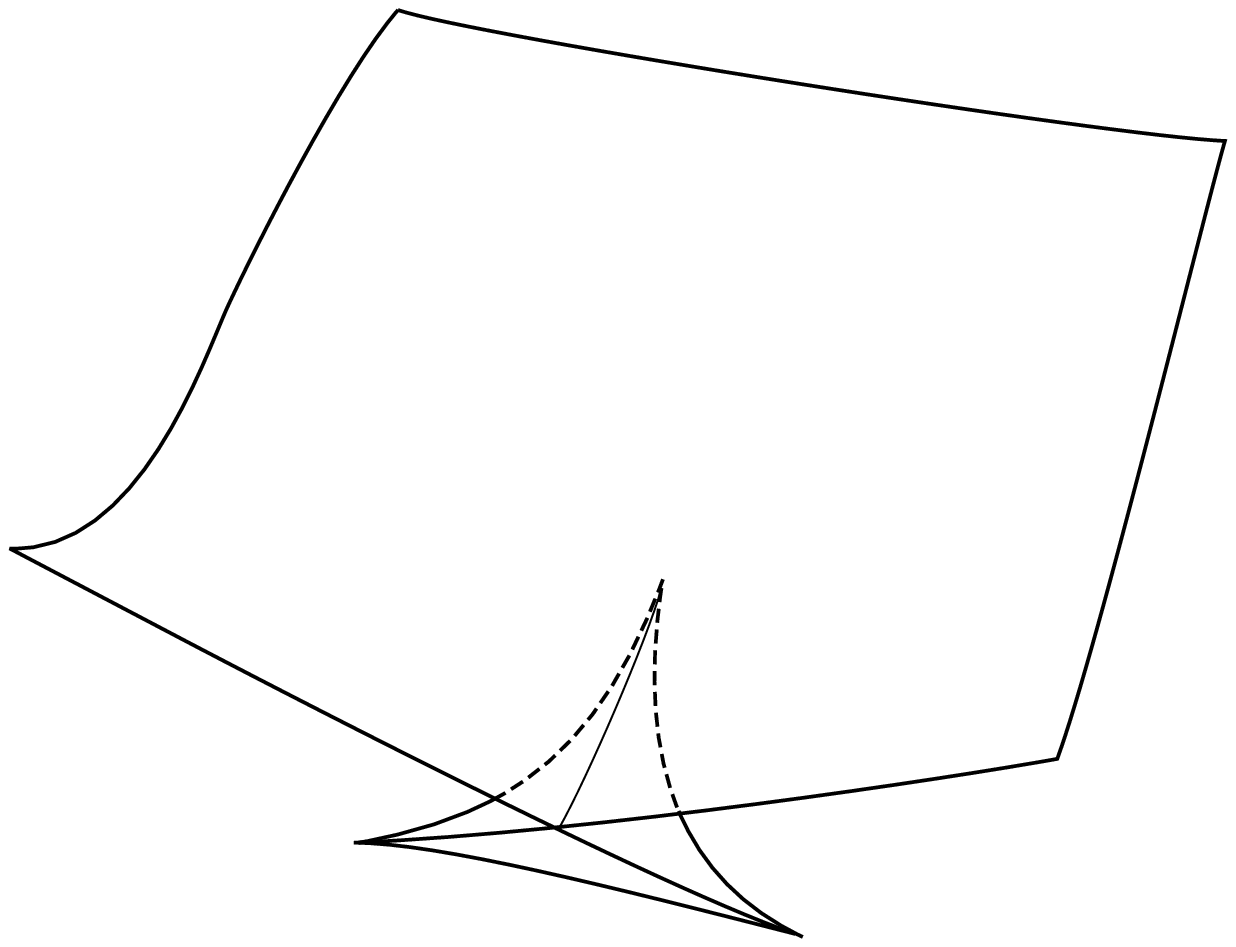}
\end{center}
 \end{minipage}
  \begin{minipage}{0.30\hsize}
  \begin{center}
   \includegraphics*[width=4cm,height=4cm]{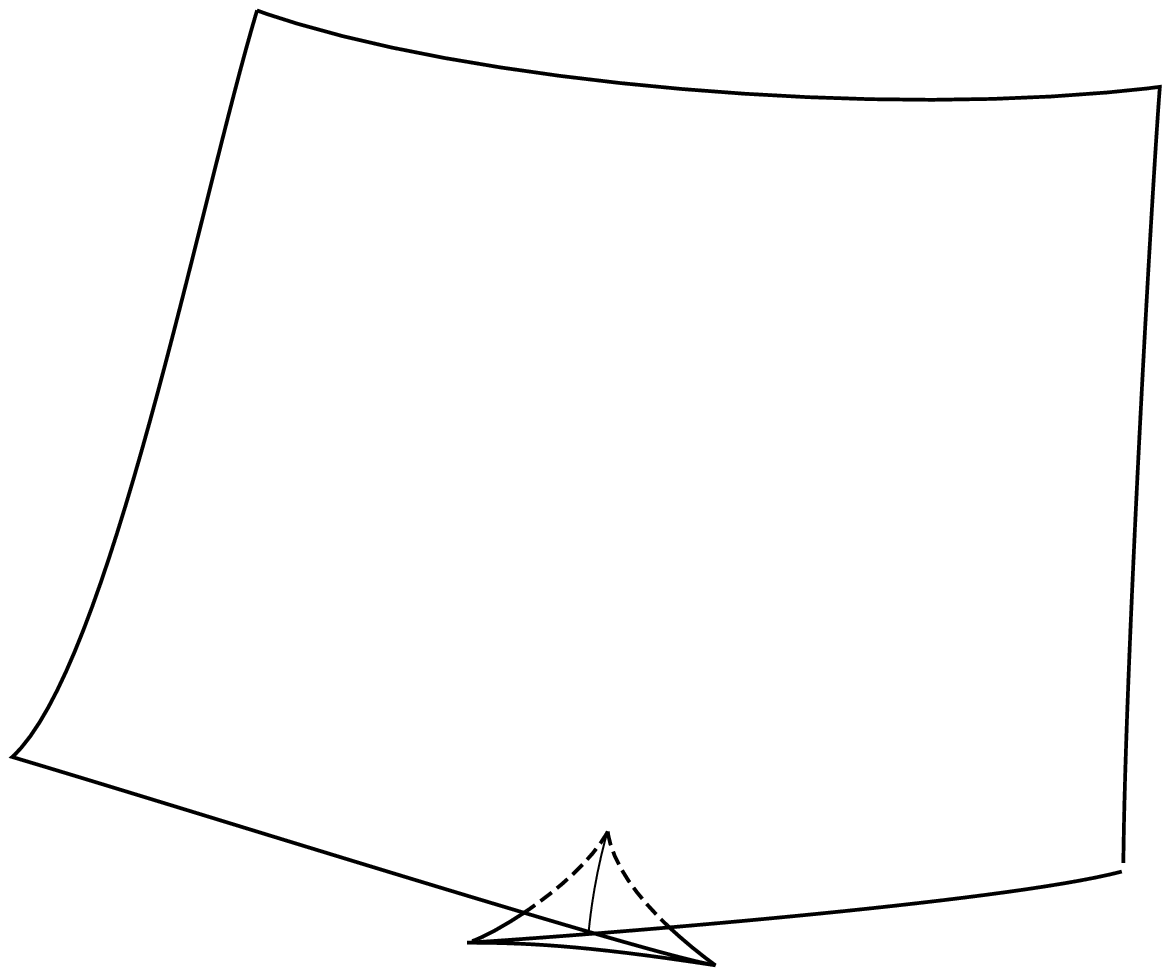}
\end{center}
 \end{minipage}
 \end{center}
 \end{figure}
 \begin{figure}[ht]
\begin{center}
 \begin{minipage}{0.30\hsize}
  \begin{center}
    \includegraphics*[width=4cm,height=4cm]{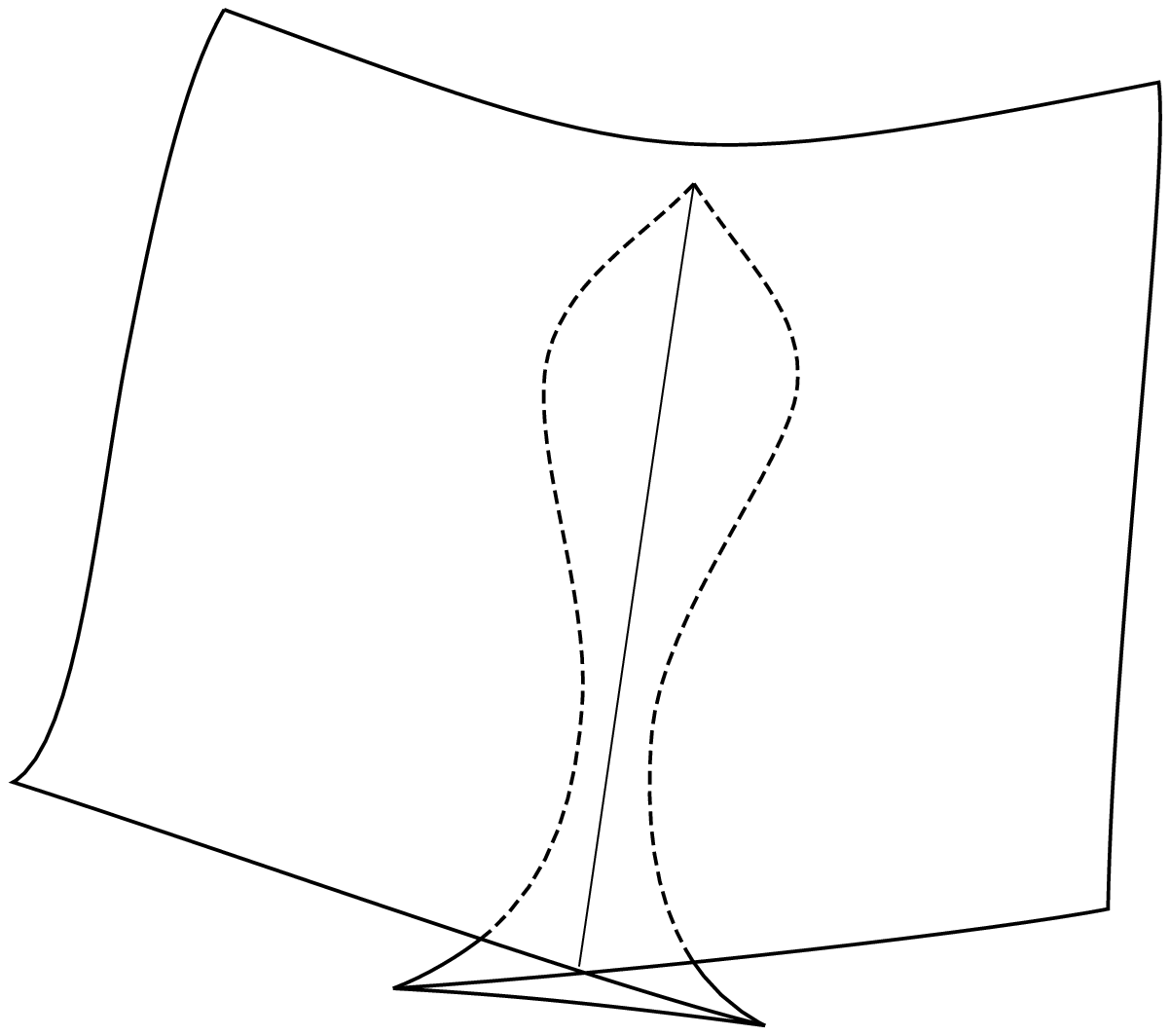} 
\end{center}
 \end{minipage}
 \begin{minipage}{0.30\hsize}
  \begin{center}
    \includegraphics*[width=4cm,height=4cm]{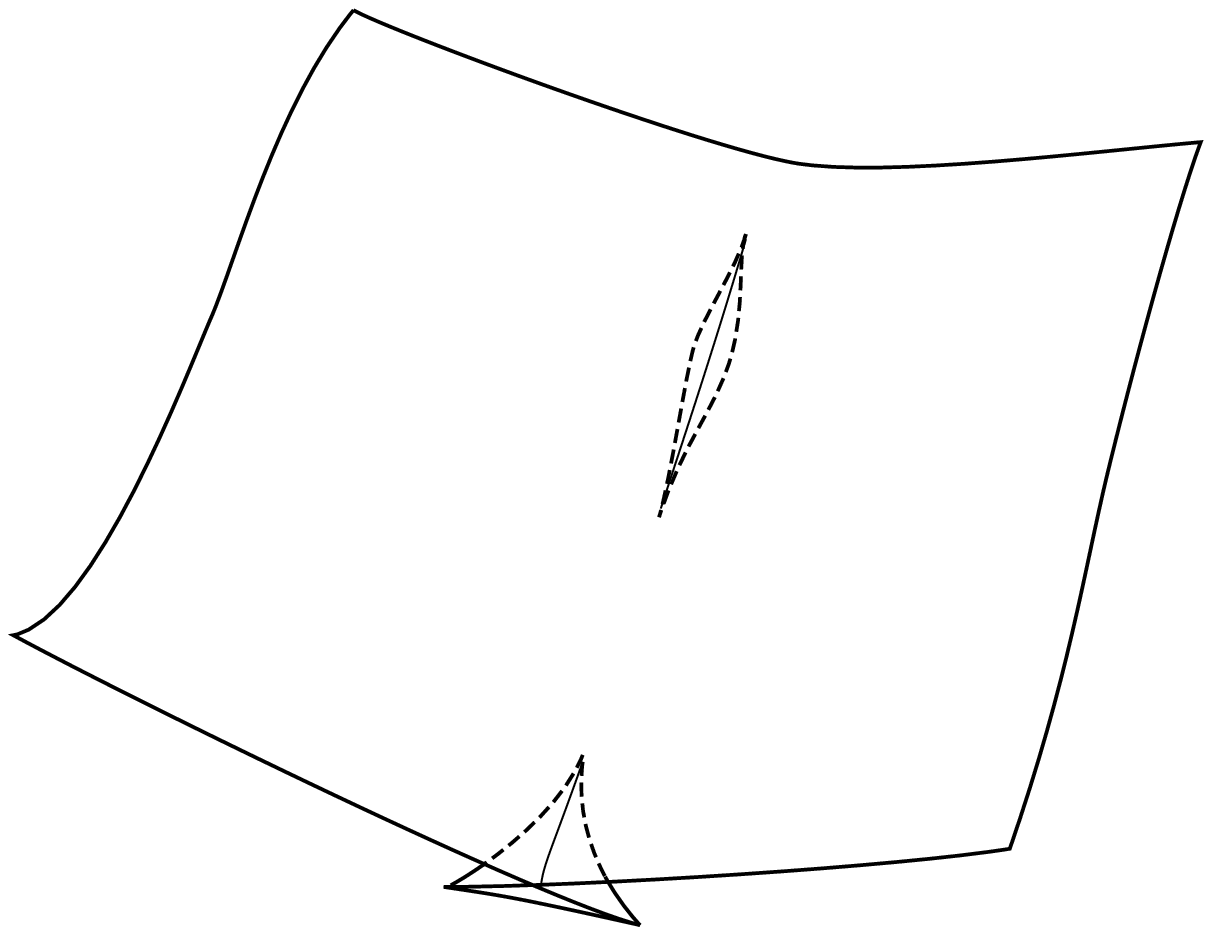}
\end{center}
 \end{minipage}
  \begin{minipage}{0.30\hsize}
  \begin{center}
   \includegraphics*[width=4cm,height=4cm]{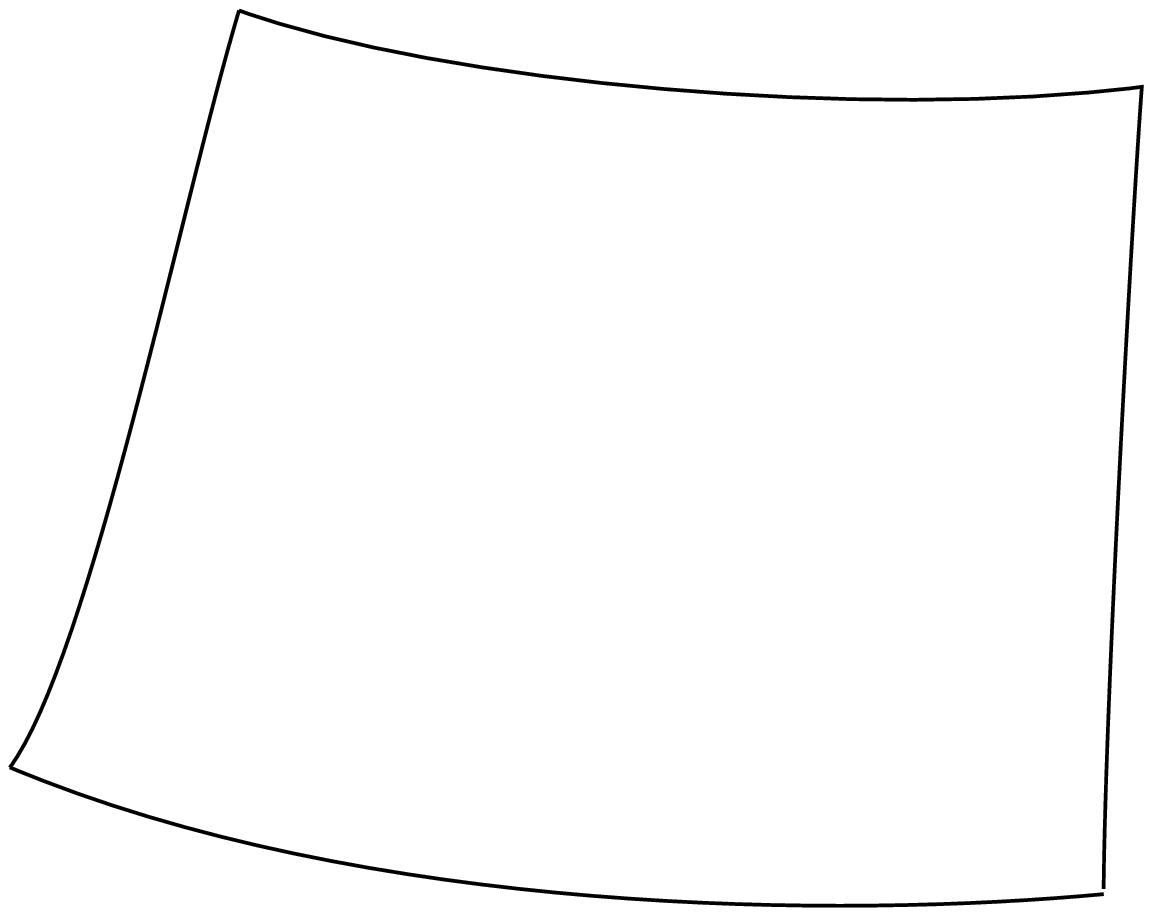}
\end{center}
 \end{minipage}
\end{center}
\caption{${}^2A_3$}
\end{figure}
\newpage
\begin{figure}[ht]
\begin{center}
 \begin{minipage}{0.30\hsize}
  \begin{center}
    \includegraphics*[width=4cm,height=4cm]{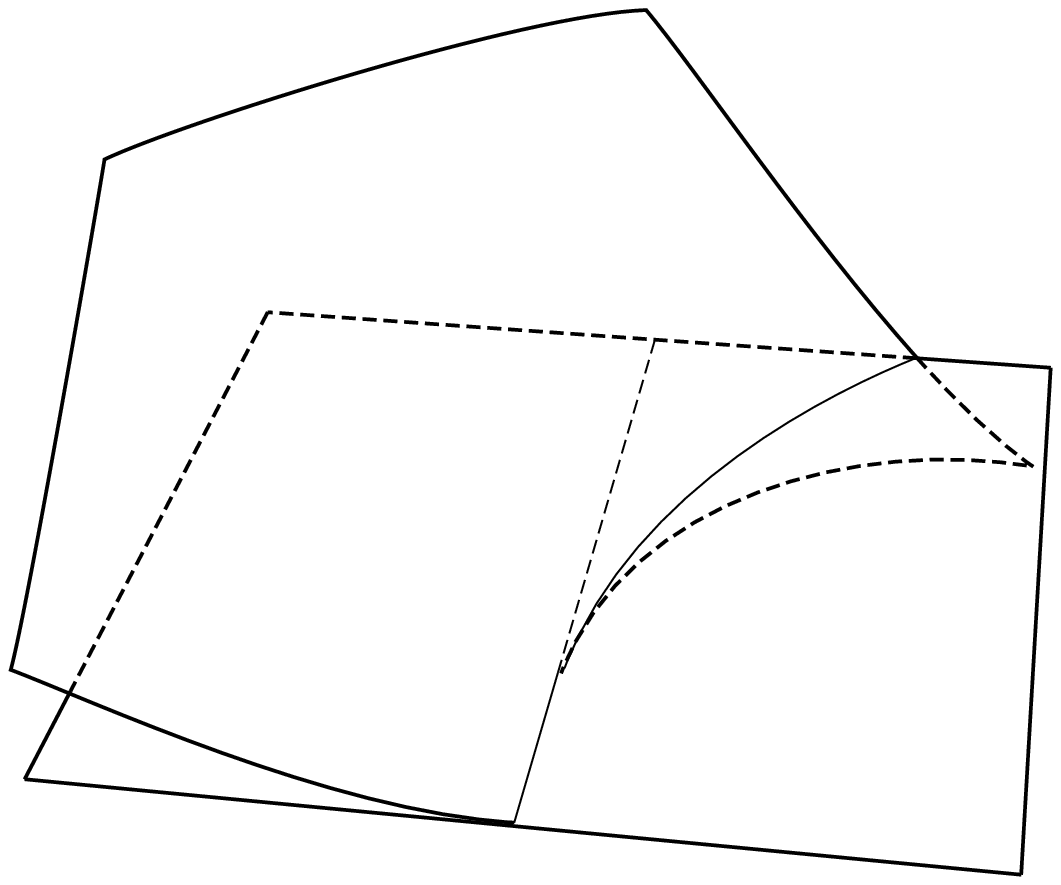} 
\end{center}
 \end{minipage}
 \begin{minipage}{0.30\hsize}
  \begin{center}
    \includegraphics*[width=4cm,height=4cm]{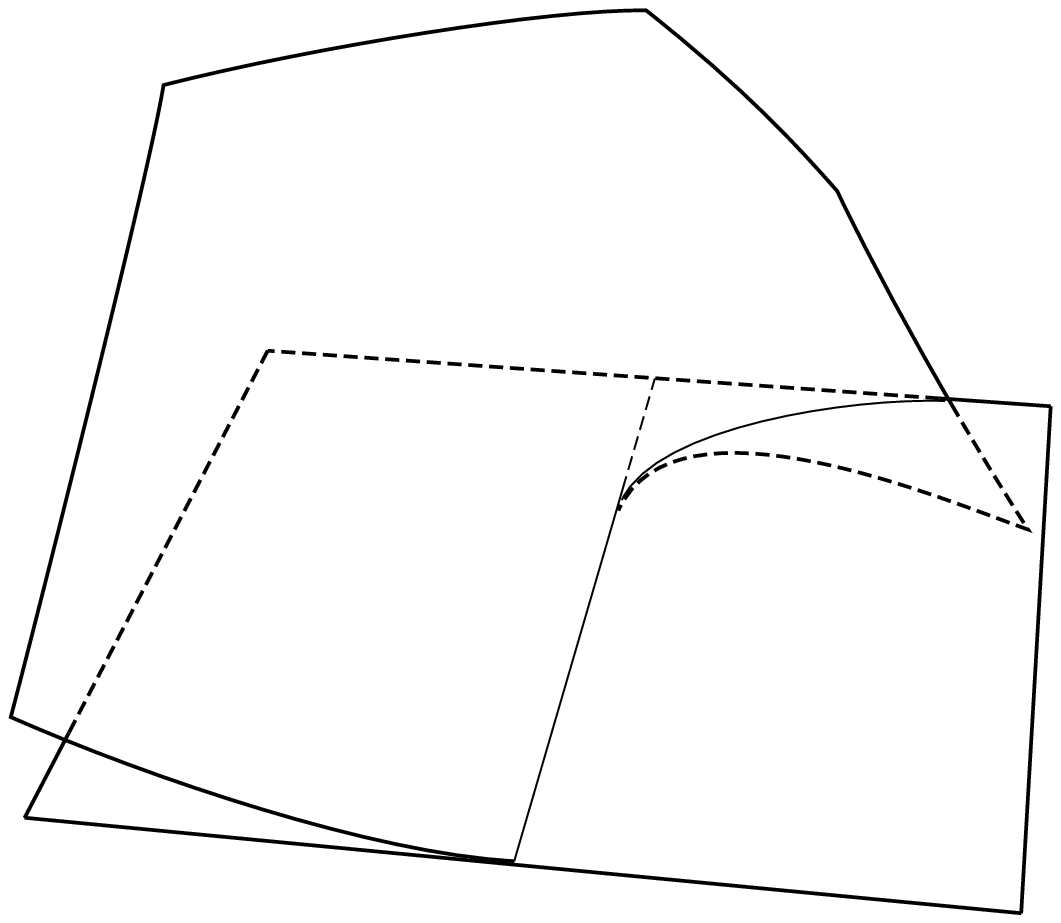}
\end{center}
 \end{minipage}
  \begin{minipage}{0.30\hsize}
  \begin{center}
   \includegraphics*[width=4cm,height=4cm]{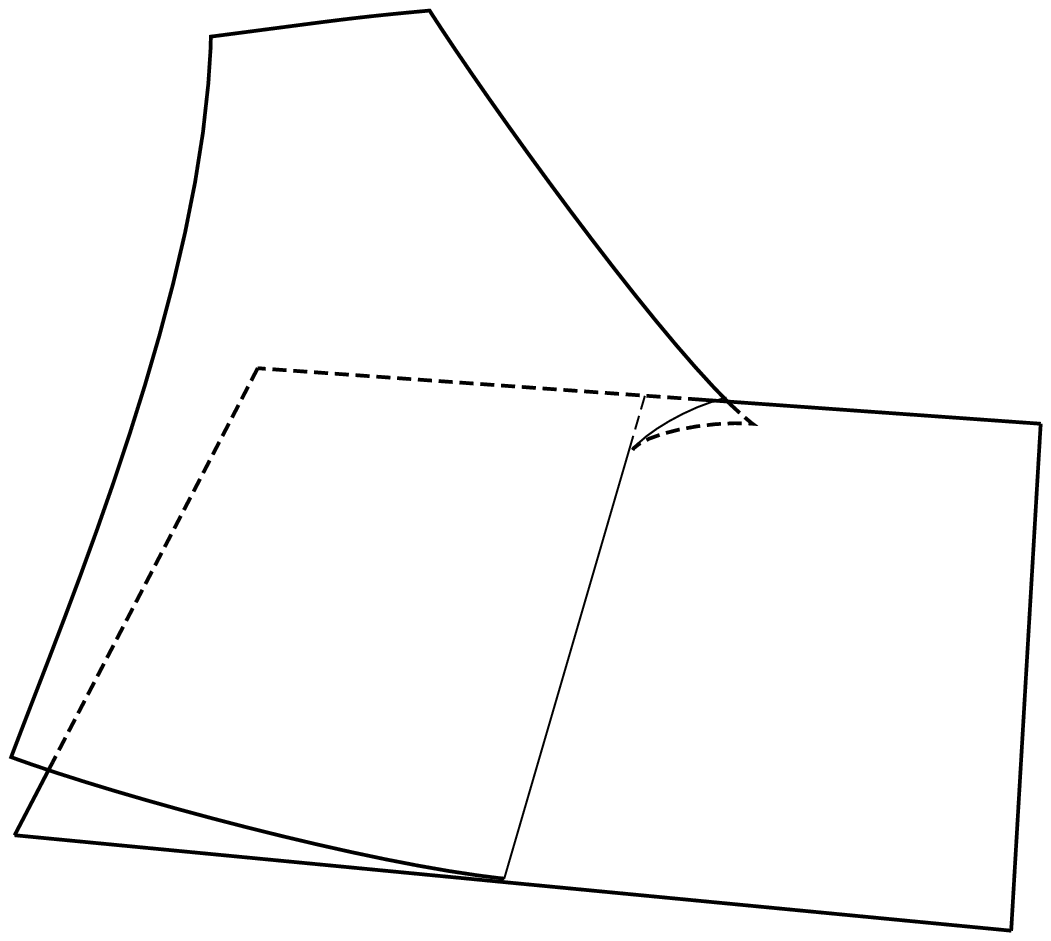}
\end{center}
 \end{minipage}
 \end{center}
 \end{figure}
 \begin{figure}[ht]
\begin{center}
 \begin{minipage}{0.30\hsize}
  \begin{center}
    \includegraphics*[width=4cm,height=4cm]{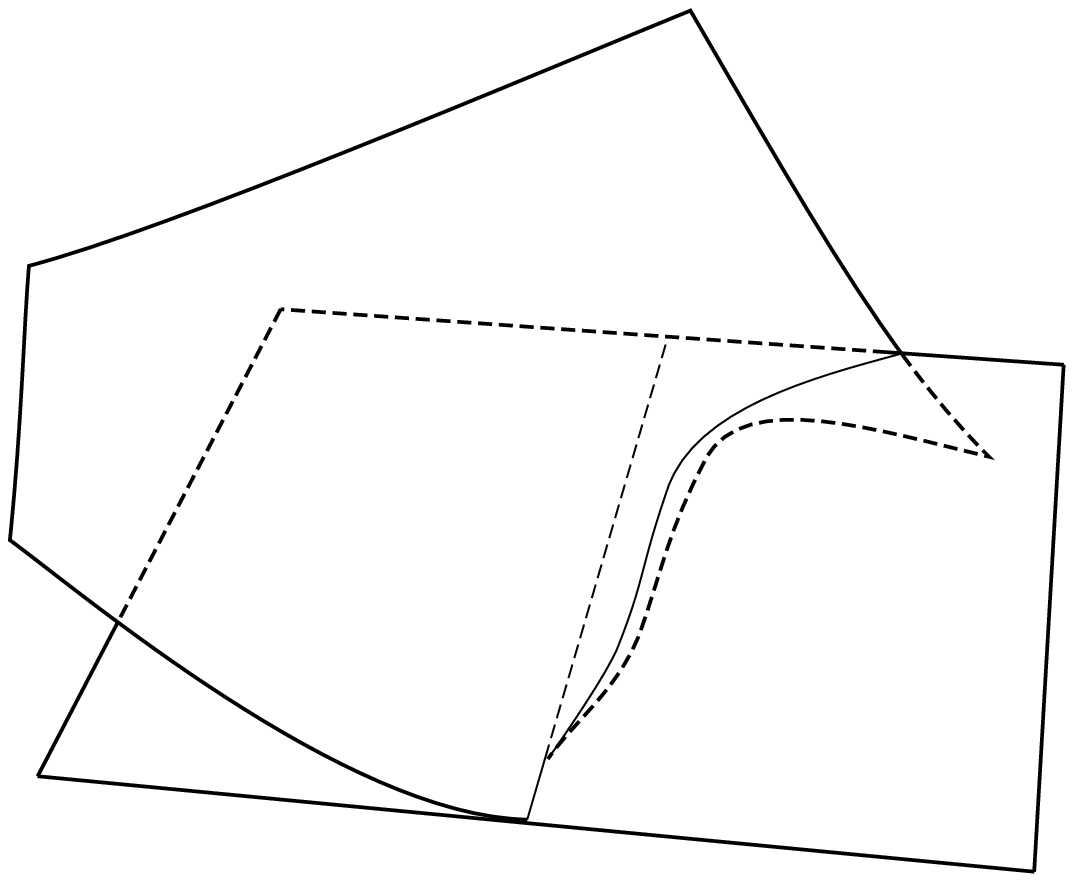} 
\end{center}
 \end{minipage}
 \begin{minipage}{0.30\hsize}
  \begin{center}
    \includegraphics*[width=4cm,height=4cm]{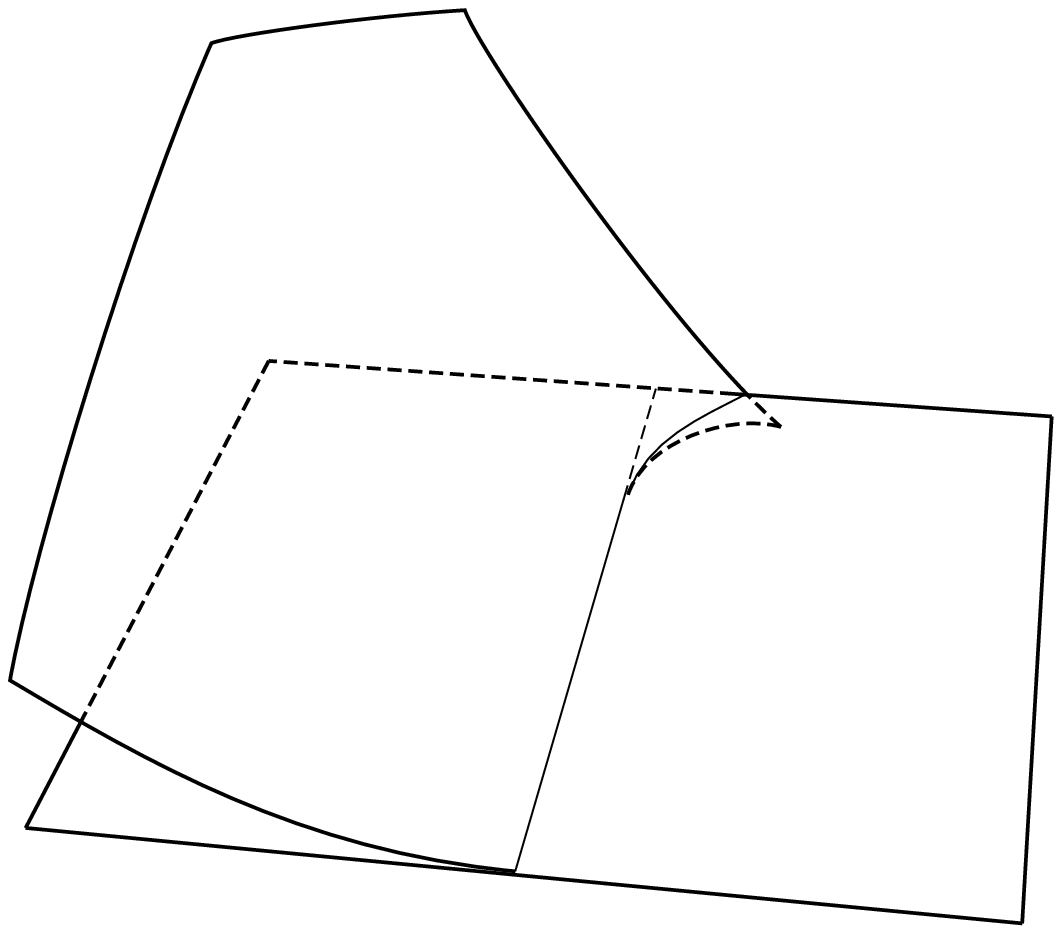}
\end{center}
 \end{minipage}
  \begin{minipage}{0.30\hsize}
  \begin{center}
   \includegraphics*[width=4cm,height=4cm]{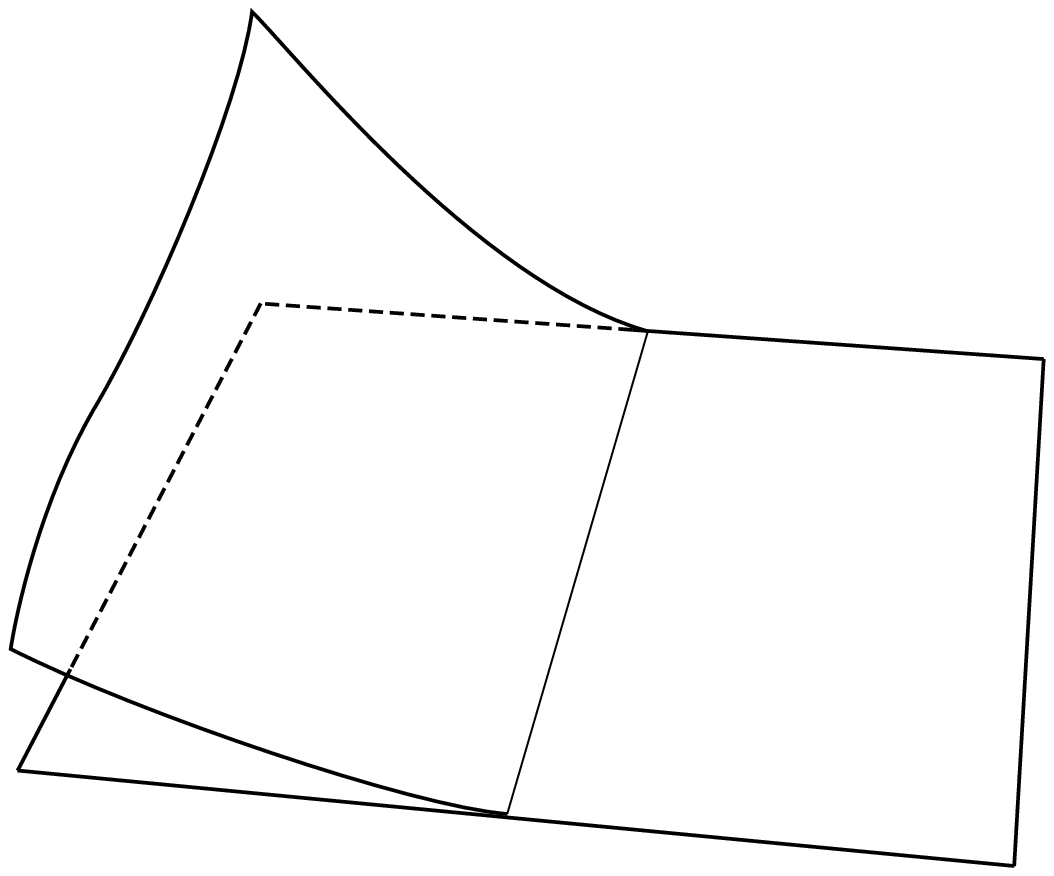}
\end{center}
 \end{minipage}
 \end{center}
 \end{figure}
 \begin{figure}[ht]
\begin{center}
 \begin{minipage}{0.30\hsize}
  \begin{center}
    \includegraphics*[width=4cm,height=4cm]{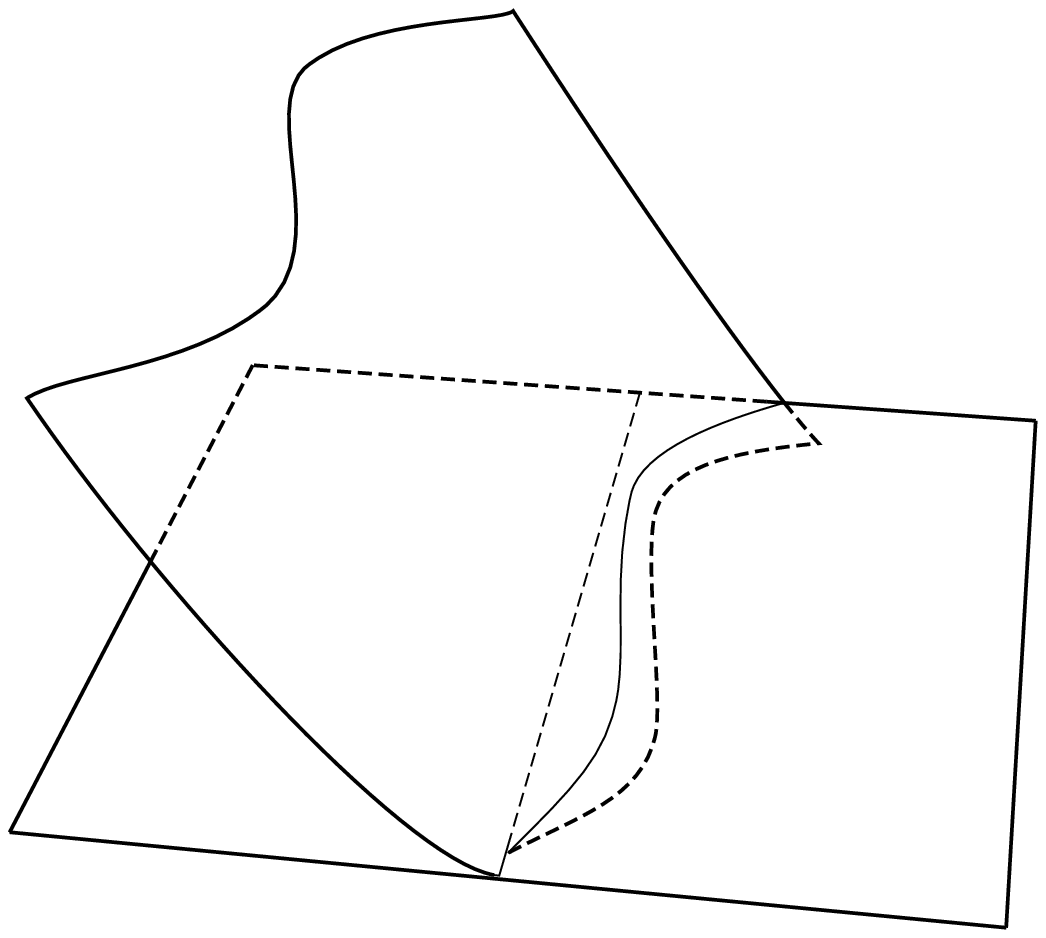} 
\end{center}
 \end{minipage}
 \begin{minipage}{0.30\hsize}
  \begin{center}
    \includegraphics*[width=4cm,height=4cm]{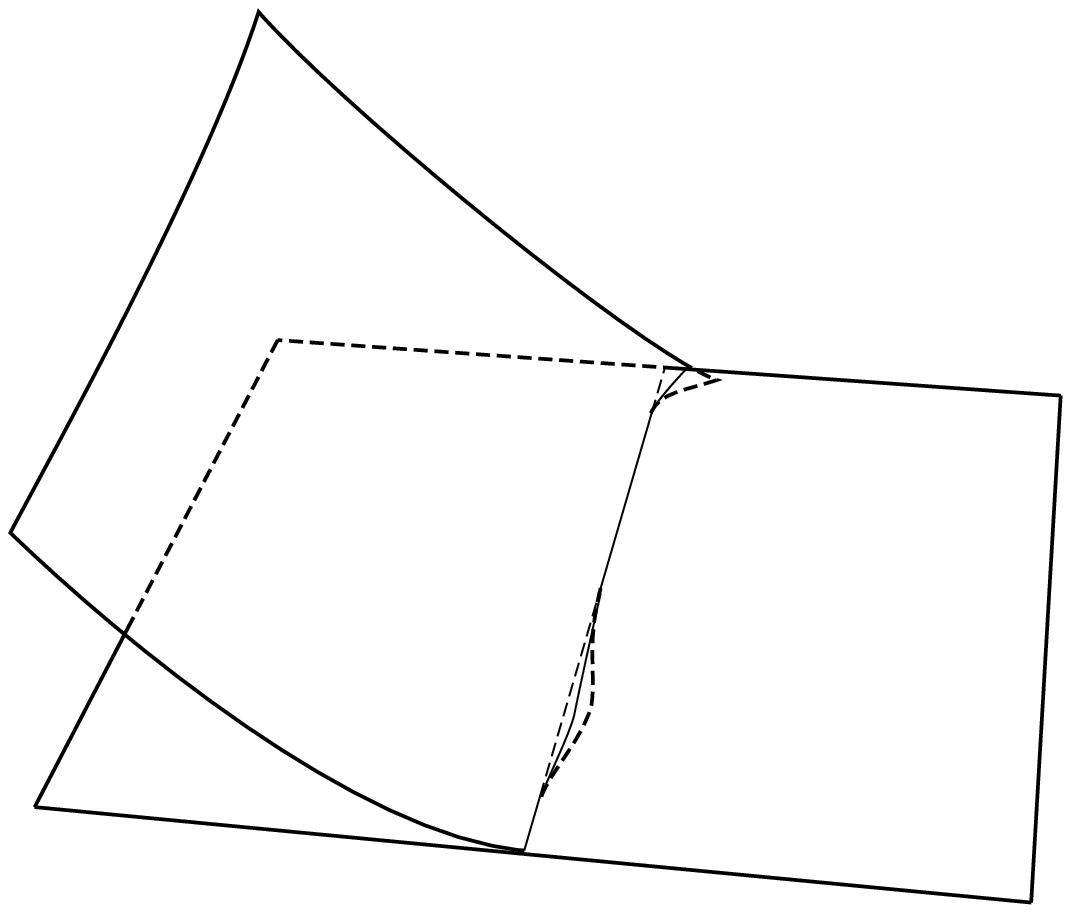}
\end{center}
 \end{minipage}
  \begin{minipage}{0.30\hsize}
  \begin{center}
   \includegraphics*[width=4cm,height=4cm]{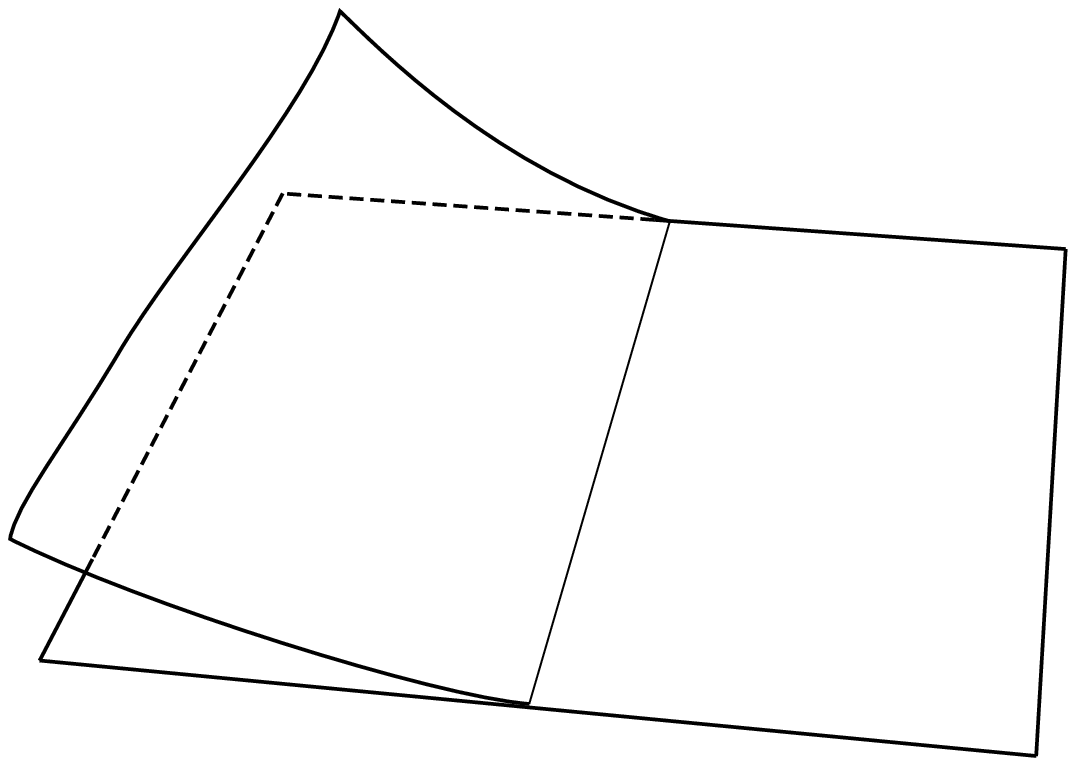}
\end{center}
 \end{minipage}
\end{center}
\caption{${}^2B_3$}
\end{figure}
\newpage
\begin{figure}[ht]
\begin{center}
 \begin{minipage}{0.30\hsize}
  \begin{center}
    \includegraphics*[width=4cm,height=4cm]{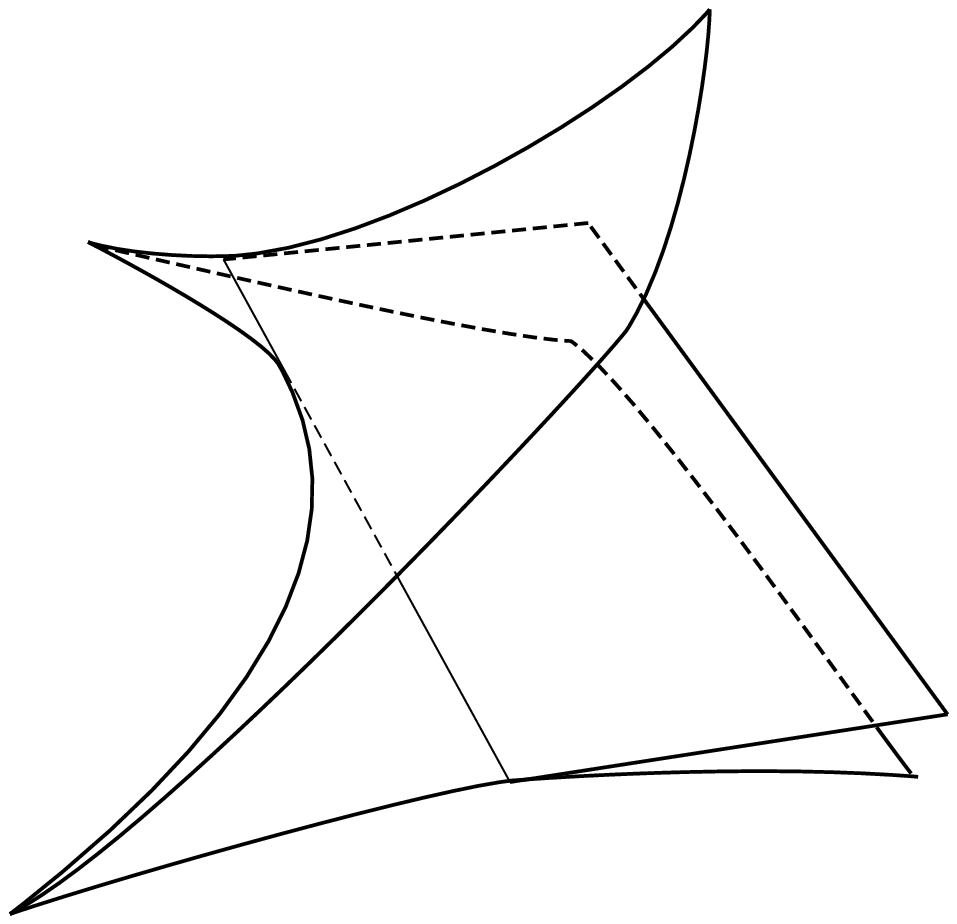} 
\end{center}
 \end{minipage}
 \begin{minipage}{0.30\hsize}
  \begin{center}
    \includegraphics*[width=4cm,height=4cm]{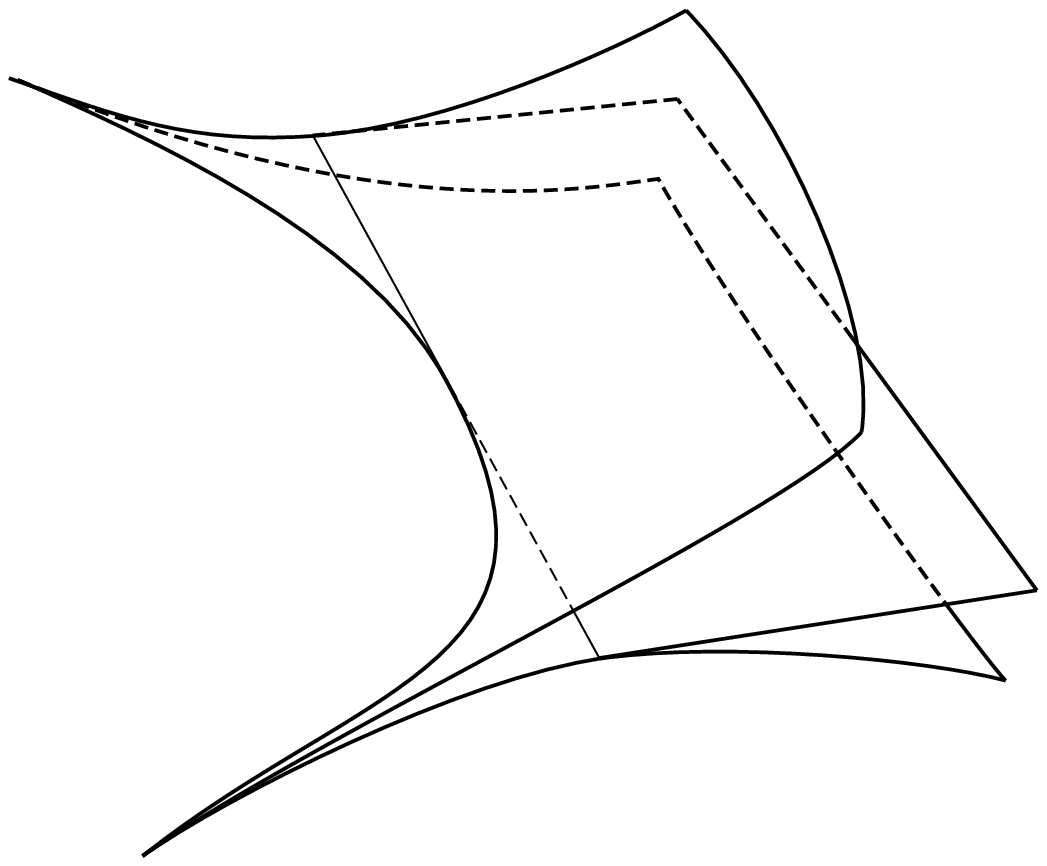}
\end{center}
 \end{minipage}
  \begin{minipage}{0.30\hsize}
  \begin{center}
   \includegraphics*[width=4cm,height=4cm]{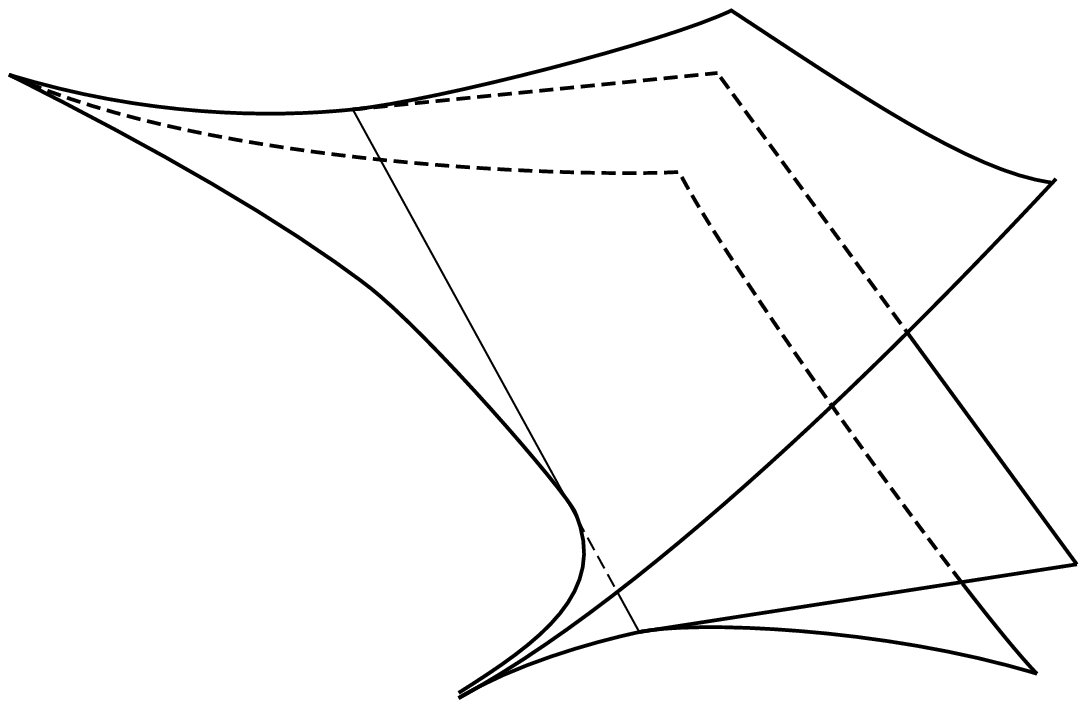}
\end{center}
 \end{minipage}
 \end{center}
 \end{figure}
 \begin{figure}[ht]
\begin{center}
 \begin{minipage}{0.30\hsize}
  \begin{center}
    \includegraphics*[width=4cm,height=4cm]{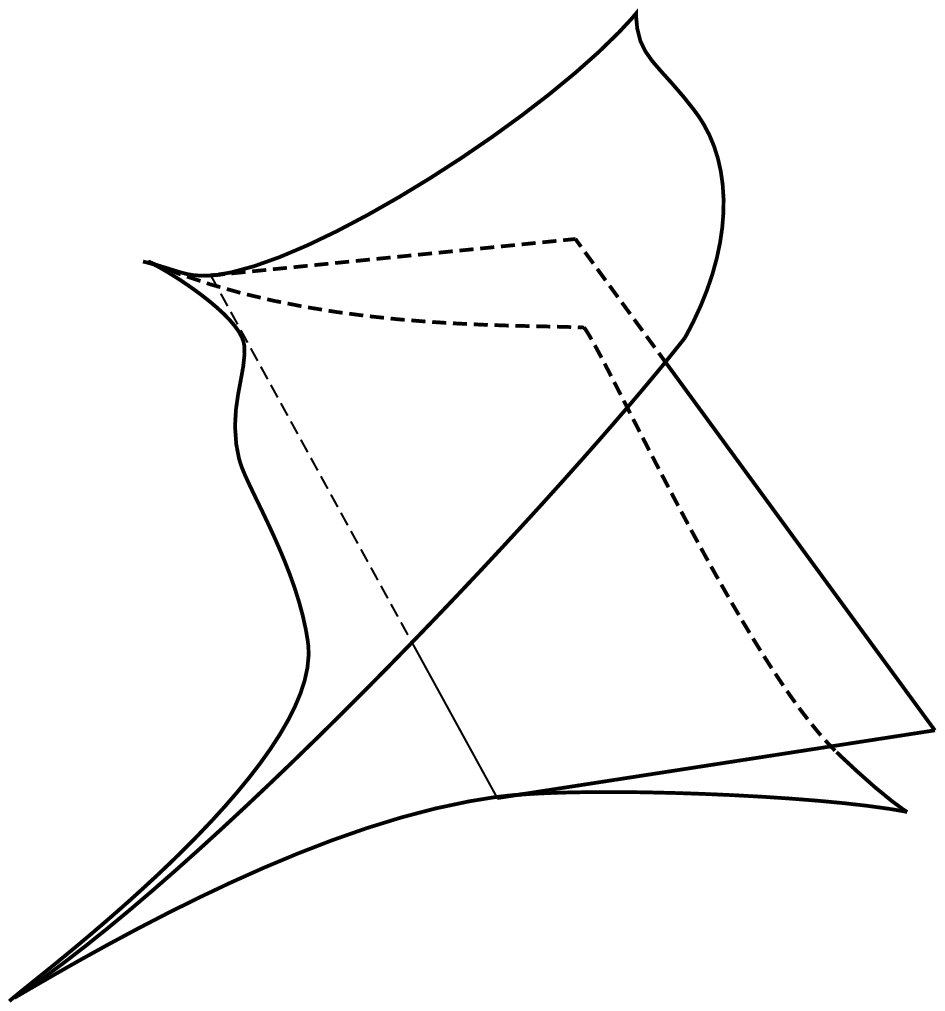} 
\end{center}
 \end{minipage}
 \begin{minipage}{0.30\hsize}
  \begin{center}
    \includegraphics*[width=4cm,height=4cm]{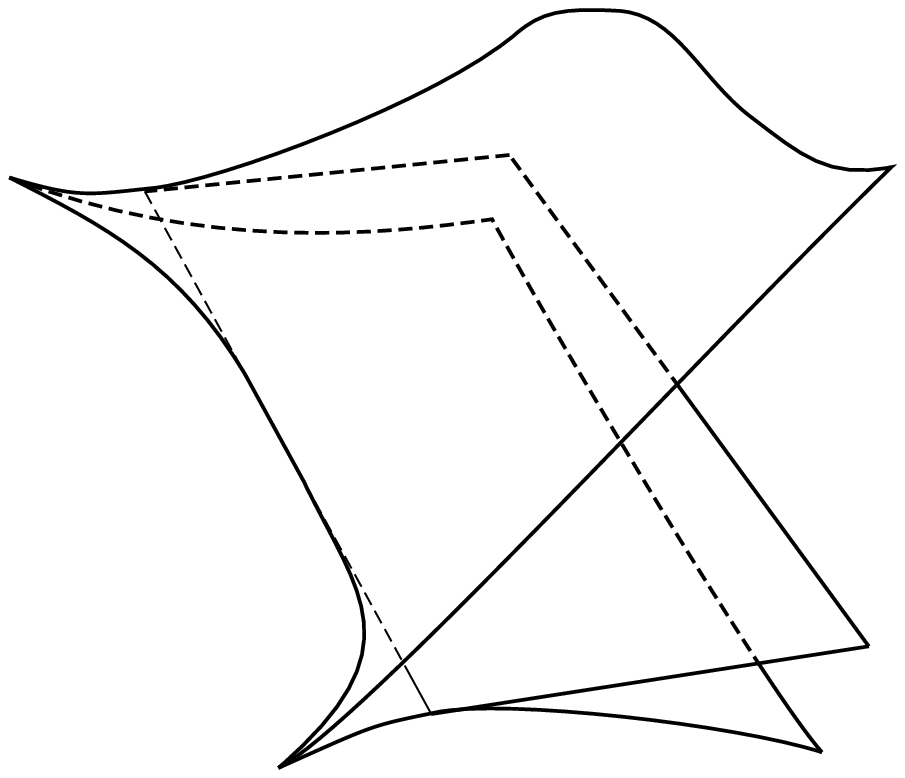}
\end{center}
 \end{minipage}
  \begin{minipage}{0.30\hsize}
  \begin{center}
   \includegraphics*[width=4cm,height=4cm]{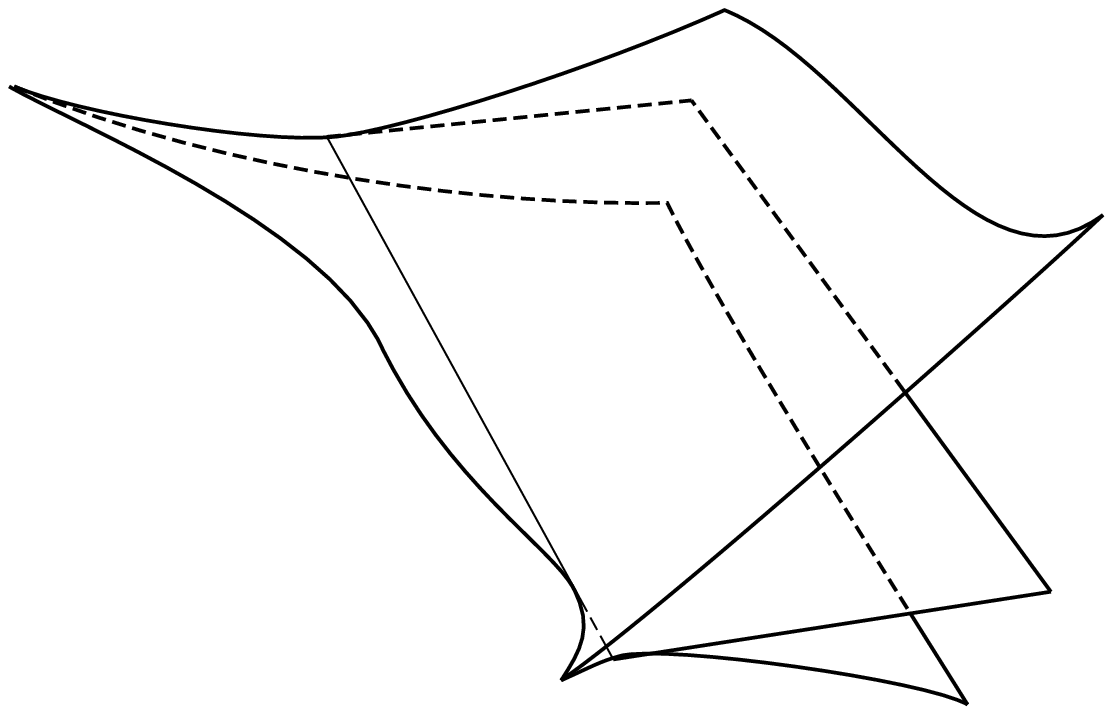}
\end{center}
 \end{minipage}
 \end{center}
 \end{figure}
 \begin{figure}[ht]
\begin{center}
 \begin{minipage}{0.30\hsize}
  \begin{center}
    \includegraphics*[width=4cm,height=4cm]{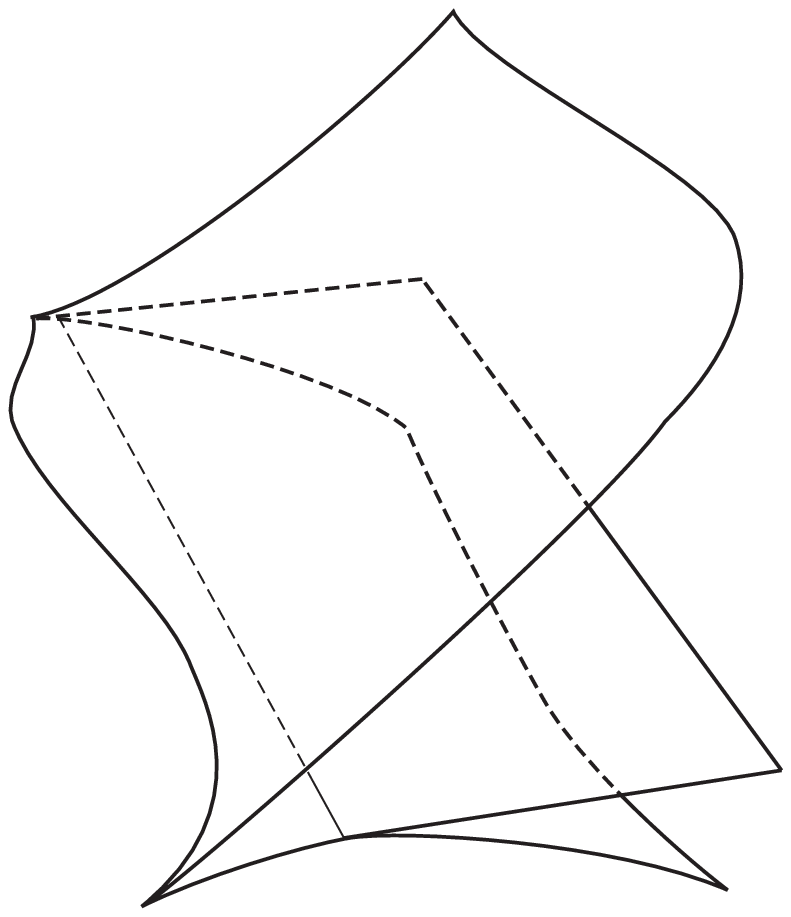} 
\end{center}
 \end{minipage}
 \begin{minipage}{0.30\hsize}
  \begin{center}
    \includegraphics*[width=4cm,height=4cm]{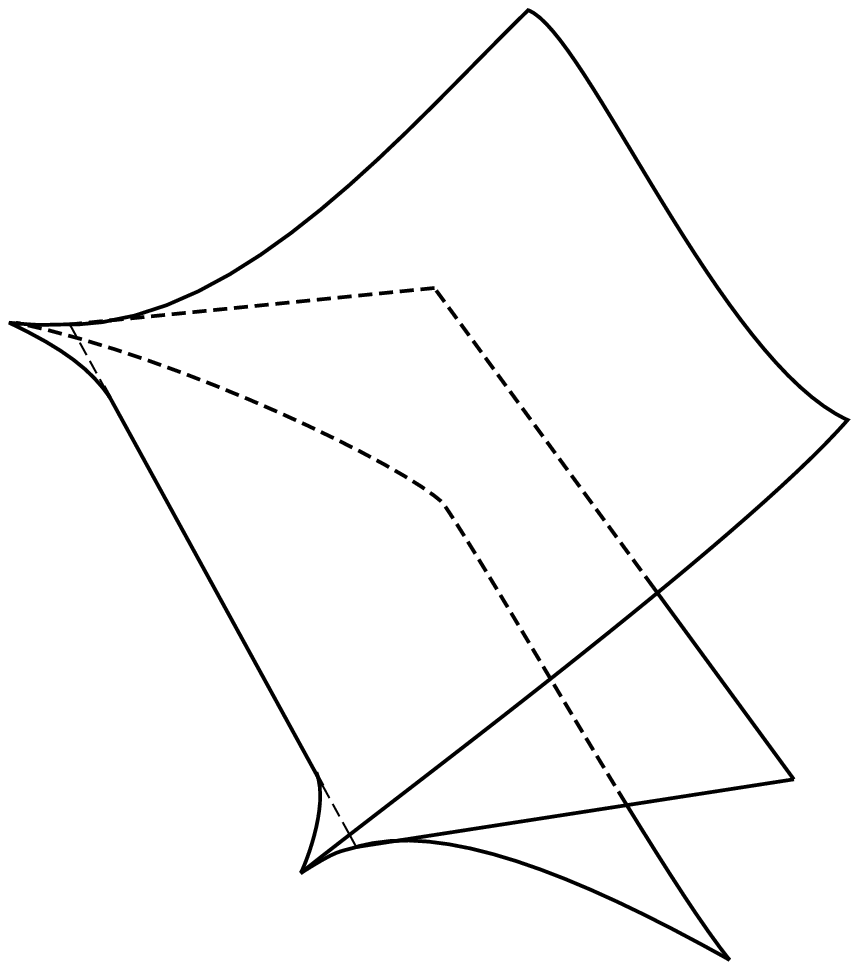}
\end{center}
 \end{minipage}
  \begin{minipage}{0.30\hsize}
  \begin{center}
   \includegraphics*[width=4cm,height=4cm]{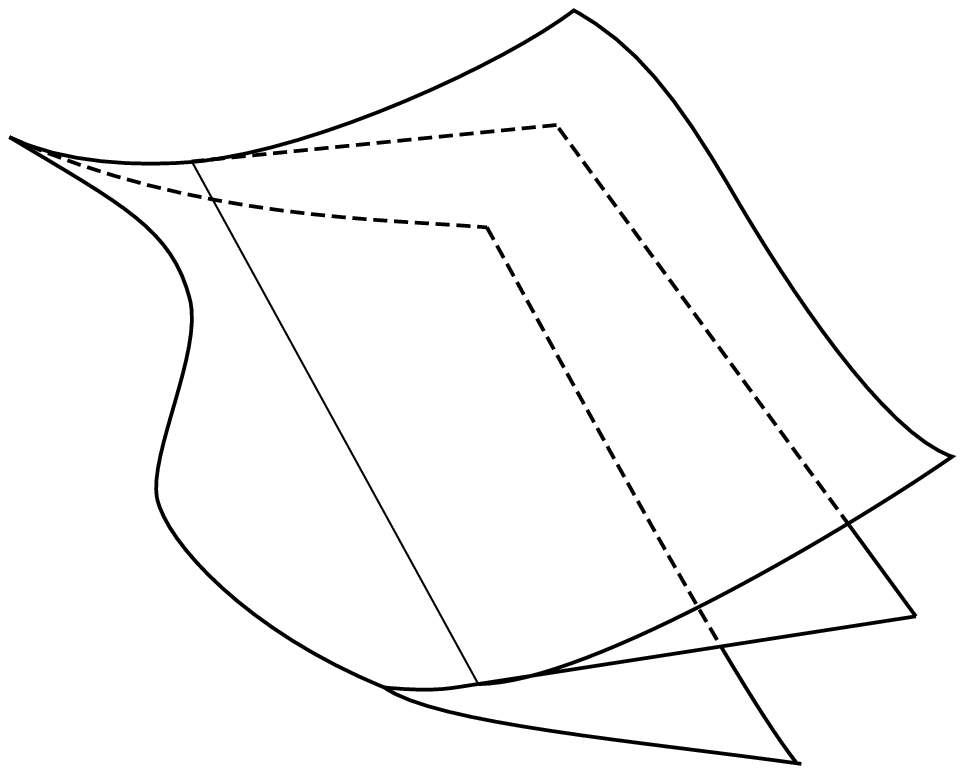}
\end{center}
 \end{minipage}
\end{center}
\caption{${}^2C_3^+$}
\end{figure}
\newpage
\begin{figure}[ht]
\begin{center}
 \begin{minipage}{0.30\hsize}
  \begin{center}
    \includegraphics*[width=4cm,height=4cm]{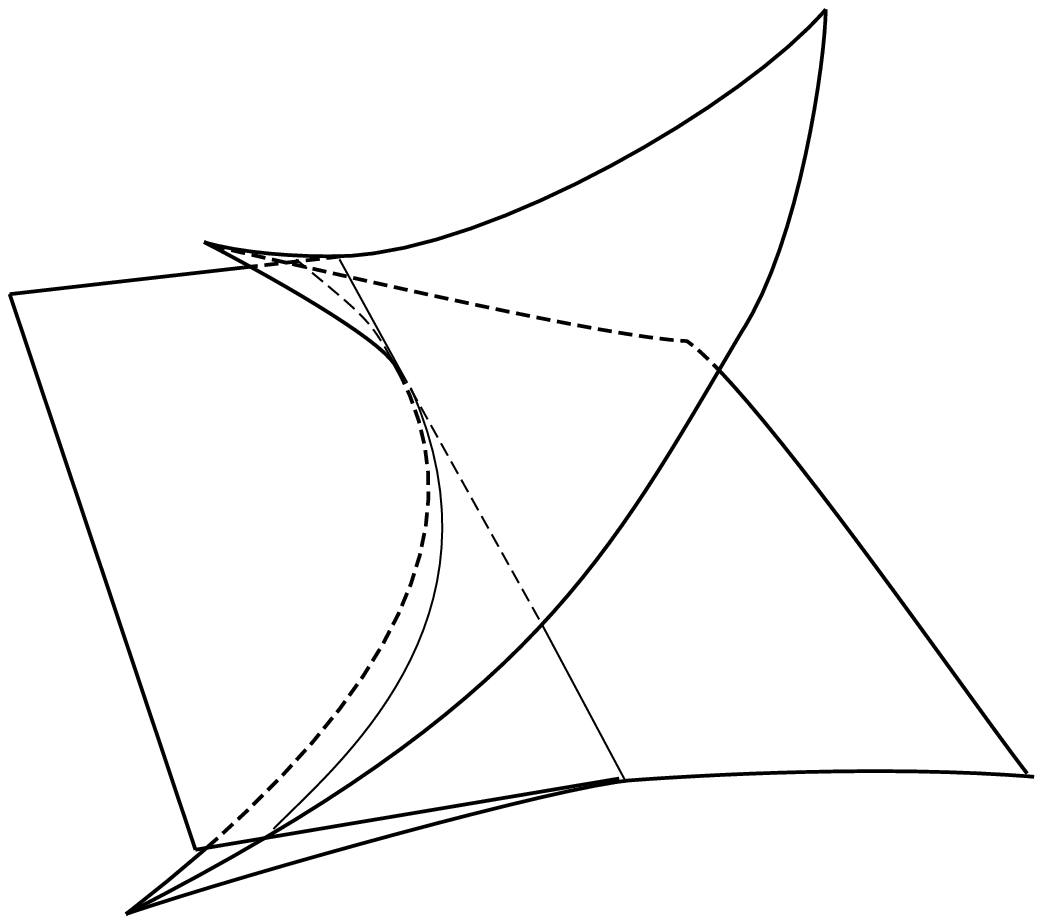} 
\end{center}
 \end{minipage}
 \begin{minipage}{0.30\hsize}
  \begin{center}
    \includegraphics*[width=4cm,height=4cm]{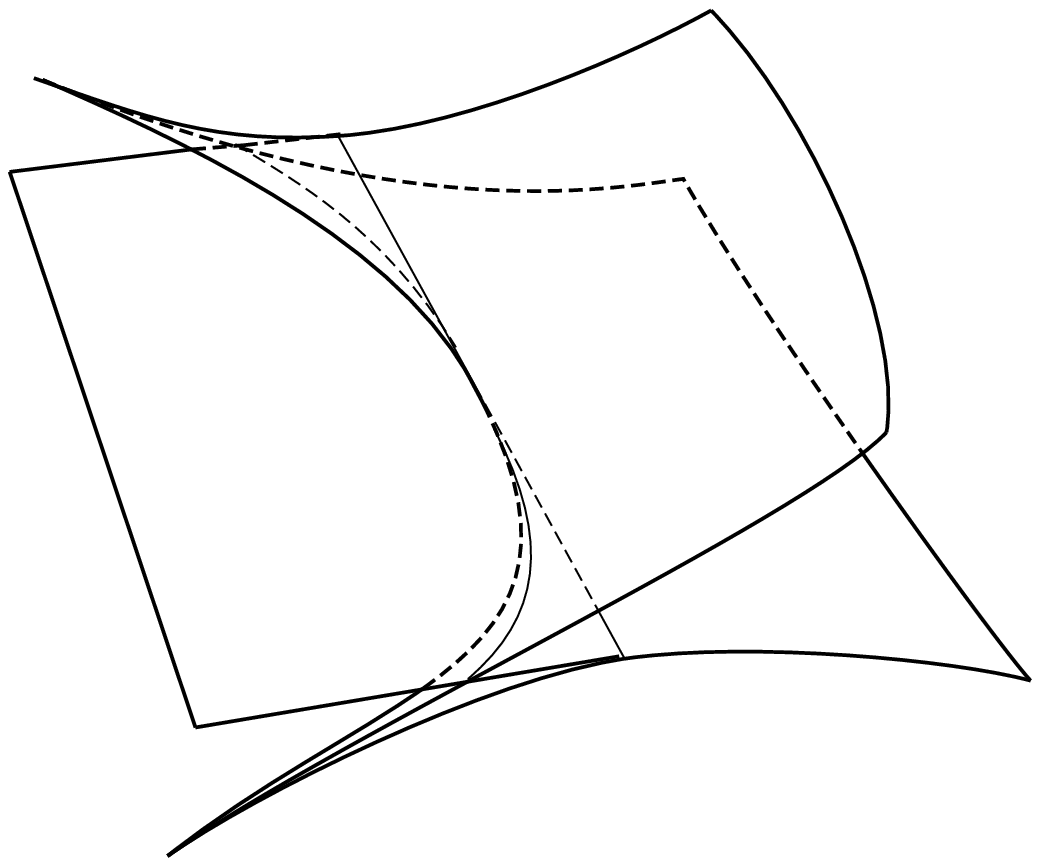}
\end{center}
 \end{minipage}
  \begin{minipage}{0.30\hsize}
  \begin{center}
   \includegraphics*[width=4cm,height=4cm]{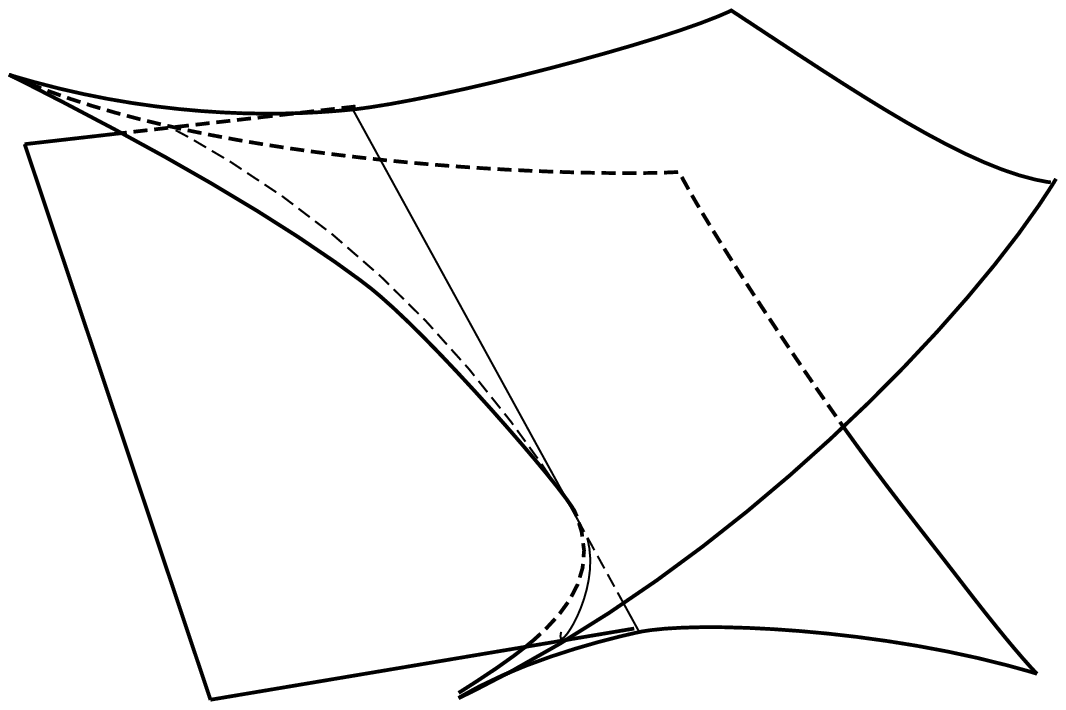}
\end{center}
 \end{minipage}
 \end{center}
 \end{figure}
 \begin{figure}[ht]
\begin{center}
 \begin{minipage}{0.30\hsize}
  \begin{center}
    \includegraphics*[width=4cm,height=4cm]{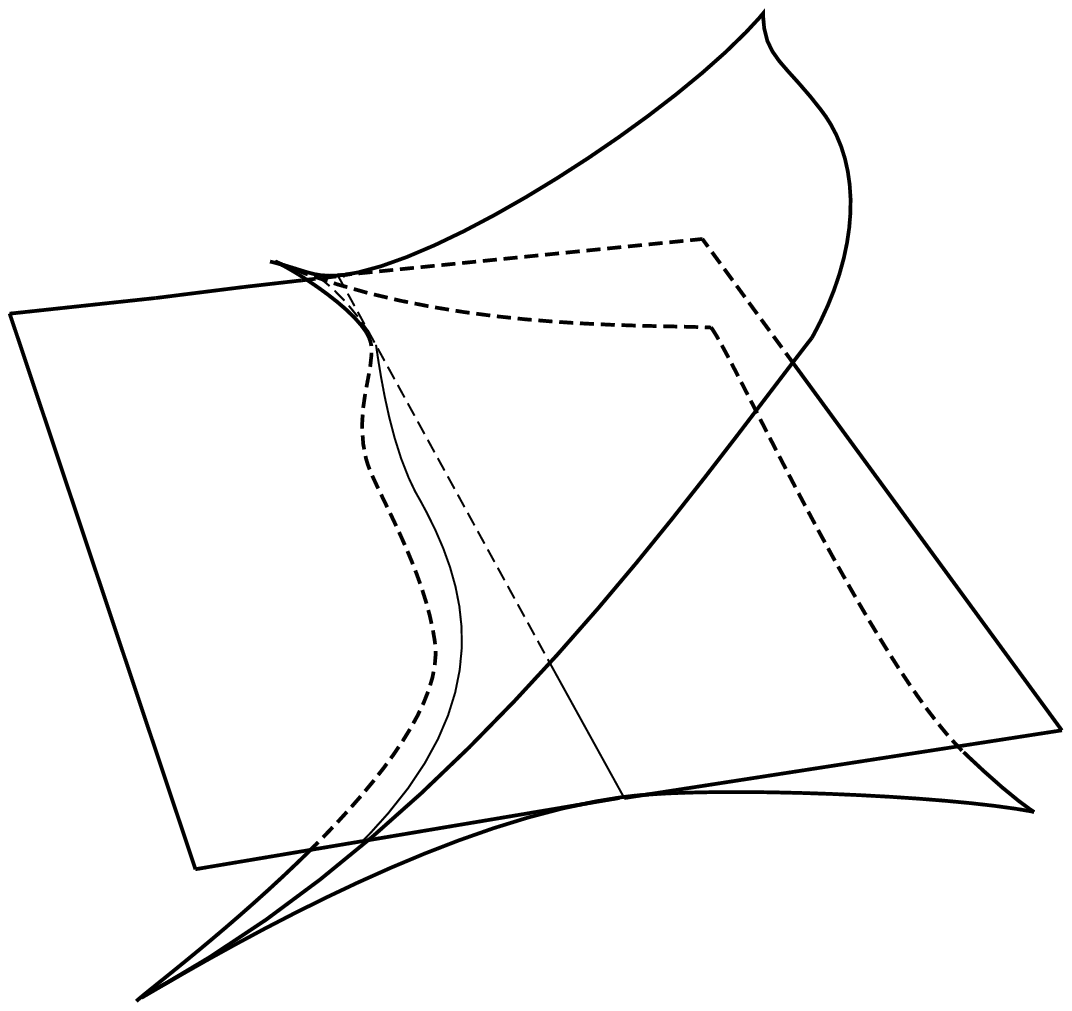} 
\end{center}
 \end{minipage}
 \begin{minipage}{0.30\hsize}
  \begin{center}
    \includegraphics*[width=4cm,height=4cm]{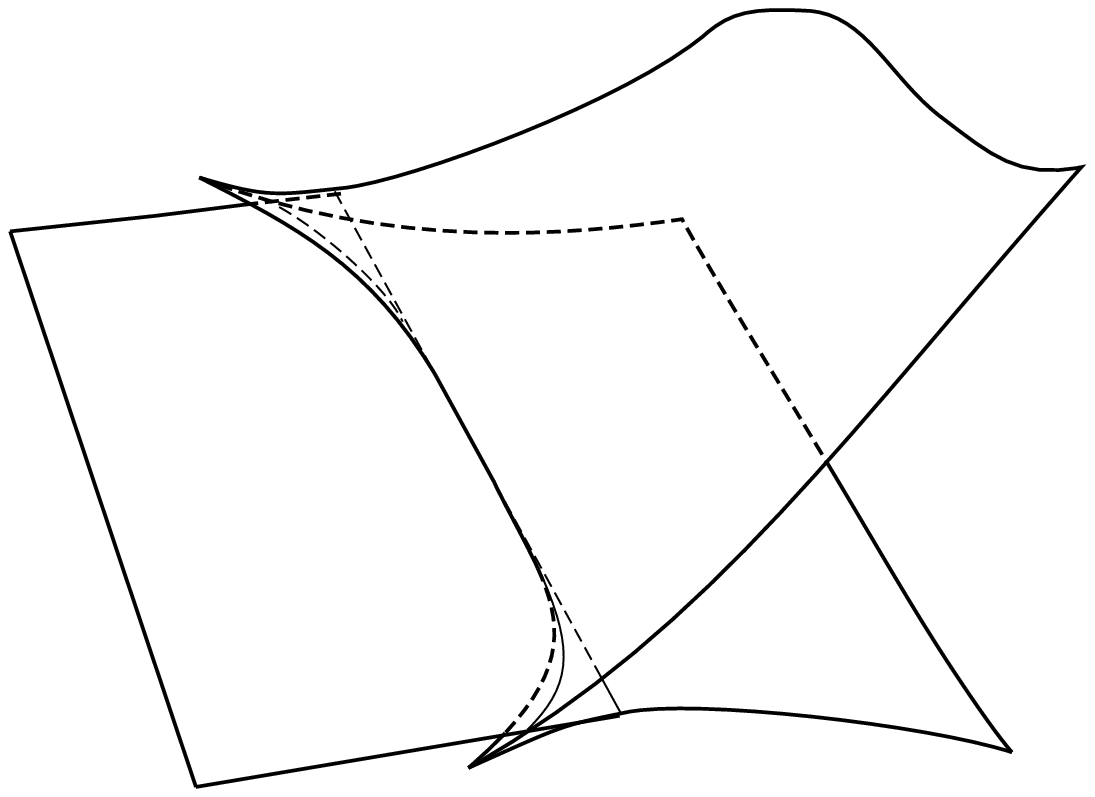}
\end{center}
 \end{minipage}
  \begin{minipage}{0.30\hsize}
  \begin{center}
   \includegraphics*[width=4cm,height=4cm]{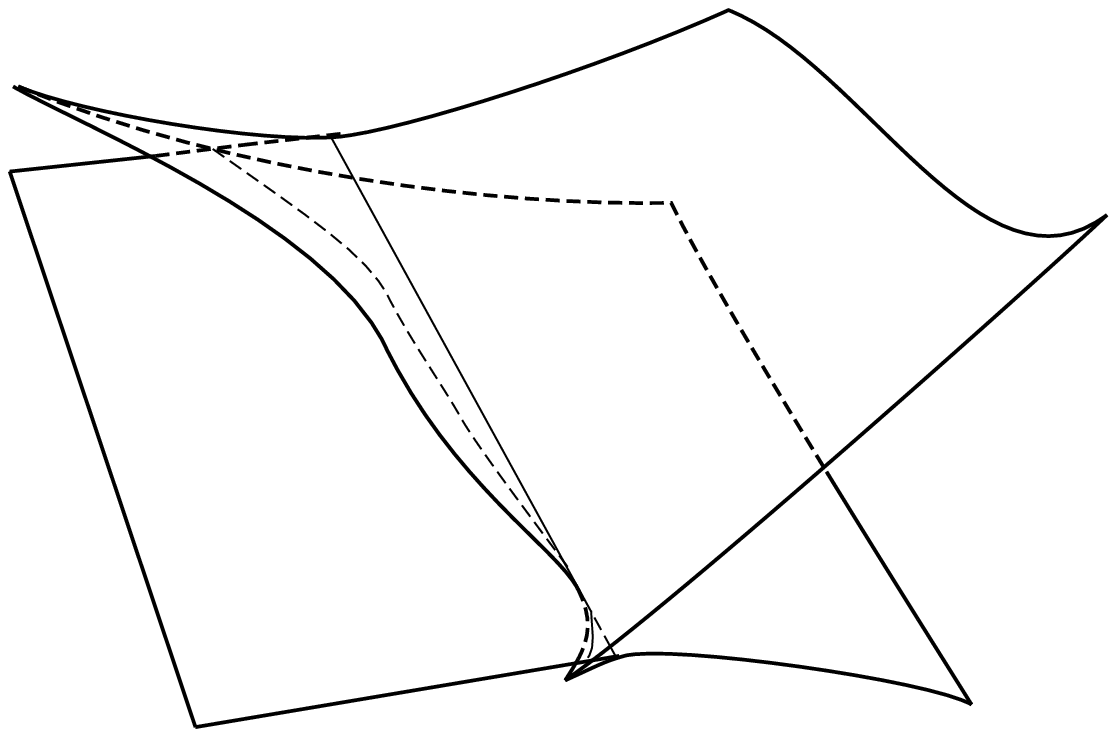}
\end{center}
 \end{minipage}
 \end{center}
 \end{figure}
 \begin{figure}[ht]
\begin{center}
 \begin{minipage}{0.30\hsize}
  \begin{center}
    \includegraphics*[width=4cm,height=4cm]{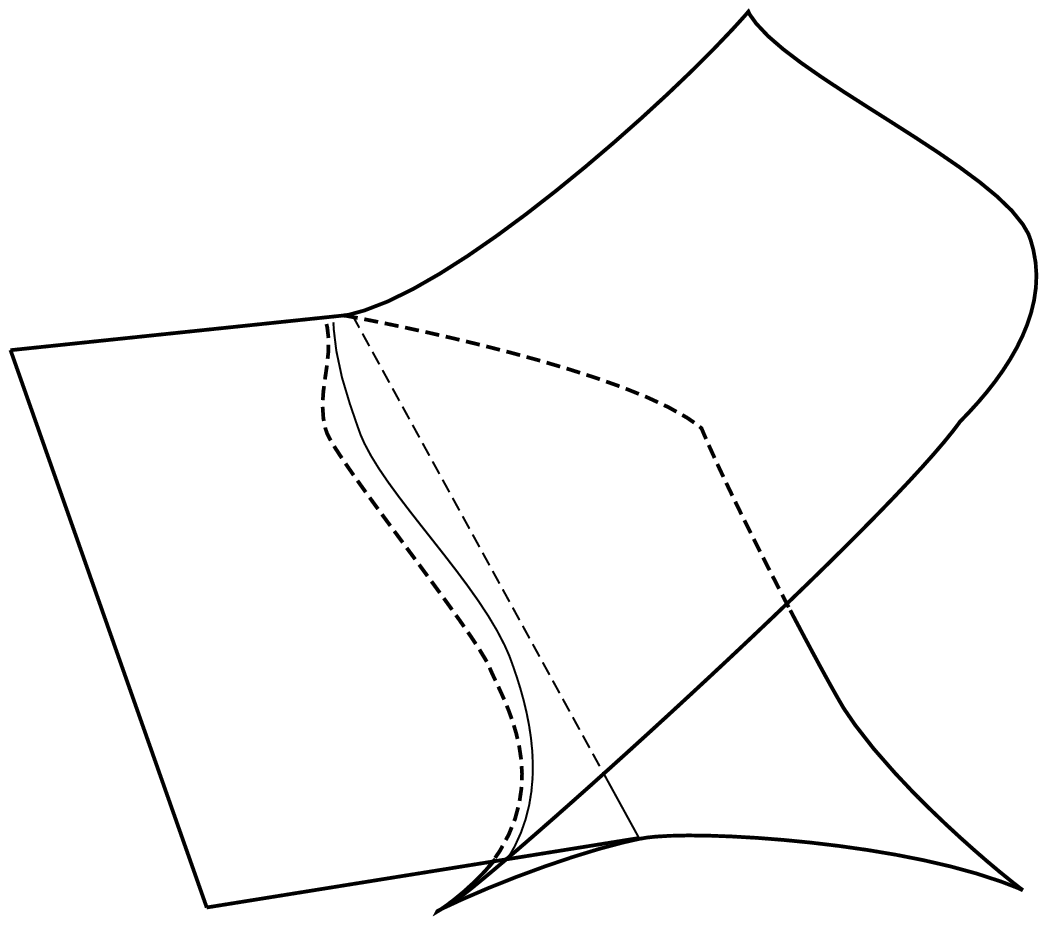} 
\end{center}
 \end{minipage}
 \begin{minipage}{0.30\hsize}
  \begin{center}
    \includegraphics*[width=4cm,height=4cm]{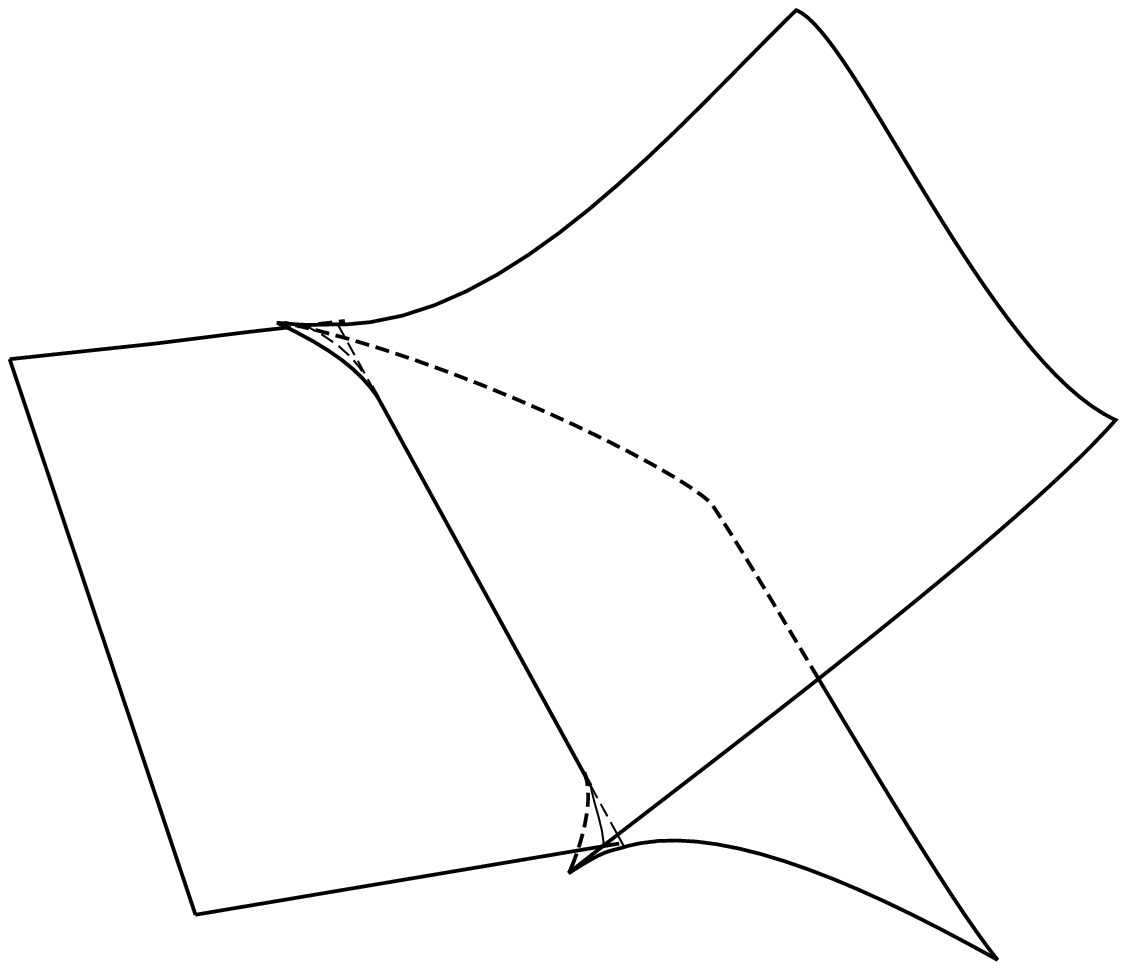}
\end{center}
 \end{minipage}
  \begin{minipage}{0.30\hsize}
  \begin{center}
   \includegraphics*[width=4cm,height=4cm]{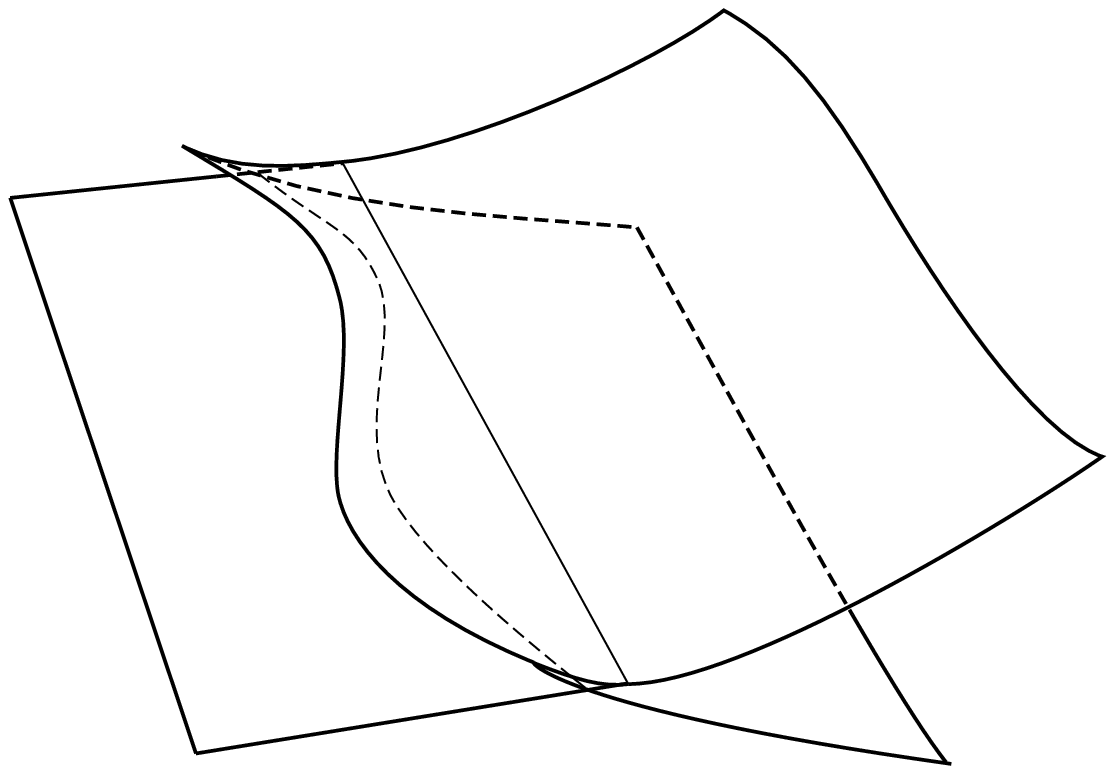}
\end{center}
 \end{minipage}
\end{center}
\caption{${}^2C_3^-$}
\end{figure}

\bibliographystyle{plain}

\begin{thebibliography}{1}

\bibitem{izumiya1}
S.~Izumiya and G.T.Kossioris.
\newblock Semi-local classification of geometric singularities for
  hamilton-jacobi equations.
\newblock {\em Journal of Differential equation}, 118:166--193, 1995.

\bibitem{retLeg}
T.~Tsukada.
\newblock Reticular legendrian singularities.
\newblock {\em The Asian J. of Math.}, 5(1):109--127, 2001.

\bibitem{tPKfunct}
T.~Tsukada.
\newblock {A generic classification of function germs with respect to the
  reticular $t$-${\cal P}$-${\cal K}$-equivalence}.
\newblock {\em Hokkaido Math. J.}, 38(1):177--203, 2009.

\bibitem{generic}
T.~Tsukada.
\newblock Genericity of caustics and wavefronts on an $r$-corner.
\newblock {\em Asian J. Math.}, 14(3):335--358, 2010.

\bibitem{bifsemi}
T.~Tsukada.
\newblock Bifurcations of wavefronts on r -corners: semi-local classification.
\newblock {\em Methods and Applications of Analysis}, 18(3):303--334, 2011.

\bibitem{spsing}
G.~Wassermann.
\newblock Stability of unfolding in space and time.
\newblock {\em Acta Mathematica}, 135(1):57--128, 1975.

\end{thebibliography}

\end{document}